\documentclass[11pt,reqno]{amsart}
\usepackage[french,english]{babel}
\usepackage[utf8]{inputenc}
\usepackage[T1]{fontenc}            
\usepackage{amsmath}
\usepackage{graphicx}
\usepackage{amssymb}
\usepackage{amsthm}
\usepackage{comment}
\usepackage{pgf, tikz}
\usepackage{amsfonts,mathrsfs}
\usepackage{thmtools}
\usepackage{geometry}
\usepackage{stmaryrd}
\usepackage{color}
\usepackage{hyperref}
\usepackage{bm}
\hypersetup{colorlinks = true , linkcolor=blue, citecolor=blue}
\numberwithin{equation}{section}
\theoremstyle{definition}

\declaretheorem[numberwithin=,title=Definition]{Def}
\declaretheorem[numberwithin=,title=Proposition]{Prop}
\declaretheorem[numberwithin=,title=Corollary]{Cor}
\declaretheorem[numberwithin=,title=Theorem]{Thm}
\declaretheorem[numberwithin=,title=Lemma]{Lemme}

\declaretheorem[numberwithin=,title=Remark]{Rq}
\declaretheorem[numberwithin=,title=Example]{Ex}

\newenvironment{dem}{\vspace{6pt} \paragraph{\textbf{Proof.} \hspace{3pt}}}{\hfill$\square$\newline }

\renewcommand{\P}{\mathbb{P}}
\newcommand{\PP}{\mathcal{P}_{\beta}(\R^d)}
\newcommand{\PPP}{\mathcal{P}}
\newcommand{\FF}{\mathcal{F}}
\newcommand{\E}{\mathbb{E}}

\newcommand{\KK}{\mathcal{K}}
\newcommand{\R}{\ensuremath{\mathbb{R}}\xspace}

\newcommand{\N}{\ensuremath{\mathbb{N}}\xspace}

\renewcommand{\S}{\ensuremath{\mathbb{S}}\xspace}
\newcommand{\eps}{\varepsilon}
\newcommand{\CC}{\mathcal{C}}

\newcommand{\NN}{\ensuremath{\mathcal{N}}\xspace}
\newcommand{\NNN}{\ensuremath{\widetilde{\mathcal{N}}}\xspace}

\newcommand{\1}{\mathbf{1}}
\newcommand{\MM}{\mathcal{M}}
\newcommand{\del}{\frac{\delta }{\delta  m}u}
\newcommand{\dell}{\frac{\delta }{\delta  m}\phi_u}
\newcommand{\dellN}{\frac{\delta }{\delta  m}\phi_{N,u}}
\newcommand{\de}{\frac{\delta }{\delta  m}}
\newcommand{\muu}{\overline{\mu}}
\newcommand{\x}{\bm{x}}

\begin{document}

	\title[Itô's formula for Poisson stochastic integrals and applications]{Itô's formula for the flow of measures of Poisson stochastic integrals and applications}

	\author{Thomas Cavallazzi}
	\address{Université de Rennes 1, CNRS, IRMAR - UMR 6625, F-35000 Rennes, France}
	\email{thomas.cavallazzi@univ-rennes1.fr}
	\date{November 27, 2022}

	\begin{abstract}
		We prove Itô's formula for the flow of measures associated with a jump process defined by a drift, an integral with respect to a Poisson random measure and with respect to the associated compensated Poisson random measure. We work in $\PPP_{\beta}(\R^d)$, the space of probability measures on $\R^d$ having a finite moment of order $\beta \in (0,2]$. As an application, we exhibit the backward Kolmogorov partial differential equation stated on $[0,T]\times \PPP_{\beta}(\R^d)$ associated with a McKean-Vlasov stochastic differential equation driven by a Poisson random measure. It describes the dynamics of the semigroup associated with the McKean-Vlasov stochastic differential equation, under regularity assumptions on it. Finally, we use the semigroup and the backward Kolmogorov equation to prove new quantitative weak propagation of chaos results for a mean-field system of interacting Ornstein-Uhlenbeck processes driven by i.i.d.\ $\alpha$-stable processes with $\alpha \in (1,2)$. 
	\end{abstract}

 \maketitle

\tableofcontents

\section{Introduction}

Let us fix $(\Omega,\FF,(\FF_t)_{t\geq0},\P)$ a filtered probability space satisfying the usual conditions and $\NN$ a Poisson random measure on $[0,T]\times \R^d$ for some finite horizon of time $T>0$, with intensity measure $dt \otimes \nu$, where $\nu$ is a Lévy measure, i.e.\ $\nu(\{0\})=0$ and $$ \int_{\R^d} \min(|z|^2,1) \, d\nu(z) < + \infty.$$ We consider the jump process $X=(X_t)_{t \in [0,T]}$ with values in $\R^d$ defined for all  $t \in [0,T]$ by \begin{equation}\label{process1}
X_t := X_0 + \int_0^t b_s \, ds + \int_0^t \int_{\{|z| <1 \}}H(s,z) \, \NNN(ds,dz) +  \int_0^t \int_{\{|z| \geq 1 \}}K(s,z) \, \NN(ds,dz), \end{equation} where $\NNN(ds,dz):=\NN(ds,dz)-ds\,d\nu(z)$ is the associated compensated Poisson random measure and $b$, $H$, $K$ are predictable processes. The distribution of $X_t$ is denoted by $\mu_t$. The assumptions made on $X_0$, $\nu$, $b$, $H$ and $K$ (see Section \ref{sectionIto}) ensure that for all $t \in [0,T],$ $\mu_t$ belongs to $\PPP_{\beta}(\R^d)$ the space of probability measures having a finite moment of order $\beta$, for some $\beta \in (0,2]$. \\

In this work, we are interested in Itô's formula for the flow of probability measures $(\mu_t)_{t \in [0,T]}.$ This formula describes the dynamics of the map $t \in [0,T] \mapsto u(\mu_t),$ for a function $ u :  \PPP_{\beta}(\R^d) \rightarrow \R$. Itô's formula for a flow of probability measures has been developed over the last decade to deal with mean-field games, which were initiated independently by Caines, Huang and Malhame in \cite{CainesLargepopulationstochastic2006} and by
Lasry and Lions in \cite{LasryMeanfieldgames2007}. An important object introduced by Lions in his lectures at Collège de France \cite{LionscourscollegedeFrance} is the  master equation related to a mean-field game, which is a Partial Differential Equation (PDE) stated on the space of probability measures.  We refer to \cite{LionscourscollegedeFrance,CardaliaguetnotesMFG2013,CarmonaProbabilisticTheoryMean2018,CarmonaProbabilisticTheoryMean2018a} for more details on mean-field games and associated master equations. As for the standard Itô formula, it allows to associate a PDE stated on the space of probability measures to a McKean-Vlasov Stochastic Differential Equation (SDE), which is the related master equation in this context. The well-posedness of these equations was studied for example in  \cite{ChassagneuxProbabilisticapproachclassical2015,BuckdahnMeanfieldstochasticdifferential2014,Cardaliaguetmasterequationconvergence2015,CrisanSmoothingpropertiesMcKeanVlasov2017,deraynal2021wellposedness}. Moreover, the problem of propagation of chaos for a mean-field interacting particle system towards the corresponding McKean-Vlasov SDE can also be tackled using the associated PDE on the space of measures. It turns out to be an efficient tool to obtain quantitative weak propagation of chaos estimates between the distribution of the solution to the McKean-Vlasov SDE $(\mu_t)_t$ and the empirical measure $(\muu^N_t)_t$ of the particle system. More precisely, a quantitative weak propagation result consists in finding an explicit rate of convergence with respect to $N$ for $\E|u(\muu^N_T) - u(\mu_T)|$ and $|\E u(\muu^N_T) - \E u(\mu_T)|$ and for $u$ belonging to a certain class of functions defined on the space of probability measures. This strategy was originally described in Chapter 5 of \cite{CarmonaProbabilisticTheoryMean2018}. It was adopted for example by Chaudru de Raynal and Frikha in \cite{deraynal2021wellposedness,frikha2021backward}, by Chassagneux, Szpruch and Tse in \cite{chassagneux2019weak} and by Delarue and Tse in \cite{delarue2021uniform}. Let us also mention that this PDE on the space of probability measures has recently been used in \cite{JourdainCentrallimittheorem2021} to prove a central limit theorem for interacting particle systems. Finally, Itô's formula for a flow of measures is also a key tool to tackle McKean-Vlasov control problems. Indeed, it induces a dynamic programming principle which describes the value function of the problem as presented in \cite{CarmonaForwardBackwardStochasticDifferential2013a} or in Chapter $6$ of \cite{CarmonaProbabilisticTheoryMean2018} (see also the references therein).\\

In the first part of this paper, we focus on Itô's formula for the flow of measures $(\mu_t)_t$ associated with \eqref{process1} and for a function $u : \PPP_{\beta}(\R^d)\rightarrow \R$ with $\beta \in (0,2]$. It states that for all $t \in [0,T]$, we have  \begin{align*}
	\notag&u(\mu_t) - u(\mu_0) \\ &=\int_0^t \E \left( \partial_v \del (\mu_s)(X_s) \cdot b_s \right) \, ds \\\notag &\quad+ \int_0^t \int_{|z|\geq 1} \E \left[ \del(\mu_s)(X_{s^-}+ K(s,z)) - \del(\mu_s)(X_{s^-})\right] \, d\nu(z) \, ds  \\ \notag&\quad +  \int_0^t \int_{|z|<1} \E \left[ \del(\mu_s)(X_{s^-}+ H(s,z)) - \del(\mu_s)(X_{s^-}) - \partial_v \del(\mu_s)(X_{s^-})\cdot H(s,z)\right] \, d\nu(z) \, ds,
\end{align*}
where $\del $ denotes the linear (functional) derivative of $u$ (see Definition \ref{deflinearderivative}). The precise assumptions  are given in Section \ref{sectionIto} (see Theorem \ref{ito}). In the Brownian setting, Itô's formula for a flow of measures has been established for example in \cite{BuckdahnMeanfieldstochasticdifferential2014} (see Theorem 6.1) or in Section $3$ of \cite{ChassagneuxProbabilisticapproachclassical2015} and Chapter 5 of \cite{CarmonaProbabilisticTheoryMean2018} (see Theorem 5.99) under less restrictive assumptions. Let us also mention \cite{dosReisItoWentzellLionsFormulaMeasure2022,CavallazziItoKrylovformula2021} for some recent extensions of Itô's formula for a flow of measures. Concerning jump processes, Itô's formula has recently been extended to flows of measures generated by càdlàg semi-martingales. It was established independently by Guo, Pham and Wei in \cite{GuoItoformulaflow2020}, who studied McKean-Vlasov control problems with jumps and by Talbi, Touzi and Zhang in \cite{talbi2021dynamic} who worked on mean-field optimal stopping problems. A common point of these two works is that the jump process considered has a finite moment of order $2$. However, this framework is not adapted when the Poisson random measure $\NN$ stems from an $\alpha$-stable Lévy process with $\alpha \in (0,2)$ since it only has finite moments of order $\beta < \alpha$. Itô's formula given in Theorem \ref{ito} can be applied for these processes since $\beta \in (0,2]$. Another difference between this work and \cite{GuoItoformulaflow2020} is that we do not assume $\partial_v \del$ to be uniformly bounded if $\beta > 1$. Moreover, the authors of \cite{talbi2021dynamic} assume that for all $\mu \in \PPP_2(\R^d)$, the function $\del (\mu)$ is of class $\CC^2$ on $\R^d$, which is not the case in this paper. It is replaced here by the $\gamma$-Hölder continuity of $\partial_v\del (\mu)$ uniformly in $\mu \in \PP$, where $\gamma$ is such that $x \mapsto |x|^{1+\gamma} \in L^1(\{|z| <1\},\nu)$. The strategy of the proof is the following. First, we localize $X=(X_t)_t$ using stopping times and we precise in Proposition \ref{propapprox} in which sense the sequence of localized processes $(X^n)_n$ approximates $X$. Then, we establish Itô's formula for the flow of measures associated with the stopped process $(X^n_t)_t$ thanks to a standard method of time discretization. We finally let $n$ tend to infinity essentially using the approximation results of Proposition \ref{propapprox}. Moreover, we extend Itô's formula for functions depending also on the time and space variables in Theorem \ref{ito2}. \\

Then, we consider a general Lévy-driven McKean-Vlasov SDE, which is assumed to be well-posed in the weak sense, of the form \begin{equation}\label{edsmckvintro}
\left\{  \begin{array}{lll}
&dX_t = b_t(X_t,\mu_t) \,dt + \sigma_t(X_{t^-},\mu_t)\, dZ_t, \quad t\in[0,T], \\ &\mu_t = \mathcal{L}(X_t),\\ &X_0 = \xi \in L^{\beta}(\Omega,\FF_0),
\end{array}\right.
\end{equation}
where $\beta \in (0,2]$ and $Z=(Z_t)_t$ is a Lévy process on $\R^d$ having the following Lévy-Itô decomposition \begin{equation*}
Z_t = \int_0^t \int_{\{|z| <1 \}} z \, \NNN(ds,dz) + \int_0^t \int_{\{|z| \geq 1 \}} z \, \NN(ds,dz).\end{equation*}  
Our aim is to describe the dynamics of the semigroup acting on functions defined on $\PP$ associated with the McKean-Vlasov SDE \eqref{edsmckvintro}. For a fixed function $u : \PPP_{\beta}(\R^d)\rightarrow \R$, the action of the semigroup on $u$ is given by $\phi_u$ defined by
\begin{equation}\label{phi_u_intro}
\phi_u : \left\{ \begin{array}{rrl}
&[0,T] \times \PPP_{\beta}(\R^d) & \rightarrow \R \\ &(t,\mu) &\mapsto u(\mathcal{L}(X_T^{T-t,\mu})),
\end{array}\right.
\end{equation}
where $\mathcal{L}(X_T^{T-t,\mu})$ denotes the distribution of the solution to \eqref{edsmckvintro} at time $T$ starting at time $T-t$ from a random variable $\xi$ with distribution $\mu.$ As in \cite{CarmonaProbabilisticTheoryMean2018}, we prove using Itô's formula for a flow of measures that under regularity assumptions on $\phi_u$, the function $(t,\mu) \in [0,T] \times \PPP_{\beta}(\R^d) \mapsto \phi_u(T-t,\mu)$ is a classical solution to a backward Kolmogorov PDE on $[0,T]\times\PPP_{\beta}(\R^d)$ (see Theorem \ref{ThmEDP}).\\

Finally, we prove new quantitative weak propagation of chaos results, in the sense previously defined, for a mean-field system of interacting Ornstein-Uhlenbeck processes defined as follows. We assume that the driving process $Z = (Z_t)_t$ is an $\alpha$-stable process on $\R^d$ with $\alpha \in (1,2)$ and with a non degenerate Lévy measure, i.e.\ satisfying \eqref{assumptionND}. More precisely, we introduce $(Z^n)_n$ an i.i.d.\ sequence of stable processes having the same distribution as $Z = (Z_t)_t,$ $A,A',B$ matrices of size $d \times d$ such that $B$ is invertible and $(X^n_0)_n$ an i.i.d.\ sequence of random variables with common distribution $\mu_0 \in \PPP_{\beta}(\R^d)$, for some $\beta \in [1,\alpha)$. The particle system is the unique solution to the following classical SDE on $(\R^{d})^N$, for $N \geq1$ \begin{equation}\label{edsparticlesOUintro}
\left\{  \begin{array}{lll}
&dX^{i,N}_t = AX^{i,N}_t\,dt + A' \displaystyle\frac{1}{N} \displaystyle\sum_{j=1}^N X^{j,N}_t \, dt +  B\,dZ^i_t, \quad t \in [0,T],  \quad \forall i \leq N,\\ &X^{i,N}_0 = X^i_0.
\end{array}\right.\end{equation}

We denote by $X=(X_t)_t$ the solution to the corresponding McKean-Vlasov SDE  \begin{equation}\label{MKVstableOUdefintro}
\left\{ \begin{array}{lll}
dX_t & = (AX_t + A' \E X_t)\,dt + B\,dZ_t,\quad t\in[0,T], \\X_0  &= \xi,
\end{array} \right.
\end{equation}
where the distribution of $\xi$ is $\mu_0$.\\

The propagation of chaos phenomenon was originally studied by McKean in \cite{McKeanPropchaos67} and later by Sznitman in \cite{SznitmanTopicspropagationchaos1991}. It states that for any $k \geq 1$, the limiting behaviour of $(X^{1,N}, \dots, X^{k,N})$, when $N$ tends to infinity, is expected to be described by $k$ independent copies of the McKean-Vlasov process \eqref{MKVstableOUdefintro}. It can be shown in the weak sense i.e.\ in distribution or in the strong sense i.e.\ at the level of paths. It has been of course widely studied in the Brownian case. Let us focus on some works dealing with mean-field systems with jumps. In \cite{GRAHAM199269}, following the approach of Sznitman \cite{SznitmanTopicspropagationchaos1991} in the Brownian case, Graham proves weak propagation of chaos under Lipschitz assumptions for a mean-field system driven by a Poisson random measure and its compensated measure. He works in the $L^1$ framework i.e.\ the Poisson random measure is associated with a Poisson process having a finite moment of order $1$. In the case of a McKean-Vlasov SDE driven by a general Lévy process having a finite moment of order $2$, we refer to Jourdain, Méléard and Woyczynski \cite{JourdainNonlinearSDEsdriven2007}. The authors prove quantitative rates of convergence for the strong propagation of chaos in $L^2$ under standard Lipschitz assumptions. Moreover, it has also been established in \cite{NeelimaWellposednesstamedEuler2020} by Neelima \& al. under relaxed assumptions. We also mention \cite{Mischlernewapproachquantitative2015} where Mischler, Mouhot and Wennberg exhibit conditions leading to quantitative rates of convergence for the propagation of chaos. As an application, they study an inelastic Boltzmann collision jump process. Finally, let us also refer to \cite{FrikhaWellposednessapproximationonedimensional2020}, where Frikha and Li study a one-dimensional McKean-Vlasov SDE driven by a compensated Poisson random measure with positive jumps. They prove quantitative rates of convergence for the strong propagation of chaos in $L^1$ under one-sided Lipschitz assumptions. \\

 Concerning our mean-field system of interacting Ornstein-Uhlenbeck processes \eqref{edsparticlesOUintro}, we are interested in quantitative weak propagation of chaos. Let us denote by $\muu^N_t$ the empirical measure of the particle system $\eqref{edsparticlesOUintro}$ and by $\mu_t$ the distribution of $X_t$. We prove in Theorem \ref{thmpropchaosOU} that there exists a constant $C>0$ such that for any $u : \PPP_{\beta}(\R^d) \rightarrow \R$ in a certain class of regular functions described in Theorem \ref{thmpropchaosOU}, one has \begin{equation}\label{eq1thmOU_intro}
	\E \left\vert u(\muu^N_T) - u(\mu_T) \right\vert \leq C\, \E W_1(\muu^N_0,\mu_0) +C\ln(N)^{\frac{1}{\alpha}}N^{\frac{1}{\alpha}-1},\end{equation} where $W_1$ is the Wasserstein distance of order $1$. We also show that if $\mu_0 \in \PPP_2(\R^d)$, one has \begin{equation}\label{eq2thmOU_intro}|\E (u(\muu^N_T) - u(\mu_T))|\leq CN^{1 - \alpha }.
 \end{equation}  Of course, the rate of convergence is better at the level of semigroup i.e.\ in \eqref{eq2thmOU_intro} than in \eqref{eq1thmOU_intro}. Notice that if we take $\alpha = 2$, which formally corresponds to the Brownian case, we recover the same rates of convergence, up to the factor $\ln(N)$ in \eqref{eq1thmOU_intro}, as in Theorem $3.6$ of \cite{frikha2021backward}, even though the drift is unbounded here. In dimension $d=1$, we recover with \eqref{eq1thmOU_intro} the same rate of convergence obtained in \cite{FrikhaWellposednessapproximationonedimensional2020}, for the strong propagation of chaos in $L^1$, since $\E W_1(\muu^N_0,\mu_0) \leq CN^{\frac{1}{\beta}-1}$ by \cite{fournier:hal-00915365}. Notice that if $\mu_0$ has enough moments, for example if it belongs to $\PPP_2(\R)$, then using \cite{fournier:hal-00915365}, we have a better rate of propagation of chaos since \eqref{eq1thmOU_intro} becomes \[	\E \left\vert u(\muu^N_T) - u(\mu_T) \right\vert \leq C\, \ln(N)^{\frac{1}{\alpha}}N^{\frac{1}{\alpha}-1}.\] \\

The method that we use relies on regularity properties and estimates on the solution to the backward Kolmogorov PDE $\phi_u$ defined in \eqref{phi_u_intro} (see Proposition \ref{ThmPCremovjumps}). This strategy was originally described in Chapter 5 of \cite{CarmonaProbabilisticTheoryMean2018} (pages $506-508$), inspired by \cite{Cardaliaguetmasterequationconvergence2015} and \cite{Mischlernewapproachquantitative2015}, and was employed for example in \cite{chassagneux2019weak,delarue2021uniform,frikha2021backward}, as mentioned above. Let us describe the main ideas. We begin by computing the time derivative of the map $ t \in [0,T) \mapsto \phi_u(T-t,\muu^N_t)$ by applying the standard Itô's formula for the empirical projection $(t,x_1,\dots,x_N)~\in ~[0,T] \times (\R^d)^N \mapsto \phi_u\left(T-t, \frac1N \sum_{k=1}^N \delta_{x_k}\right)$ and for the particle system. Noting that $t \in[0,T] \mapsto \phi_u(T-t,\mu_t)$ is constant, we naturally expect that the time derivative previously computed tends to $0$ as $N$ converges to infinity. This convergence has to be shown with an explicit rate of convergence using the PDE satisfied by $\phi_u$ and some estimates on $\phi_u$. Finally, we express $u(\muu^N_T) - u(\mu_T) = \phi_u(0,\muu^N_T) - \phi_u(0,\mu_T)$ as the sum of $ \phi_u(T,\muu^N_0) - \phi_u(T,\mu_0)$ plus a remainder term related to the time derivative previously estimated. Since the initial data are i.i.d., the first term is controlled by standard estimates, for example those in \cite{fournier:hal-00915365}. \\

An important difference in the jump case in comparison to the Brownian case is that the Kolmogorov PDE satisfied by $\phi_u$ does not directly appear when we apply Itô's formula for the empirical projection of $\phi_u$ and for the particle system. In order to use the backward Kolmogorov PDE, we thus have to control the corresponding error term. To do this and because of the unboundedness of the drift, we need to remove the big jumps of the driving processes in a first step. The key point is that we consider the solutions to \eqref{edsparticlesOUintro} and \eqref{MKVstableOUdefintro} driven by truncated noises for which we remove the jumps bigger than the number of particles in interaction $N$. We need to perform this truncation procedure to control the error term mentioned above. This is essentially because the moment of order $2$ of the Lévy measure $\nu$ appears in our computations due to the unboundedness of the drift. We refer to Remark \ref{rqbigjumps} for more details. The presence of the factor $\ln(N)$ in \eqref{eq1thmOU_intro} comes from this procedure.\\
 
  Let us eventually emphasize that this mean-field system of interacting Ornstein-Uhlenbeck  processes is a first application. We aim at showing that this method, already used in the Brownian setting, and leading to quantitative weak propagation of chaos results, can be generalized in the context of jump processes. 
  Our next goal is to follow the same strategy to establish quantitative weak propagation of chaos results for a general class of Lévy-driven McKean-Vlasov SDEs, under mild assumptions on the coefficients, as done in \cite{frikha2021backward}.\\

The paper is organized as follows. In Section \ref{sectionIto}, we state and prove Itô's formula in Theorem \ref{ito} and Theorem \ref{ito2} for the flow of measure of Poisson stochastic integrals. In Section \ref{sectionPDE}, and more precisely in Theorem \ref{ThmEDP}, we use our Itô's formula to derive the backward Kolmogorov PDE on the space of measures describing the dynamics of the semigroup associated with a general Lévy-driven McKean-Vlasov SDE. Then, in Section \ref{sectionOU}, we study the mean-field system of interacting Ornstein-Uhlenbeck processes and we prove in Theorem \ref{thmpropchaosOU} quantitative weak propagation of chaos. Appendix~ \ref{appendixdensityestimates} is devoted to the proof of moment estimates, which are uniform with respect to the truncation of the big jumps, on the density and its derivatives of a stable Ornstein-Uhlenbeck process. Finally, in Appendix~~\ref{proofpropositionreg}, we prove the regularity properties and the controls required in Section \ref{sectionOU} on the solution to the backward Kolmogorov PDE associated with \eqref{MKVstableOUdefintro}.\\

 Let us finally introduce some notations used several times in the article.\\

\textbf{Notations}

\begin{enumerate}

	\item[-] $\PPP_{\beta}(\R^d)$ denotes the set of probability measures $\mu$ on $\R^d$ such that $\int_{\R^d} |x|^\beta \, d\mu(x) < + \infty,$ for $\beta >0.$ It is equipped with the Wasserstein metric of order $\beta$ denoted by $W_\beta,$ which makes it complete. Denoting by $\Pi(\mu,\nu)$ the set of couplings between two probability measures $\mu, \nu \in \PPP_{\beta}(\R^d)$, the metric $W_\beta$ is defined by $$W_{\beta}(\mu,\nu) = \inf_{\pi \in\Pi(\mu,\nu) } \left(\int_{\R^d\times \R^d} |x-y|^\beta \, d\pi(x,y)\right)^{\frac{1}{\beta}} \quad \text{if $\beta \geq 1$,}$$ and by  $$W_{\beta}(\mu,\nu) = \inf_{\pi \in \Pi(\mu,\nu)} \left(\int_{\R^d\times \R^d} |x-y|^\beta \, d\pi(x,y)\right)\quad \text{if $\beta < 1$.}$$ 
	\item[-]  $\muu^N_{\bm{x}}:= \frac1N \sum_{k=1}^N \delta_{x_k}$ denotes the empirical measure, for $\bm{x} = (x_1,\dots,x_N) \in (\R^d)^N$.
	\item[-]  $\bm{\tilde{z}_k} := (0, \dots,z,\dots,0) \in (\R^d)^N$ for $z \in \R^d,$ where $z$ appears in the $k$-th position.

\item[-] $B_r$ is the open ball centered at $0$ and of radius $r$ in $\R^d$ for the euclidean norm.

\item[-] $B_r^c$ denotes the complementary of $B_r$.

	\item[-] $p'$ is the conjugate exponent of $p \in [1,\infty]$ defined by $ \frac1p + \frac{1}{p'} =1.$
	\item[-] $a\wedge b$ denotes the minimum between $a$ and $b$.
	\item[-] $a \vee b $ denotes the maximum between $a$ and $b$.
	
		\item[-] $C$ is a generic constant that may depend only on the fixed parameters of the problem and which may change from line to line.
	
\end{enumerate}

\section{Itô's formula for the flow of measures of Poisson stochastic integrals}\label{sectionIto}

\subsection{Assumptions and Itô's formula for a flow of measure}

 Let $\nu$ be a Lévy measure on $\R^d$ i.e.\ $\nu(\{0\})=0$ and $$ \int_{\R^d} |z|^2\wedge 1 \, d\nu(z) < + \infty.$$ We introduce a $(\FF_t)_{t\in[0,T]}$-Poisson random measure $\NN$ on $[0,T]\times \R^d\backslash\{0\}$ with intensity measure $dt \otimes \nu,$ and we denote by $\NNN(ds,dz):=\NN(ds,dz) - ds\,d\nu(z)$ the compensated Poisson random measure associated with $N.$ Let us define the jump process $X=(X_t)_{t \in [0,T]}$ by \begin{equation}\label{process}
 	\forall t \in [0,T], \, X_t := X_0 + \int_0^t b_s \, ds + \int_0^t \int_{B_1}H(s,z) \, \NNN(ds,dz) +  \int_0^t \int_{B_1^c}K(s,z) \, \NN(ds,dz), \end{equation} where $b: [0,T] \times \Omega \rightarrow \R^d,$ $H : [0,T]\times B_1 \times \Omega \rightarrow \R^d,$ and  $K : [0,T]\times B_1^c \times \Omega \rightarrow \R^d$ are predictable processes. The distribution of $X_t$ is denoted by $\mu_t.$\\

Let us fix our assumptions on the jump process $(X_t)_t$ in order to establish Itô's formula for the flow $(\mu_t)_t.$ We assume that there exists two constants $\beta \in (0,2]$ and $\gamma \in [0,1]$ such that $\beta \leq 1 + \gamma$ and satisfying the following properties.

\begin{enumerate}
	
	\item[] \textbf{(M)} The random variable $X_0$ belongs to $L^{\beta}(\Omega;\R^d)$ and we have  \begin{equation}\label{assumptM}
		\E\int_0^T |b_s|^{\beta\vee 1} \, ds < + \infty.	\end{equation}
	\item[] \textbf{(J1)} There exists a predictable process $(\widetilde{H}_s)_{s \in [0,T]}$ which is assumed to be almost surely locally bounded and to satisfy  \begin{equation}\label{assumptJ1}  \text{a.s.} \, \forall s \in [0,T],\, \forall z \in B_1,\, |H(s,z)| \leq |\widetilde{H}_s||z| \quad \text{and}\quad \E \int_0^T \int_{B_1} (|\widetilde{H}_s||z|)^{1+\gamma} \, d\nu(z) \, ds < + \infty.	\end{equation}
		\item[] \textbf{(J2)} We have  \begin{equation}\label{assumptJ2} \E \int_0^T \int_{B_1^c} |K(s,z)|^{\beta} \, d\nu(z) \, ds < + \infty.	\end{equation}

\end{enumerate}

\begin{Def}[Linear derivative] \label{deflinearderivative}A function $u : \PPP_{\beta}(\R^d) \rightarrow \R$ is said to have a linear derivative if there exists a function $\del \in \CC^0(\PPP_{\beta}(\R^d)\times \R^d;\R)$ satisfying the following properties.
	\begin{enumerate}
		\item For all compact subset $\KK \subset \PPP_{\beta}(\R^d),$ there exists a constant $C_{\KK}>0$ such that $$ \forall \mu \in \KK,\, \forall v \in \R^d, \, \left\vert \del (\mu)(v) \right\vert \leq C_{\KK} (1+|v|^{\beta}).$$
		\item For all $\mu,\nu \in \PPP_{\beta}(\R^d),$ we have $$ u(\mu)-u(\nu) = \int_0^1 \int_{\R^d} \del (t \mu + (1-t)\nu)(v)\, d(\mu-\nu)(v) \, dt.$$
	\end{enumerate}
\end{Def}

We now state Itô's formula for the flow of measures $(\mu_t)_{t \in [0,T]}$.

\begin{Thm}[Itô's formula]\label{ito} Suppose that Assumptions $(\textbf{M}), \textbf{(J1)},$ and $\textbf{(J2)}$ are satisfied. Let $u: \PP \rightarrow \R$ be a function having a linear derivative $\del$ which satisfies the following properties. \begin{enumerate}
		\item For any $\mu \in \PP,$ the function $\del(\mu) \in \CC^1(\R^d;\R)$ and $\partial_v \del$ is continuous on $\PP \times \R^d.$
		
			\item If $\gamma>0,$ for any compact $\KK \subset \PP,$ there exists $C_{\KK}>0$ such that $$  \forall \mu \in \KK, \, \forall x,y \in \R^d,\, \displaystyle\left\vert \partial_v\del(\mu)(x) - \partial_v\del(\mu)(y)  \right\vert \leq C_{\KK} |x-y|^{\gamma}.$$

		\item For any compact $\KK \subset \PP,$ we have  $$\left\{ \begin{array}{lll}  &\displaystyle\sup_{\mu \in \KK} \int_{\R^d} \left\vert \partial_v \del (\mu)(v) \right\vert^{\beta'} \, d\mu(v) < + \infty\quad \text{if ${\beta >1},$} \\ &  \displaystyle\sup_{v \in \R^d}\sup_{\mu\in \KK} \left\vert \partial_v \del (\mu)(v) \right\vert < + \infty \quad\text{if $\beta \leq1$}. \end{array} \right.$$
		\end{enumerate}
		
		Then, we have for all $t \in [0,T]$ \begin{align}\label{itoformula}
	\notag&u(\mu_t) - u(\mu_0) \\ &=\int_0^t \E \left( \partial_v \del (\mu_s)(X_s) \cdot b_s \right) \, ds \\\notag &\quad+ \int_0^t \int_{B_1^c} \E \left[ \del(\mu_s)(X_{s^-}+ K(s,z)) - \del(\mu_s)(X_{s^-})\right] \, d\nu(z) \, ds  \\ \notag&\quad +  \int_0^t \int_{B_1} \E \left[ \del(\mu_s)(X_{s^-}+ H(s,z)) - \del(\mu_s)(X_{s^-}) - \partial_v \del(\mu_s)(X_{s^-})\cdot H(s,z)\right] \, d\nu(z) \, ds. 
		\end{align}

\end{Thm}

\begin{Rq}\label{rqgrowth}
	 Assumption $(3)$ is implied by the following stronger assumption: for all compact $\KK \subset \PP $, there exists $ C_{\KK}>0$ such that  $$\forall v \in \R^d, \, \sup_{\mu \in \KK} \left\vert \partial_v \del(\mu)(v)\right\vert \leq C_{\KK}(1+|v|^{\beta-1}\1_{\beta>1}).$$ When $\beta >1,$ it follows from Hölder's inequality.  
\end{Rq}

We specify in the three following corollaries the particular case where the Poisson random measure is associated with an $\alpha$-stable Lévy process $Z=(Z_t)_t$ with $\alpha \in (0,2)$. Consider $\sigma : [0,T]\times \Omega \rightarrow \R^{d \times d}$ a predictable process such that $$ \E \int_0^T |\sigma_s|^{1+\gamma} \, ds < + \infty$$ for some $\gamma \in [0,1]$ which will be specified. Fix also $X_0$ and $b$ satisfying Assumption \textbf{(M)} for some $\beta \leq 1 + \gamma$ to be specified. In this example, the process $X$ considered is defined for all $t \in [0,T] $ by $$ X_t = X_0 + \int_0^t b_s \,ds + \int_0^t \sigma_s \, dZ_s.$$

\begin{Cor}[Sub-critical case]\label{coralpha>1}Assume that $\alpha \in (1,2),$ $\beta \in (1,\alpha)$ and $\gamma \in (\alpha -1,1].$ Then, the process $X=(X_t)_t$ satisfies Assumptions \textbf{(M)}, \textbf{(J1)} and \textbf{(J2)} with $ K(s,z)=H(s,z)=\sigma_s z.$ Let $u: \PP \rightarrow \R$ be a function having a linear derivative $\del$ which satisfies the following properties. \begin{enumerate}
		\item For any $\mu \in \PP,$ the function $\del(\mu) \in \CC^1(\R^d;\R)$ and $\partial_v \del$ is continuous on $\PP \times \R^d.$
		
		\item For any compact $\KK \subset \PP,$ there exists $C_{\KK}>0$ such that $$  \forall \mu \in \KK, \, \forall x,y \in \R^d,\, \displaystyle\left\vert \partial_v\del(\mu)(x) - \partial_v\del(\mu)(y)  \right\vert \leq C_{\KK} |x-y|^{\gamma}.$$

		\item For any compact $\KK \subset \PP,$ we have  $$\displaystyle\sup_{\mu \in \KK} \int_{\R^d} \left\vert \partial_v \del (\mu)(v) \right\vert^{\beta'} \, d\mu(v) < + \infty.$$
	\end{enumerate} Then, Itô's formula \eqref{itoformula} holds true for $X$ and $u.$
\end{Cor}

\begin{Cor}[Super-critical case]\label{coralpha<1}Assume that $\alpha \in (0,1),$ $\beta \in (0,\alpha)$ and $\gamma =0.$ Then, the process $X=(X_t)_t$ satisfies Assumptions \textbf{(M)}, \textbf{(J1)} and \textbf{(J2)} with $ K(s,z)=H(s,z)=\sigma_s z.$ Let $u: \PP \rightarrow \R$ be a function having a linear derivative $\del$ which satisfies the following properties. \begin{enumerate}
		\item For any $\mu \in \PP,$ the function $\del(\mu) \in \CC^1(\R^d;\R)$ and $\partial_v \del$ is continuous on $\PP \times \R^d.$

		\item  For any compact $\KK \subset \PP,$ we have  $$\forall v \in \R^d,\, \displaystyle\sup_{v \in \R^d}\sup_{\mu\in \KK} \left\vert \partial_v \del (\mu)(v) \right\vert < + \infty.$$
	\end{enumerate} Then, Itô's formula \eqref{itoformula} holds true for $X$ and $u.$
\end{Cor}

\begin{Cor}[Critical case]\label{coralpha=1}Assume that $\alpha=1$ $\beta \in (0,1)$ and $\gamma \in (0,1].$  Then, the process $X=(X_t)_t$ satisfies Assumptions \textbf{(M)}, \textbf{(J1)} and \textbf{(J2)} with $ K(s,z)=H(s,z)=\sigma_s z.$ Let $u: \PP \rightarrow \R$ be a function having a linear derivative $\del$ which satisfies the following properties. \begin{enumerate}
		\item For any $\mu \in \PP,$ the function $\del(\mu) \in \CC^1(\R^d;\R)$ and $\partial_v \del$ is continuous on $\PP \times \R^d.$
		
		\item For any compact $\KK \subset \PP,$ there exists $C_{\KK}>0$ such that $$  \forall \mu \in \KK, \, \forall x,y \in \R^d,\, \displaystyle\left\vert \partial_v\del(\mu)(x) - \partial_v\del(\mu)(y)  \right\vert \leq C_{\KK} |x-y|^{\gamma}.$$

		\item  For any compact $\KK \subset \PP,$ we have  $$\forall v \in \R^d,\, \displaystyle\sup_{v \in \R^d}\sup_{\mu\in \KK} \left\vert \partial_v \del (\mu)(v) \right\vert < + \infty.$$
	\end{enumerate} Then, Itô's formula \eqref{itoformula} holds true for $X$ and $u.$
\end{Cor}

Let us give an example of a function $u$ which satisfies the assumptions of the three preceding corollaries. We take $d=1$ to simplify the computations.

\begin{Ex} Let $\beta \in (0,2).$ For $\eps >0,$ consider a function $\chi_{\eps}$ of class $\CC^2$ on $\R$ such that \begin{enumerate}
	\item[-] $\forall x \in \R, \, \chi_{\eps}(x) \in [0,1]$,
	\item[-]$\chi_{\eps}$ is equal to $1$ on $[-2\eps,2\eps]^c$ and to $0$ on $[-\eps,\eps].$
\end{enumerate}
	
	Define the function $u$ by $$ \forall \mu \in \mathcal{P}(\R), \, u(\mu):= \int_{\R} |x|^{\beta} \chi_{\eps}(x) \, d\mu(x).$$
	
Then, $u$ satisfies the assumptions of Corollaries \ref{coralpha>1}, \ref{coralpha<1} and \ref{coralpha=1} with $\gamma =1.$
\end{Ex}

\begin{proof}
	An easy computation shows that $u$ has a linear derivative given, for all $\mu \in \PPP_{\beta}(\R^d)$ and $ v \in \R^d$, by $$\del(\mu)(v) =  |v|^{\beta} \chi_{\eps}(v),$$ which is clearly of class $\CC^2.$ Moreover, one has for all $\mu \in \PPP_{\beta}(\R^d)$ and $ v \in \R^d$ $$ \forall \mu \in \mathcal{P}_{\beta}(\R),\, \forall v \in \R,\, \partial_v \del(\mu)(v) =\beta \text{sgn}(v)|v|^{\beta -1}\chi_{\eps}(v) + |v|^{\beta}  \chi_{\eps}' (v).$$
	Thus, the condition in Remark \ref{rqgrowth} is satisfied, i.e.\ for all $\mu \in \mathcal{P}_{\beta}(\R)$ and $v \in \R$ $$ \left\vert \partial_v \del (\mu)(v)\right\vert \leq C (1 + |v|^{\beta-1} \1_{\beta>1}).$$ Moreover, we have for all $\mu \in \PPP_{\beta}(\R^d)$ and $ v \in \R^d$ $$  \partial^2_v \del(\mu)(v) =\beta(\beta-1) \text{sgn}(v)|v|^{\beta -2}\chi_{\eps}(v) + |v|^{\beta} \chi_{\eps}'' (v) + 2\beta \text{sgn}(v)|v|^{\beta -1}\chi_{\eps}'(v).$$
	The right-hand side term is clearly bounded on $\R,$ and thus Assumption $(2)$ in Theorem \ref{ito} is satisfied with $\gamma =1$ by the mean value theorem.
\end{proof}

We can easily deduce the following extension of Itô's formula in the case where the function $u$ depends also on the time and space variables. Let us introduce $\NN^2$ another Poisson random measure on $[0,T] \times \R^d$ with Lévy measure $\nu^2$ such that $\int_{B_1}|z|^{1+\gamma'} \, d\nu^2(z) < + \infty$, for some $\gamma' \in (0 ,1]$. Then we define \[Y_t := Y_0 + \int_0^t b^2_s \, ds + \int_{0}^t \int_{B_1} H^2(s,z) \, \NNN^2(ds,dz) + \int_0^t \int_{B_1^c} K^2(s,z) \, \NN^2(ds,dz), \]
where $Y_0$ is $\FF_0$-measurable, $b^2,H^2$ and $K^2$ are predictable processes with $$\int_0^T |b^2_s| \, ds +  \int_0^T \int_{B_1} |H^2(s,z)|^2 \, d\nu^2(z) \, ds < + \infty \quad \text{a.s.}$$
 
\begin{Thm}[Extension of Itô's formula]\label{ito2}
		Instead of \eqref{assumptM}, \eqref{assumptJ1} and \eqref{assumptJ2} in Assumptions \textbf{(M)}, \textbf{(J1)} and \textbf{(J2)}, we assume that \begin{equation}\label{assumption_extension_ito}
		\E|X_0|^\beta +\E \sup_{t\leq T} |b_s|^{\beta\vee 1 } + \int_{B_1^c} \E\sup_{t\leq T} |K(t,z)|^\beta \, d\nu(z)  + \int_{B_1} \E\sup_{t \leq T} (|\widetilde{H}_t||z|)^{1 + \gamma} \, d\nu(z) < + \infty,
	\end{equation}
 with $\beta \leq 1 + \gamma$. Moreover, we assume that for all $t \in [0,T]$ and $z \in \R^d$, $K(\cdot,z)$ and $H(\cdot,z)$ is almost surely continuous at $t$.\\

Let $u:[0,T]\times \R^d \times \PP \rightarrow \R$ be a continuous function satisfying the following properties. \begin{enumerate}
		\item  For any $\mu \in \PP,$ $x\in \R^d$ $u(\cdot,x,\mu)$ is continuously differentiable and $\partial_t u $ is continuous on $[0,T]\times \R^d\times \PP$. Moreover for any $t \in [0,T],$ $\mu \in \PP$, $u(t,\cdot,\mu)$ is of class $\CC^1$ on $\R^d$ with $\partial_x u$ continuous on $[0,T]\times \R^d \times \PP$ and such that $\partial_x u(t,\cdot,\mu)$ is $\gamma'$-Hölder continuous uniformly with respect to $t$ and $\mu$.
		\item For any $t \in [0,T],$ $x \in \R^d$, the function $u(t,x,\cdot)$ admits a linear derivative $\del(t,x,\cdot)(\cdot)$ such that $\del$ is continuous on $[0,T]\times \PP \times \R^d$ and for all compact $\KK \subset \R^d \times\PP,$ there exists $C_{\KK} >0$ such that  $$ \forall t \in [0,T], \, \forall (x,\mu) \in \KK,\, \forall v\in \R^d, \, \left\vert \del (t,x,\mu)(v) \right\vert \leq C_{\KK}(1+|v|^\beta).$$ Moreover, we assume that for any $ t \in [0,T],$ $x\in \R^d$, $\mu \in \PP,$ $\del(t,x,\mu)$ is differentiable with respect to the space variable and such that $\partial_v \del$ is continuous on $[0,T]\times \R^d \times\PP \times \R^d.$
		
		\item If $\gamma>0,$ for any compact $\KK \subset \R^d \times\PP,$ there exists $C_{\KK}>0$ such that $$  \forall t\in [0,T],\, \forall (x,\mu) \in \KK, \, \forall v,v' \in \R^d,\, \displaystyle\left\vert \partial_v\del(t,x,\mu)(v) - \partial_v\del(t,x,\mu)(v')  \right\vert \leq C_{\KK} |v-v'|^{\gamma}.$$

		\item For any compact $\KK \subset \R^d \times \PP,$ we have  $$\left\{ \begin{array}{lll}  &\displaystyle \sup_{t \in[0,T]}\sup_{(x,\mu) \in \KK} \int_{\R^d} \left\vert \partial_v \del (t,x,\mu)(v) \right\vert^{\beta'} \, d\mu(v) < + \infty \quad\text{if ${\beta >1},$} \\&  \displaystyle \sup_{t \in [0,T]}\sup_{v \in \R^d}\sup_{(x,\mu)\in \KK} \left\vert \partial_v \del (t,x,\mu)(v) \right\vert < + \infty \quad\text{if $\beta \leq1.$}  \end{array} \right.$$
	\end{enumerate}

Then, the function $(t,x) \in [0,T]\times \R^d \mapsto u(t,x,\mu_t)$ is of class $\CC^{1}$, with $\partial_x u(t,\cdot,\mu_t)$ $\gamma'$-Hölder continuous uniformly with respect to $t$. Moreover, we have almost surely for all $t \in [0,T]$ \begin{align}\label{itoformula2}
	\notag&u(t,Y_t,\mu_t) - u(0,Y_0,\mu_0) \\ \notag&=  \int_0^t\partial_t u(s,Y_s,\mu_s) \, ds  + \int_0^t \overline{\E} \left( \partial_v \del (s,Y_s,\mu_s)(\overline{X}_s) \cdot \overline{b}_s \right)\, ds  \\\notag &\quad+  \int_0^t\int_{B_1^c} \overline{\E} \left[ \del(s,Y_s,\mu_s)(\overline{X}_{s^-}+ \overline{K}(s,z)) - \del(s,Y_s,\mu_s)(\overline{X}_{s^-})\right] \, d\nu(z) \,ds \\ \notag&\quad + \int_0^t \int_{B_1} \overline{\E} \left[ \del(s,Y_s,\mu_s)(\overline{X}_{s^-}+ \overline{H}(s,z)) - \del(s,Y_s,\mu_s)(\overline{X}_{s^-}) \right.\\ &\left. \hspace{7cm}- \partial_v \del(s,Y_s,\mu_s)(\overline{X}_{s^-})\cdot \overline{H}(s,z)\right] \, d\nu(z) \,ds \\ &\notag \quad+\int_0^t \partial_x u(s,Y_s,\mu_s)\cdot b^2_s\,ds + \int_0^t \int_{B_1^c} u(s,Y_{s^-} + K^2(s,z),\mu_s) - u(s,Y_{s^-},\mu_s) \, \NN^2(ds,dz) \\ &\notag \quad+ \int_0^t \int_{B_1} u(s,Y_{s^-} + H^2(s,z),\mu_s) - u(s,Y_{s^-},\mu_s) \, \NNN^2(ds,dz) \\ &\notag \quad+ \int_0^t \int_{B_1} u(s,Y_{s^-} + H^2(s,z),\mu_s) - u(s,Y_{s^-},\mu_s) - \partial_x u(s,Y_{s^-},\mu_s)\cdot H^2(s,z) \, d\nu^2(z) \, ds,
	\end{align}
where $(\overline{\Omega},\overline{\FF},\overline{\P})$ is an independent copy of $(\Omega,\FF,\P)$ and $(\overline{b},\overline{H},\overline{K},\overline{X})$ is a copy of $(b,H,K,X)$.
	
\end{Thm}

Before proving Itô's formulas of Theorem \ref{ito} and Theorem \ref{ito2}, we gather in the next section some useful properties of the process $X$ and on its flow $(\mu_t)_t.$ We also introduce the localization of $X$ that is at the heart of our proof. 

\subsection{Preliminary study of the process and localization}

We introduce the driving Lévy process $Y$ defined for $t \in [0,T]$ by \begin{equation}\label{defLevy}
	 Y_t :=  \int_0^t \int_{B_1} z \, \NNN(ds,dz) + \int_0^t \int_{B_1^c} z \,\NN(ds,dz).\end{equation} The jumping times of $Y$ of size greater than $1$ are denoted by $(\tilde{T}_n)_n$ and the associated jumps $(Z_n)_n$ are defined for all $n \geq 1$ by $$Z_n := \Delta Y_{\tilde{T}_n} := Y_{\tilde{T}_n} - Y_{\tilde{T}_n^-}.$$ It is well-known that $(\tilde{T}_{n+1}-\tilde{T}_n)_n$ are i.i.d.\ with common distribution exponential of parameter $\nu(B_1^c)$ and that $(Z_n)_n$ is independent of $(\tilde{T}_n)_n$ with common distribution $\frac{\nu_{|B_1^c}}{\nu(B_1^c)}.$ We can write for all $t \in [0,T]$  $$  \int_0^t \int_{B_1^c} z \,\NN(ds,dz)= \sum_{n=1}^{\infty} Z_n \1_{\tilde{T}_n \leq t} = \sum_{n=1}^{N_t} Z_n,$$ where $N_t = \displaystyle\sum_{k=1}^{\infty} \1_{\tilde{T}_k\leq t}$ is the associated Poisson process. Let us begin with the continuity of the flow $(\mu_t)_t$

\begin{Prop}\label{momentXprop}
	There exists a constant $C>0$ such that for all stopping time $\tau$ and for all $0 \leq s\leq t\leq T $, one has \begin{align}\label{controlcontinuityflow} \notag&\E |X_{t \wedge \tau} - X_{s \wedge \tau}|^{\beta}\\ & \leq C\left[ \left( \E \int_s^t|b_u|^{\beta \vee 1} \ du \right)^{\frac{\beta}{\beta \vee 1}}+  \left( \E\int_s^t \int_{B_1} |H(r,z)|^{1 + \gamma} \, d\nu(z) \, dr\right)^{\frac{\beta}{1 + \gamma}} + \E\int_s^t \int_{B_1^c} |K(r,z)|^{\beta} \, d\nu(z) \, dr\right].
	\end{align}Thus, the map $t \in [0,T] \mapsto \mu_t \in \PP$ is continuous. Moreover, we have the following uniform in time moment estimate \begin{equation}\label{momentX}
	 \E \sup_{t\leq T} |X_t|^{\beta} < + \infty.\end{equation}
\end{Prop}

\begin{proof} For the drift term, we work in $L^{\beta \vee 1}(\Omega).$ In fact, Jensen's inequality ensures that $$ \E \left\vert \int_{s \wedge \tau}^{t \wedge \tau} b_u \, du \right\vert^{\beta\vee 1} \leq C\, \E \int_{s \wedge \tau}^{t \wedge \tau}|b_u|^{\beta\vee 1} \, du.$$ Then, a discussion on the relative position of $\tau$ with respect to $s$ and $t$ yields $$ \E \left\vert \int_{s \wedge \tau}^{t \wedge \tau} b_u \, du \right\vert^{\beta\vee 1} \leq C\, \E \int_{s }^{t}|b_u|^{\beta\vee 1} \, du.$$ The compensated Poisson integral is treated in $L^{1+\gamma}(\Omega)$ since $\beta \leq 1 + \gamma.$ The Burkholder-Davis-Gundy (BDG) inequality together with the fact that $1 +\gamma \leq 2$ imply that \begin{align}\label{BDGsmalljump}
	\notag\E \left|\int_0^{t \wedge \tau} \int_{B_1} H(r,z) \, \NNN(dr,dz) - \int_0^{s \wedge \tau} \int_{B_1} H(r,z) \, \NNN(dr,dz) \right|^{1+ \gamma}&\leq C\, \E \left(\int_{s \wedge \tau}^{t \wedge \tau} \int_{B_1}|H(r,z)|^2 \, \NN(dr,dz) \right)^{\frac{1+\gamma}{2}} \\ \notag&= C\, \E \left( \sum_{s\wedge \tau < r \leq t\wedge \tau} |H(r,\Delta Y_r)|^2\right)^{\frac{1+\gamma}{2}} \\ &\leq C\, \E \sum_{s\wedge \tau< r \leq t\wedge \tau} |H(r,\Delta Y_r)|^{1 +\gamma} \\ 	\notag&= C\, \E \int_{s \wedge \tau}^{t \wedge \tau} \int_{B_1} |H(r,z)|^{1+\gamma} \, \NN(dr,dz) \\ 	\notag&=C\, \E \int_{s \wedge \tau}^{t \wedge \tau} \int_{B_1} |H(r,z)|^{1+\gamma} \, d\nu(z) \, dr\\	\notag &\leq C\, \E \int_{s}^{t} \int_{B_1} |H(r,z)|^{1+\gamma} \, d\nu(z) \, dr.
\end{align}
 Finally, for the Poisson integral, we work in $L^{\beta}(\Omega)$ and we begin with the case $\beta \leq 1.$ In this case, the same computations as done for the compensated Poisson random integral yield \begin{align}\label{BDGbigjumpbetaleq1}
	\notag\E \left\vert \int_{s \wedge \tau}^{t \wedge \tau} \int_{B_1^c} K(r,z) \, \NN(dr,dz)\right\vert^{\beta} & = \E \left\vert \sum_{n=1}^{\infty} K(\tilde{T}_n,Z_n) \1_{s\wedge \tau <\tilde{T}_n \leq t\wedge \tau} \right\vert^{\beta}  \\ &\leq C\, \E \int_s^t \int_{B_1^c}|K(r,z)|^{\beta}  \, d\nu(z)\,dr. 
\end{align}
Let us now deal with the case $\beta \in (1,2].$ Writing artificially the Poisson integral as a compensated Poisson integral and using BDG's and Jensen's inequalities, we deduce that \begin{align}\label{BDGbigjumpsbeta>1}
	\notag&\E \left\vert\int_{s \wedge \tau}^{t \wedge \tau} \int_{B_1^c} K(r,z) \, \NN(dr,dz)\right\vert^{\beta} \\ \notag&\leq C \left[\E \left\vert \int_{s \wedge \tau}^{t \wedge \tau} \int_{B_1^c} K(r,z) \, \NNN(dr,dz) \right\vert^{\beta} + \E \left\vert \int_{s \wedge \tau}^{t \wedge \tau} \int_{B_1^c} K(r,z) \, d\nu(z)\,dr \right\vert^{\beta} \right] \\ &\leq  C \left[\E \left\vert \int_{s \wedge \tau}^{t \wedge \tau}\int_{B_1^c} |K(r,z)|^2 \, \NN(dr,dz) \right\vert^{\frac{\beta}{2}} + \E \int_{s \wedge \tau}^{t \wedge \tau} \int_{B_1^c} |K(r,z)|^{\beta} \, d\nu(z)\,dr \right] \\ \notag&\leq C \left[\E \int_{s \wedge \tau}^{t \wedge \tau} \int_{B_1^c} |K(r,z)|^{\beta} \, \NN(dr,dz)  + \E \int_{s \wedge \tau}^{t \wedge \tau} \int_{B_1^c} |K(r,z)|^{\beta} \, d\nu(z)\,dr  \right] \\\notag &\leq C\, \E \int_{s \wedge \tau}^{t \wedge \tau} \int_{B_1^c} |K(r,z)|^{\beta} \, d\nu(z)\,dr \\ \notag&\leq  C\, \E \int_{s}^{t} \int_{B_1^c} |K(r,z)|^{\beta} \, d\nu(z)\,dr. 
\end{align}

We have thus proved that \eqref{controlcontinuityflow} holds true. From this, we deduce that the map $t \mapsto \mu_t \in \PP$ is continuous using the dominated convergence theorem and Assumptions \textbf{(M)}, \textbf{(J1)} and \textbf{(J2)}.\\

For the moment estimate \eqref{momentX}, Jensen's inequality gives that $$ \sup_{t \leq T} \left\vert \int_0^t b_s \, ds \right\vert^{\beta \vee 1} \leq C\, \E \int_0^T |b_s|^{\beta \vee 1} \, ds.$$ The upper-bound is finite owing to Assumption \textbf{(M)}. For the compensated Poisson integral, we work in $L^{1 +\gamma}(\Omega)$ since $\beta \leq 1 + \gamma.$ Reasoning as in \eqref{BDGsmalljump}, it follows from BDG's inequality that \begin{align*}
	\E \sup_{t \leq T} \left\vert \int_0^t \int_{B_1} H(s,z) \, \NNN(ds,dz) \right\vert^{1+\gamma} &\leq C\, \E \left(\int_0^T \int_{B_1} |H(s,z)|^2 \, \NN(ds,dz) \right)^{\frac{1+\gamma}{2}} \\ &\leq C\, \E \int_0^T \int_{B_1} |H(s,z)|^{1 + \gamma} \, d\nu(z) \, ds,
	\end{align*}
which is finite thanks to Assumption \textbf{(J1)}. Finally, for the Poisson integral, we have $$ \E \sup_{t \leq T} \left\vert \int_0^t \int_{B_1^c} K(s,z) \, \NN(ds,dz) \right\vert^{\beta} \leq \E \left(\int_0^T\int_{B_1^c} |K(s,z)| \, \NN(ds,dz) \right)^{\beta}.$$ Reasoning as in \eqref{BDGbigjumpbetaleq1} if $\beta \leq1$ and as in \eqref{BDGbigjumpsbeta>1} if $\beta >1,$ one deduces that $$ \E \sup_{t \leq T} \left\vert \int_0^t \int_{B_1^c} K(s,z) \, \NN(ds,dz) \right\vert^{\beta} \leq C\, \E \int_0^T \int_{B_1^c} |K(s,z)|^{\beta} \, d\nu(z) \,ds.$$ The upper-bound is finite by Assumption \textbf{(J2)}.

\end{proof}
We now localize the process $X$. \begin{Def}[Localized process and approximate flow]\label{localizedprocessdef}
Let us introduce the sequence of localizing stopping times $(T_n)_n$ defined by $$ T_n := \inf\{t \in [0,T], \, |X_t| \geq n \, \text{or} \, |\widetilde{H}_s| \geq n\}\wedge T,$$ where $\widetilde{H}$ was introduced in Assumption $\textbf{(J2)}$. Then, we define the localized process $X^n$ and its flow of measures, for all $t \in [0,T],$ by \begin{equation}\label{localizedprocess}
\forall t \in [0,T], \, X^n_t := X_{t \wedge T_n} \quad \text{and} \quad \mu^n_t := \mathcal{L}(X^n_t).
\end{equation}
\end{Def}

\begin{Rq}\label{rqstoppingtime}
	Since $X$ is a càdlàg process and $\widetilde{H}$ is almost surely locally bounded, it is obvious that almost surely $T_n=T$ for $n$ bigger than a random constant. 
\end{Rq}

Let us start with the continuity of the sequence of approximate flows $(\mu^n)_n$ and a uniform moment estimate as in Proposition \ref{momentXprop}.

\begin{Prop}\label{momentlocalizedprocess}We have \begin{equation}\label{unifcontXn} \sup_n \E |X^n_t - X^n_s|^{\beta} \underset{s \rightarrow t}{\longrightarrow}0, 
	\end{equation}
which yields the continuity of $t\in [0,T] \mapsto \mu^n_t \in \PP$ for all $n\geq 1.$ Moreover, the following uniform moment estimate holds true \begin{equation}\label{momentXn}\E \sup_{n \geq 1}\sup_{t \leq T} |X^n_t|^{\beta} < + \infty. \end{equation}
\end{Prop}

\begin{proof}
	It follows immediately from Proposition \ref{momentXprop}.
\end{proof}

We can now state the main approximation result that we will use to prove Theorem \ref{itoformula}.

\begin{Prop}[Approximation]\label{propapprox}
	We have almost surely $\displaystyle\sup_{t \leq T} | X^n_t - X_t| \underset{n \rightarrow + \infty }{\longrightarrow}0.$ Moreover, for all $t \in [0,T],$ $\mu^n_t \overset{\PP}{\longrightarrow}\mu_t.$
		 	
\end{Prop}

\begin{proof}
The first point is a direct consequence of the fact that $T_n=T$ for $n$ large enough (see Remark \ref{rqstoppingtime}). For the second convergence, we show that $\E|X^n_t - X_t|^{\beta} \rightarrow 0$ for any $t\in [0,T].$ For the drift term, we work in $L^{\beta \vee 1}(\Omega)$ and apply the dominated convergence theorem. Indeed, $(T_n)_n$ converges to $T$ and the domination is a consequence of Assumption \textbf{(M)} by Jensen's inequality. We now deal with the Poisson integral which can be written as $$ \sum_{k=1}^{\infty} K(\tilde{T}_k,Z_k)\1_{\tilde{T}_k \leq t\wedge T_n}.$$ Since the sum is almost surely finite, it is clear that for $t \leq T$ $$ \sum_{k=1}^{\infty} K(\tilde{T}_k,Z_k) \1_{\tilde{T}_k \leq t\wedge T_n} \underset{n \rightarrow +\infty}{\longrightarrow} \sum_{k=1}^{\infty} K(\tilde{T}_k,Z_k) \1_{\tilde{T}_k \leq t}.$$ We conclude with the dominated convergence theorem. For the domination, one has \begin{align*}
 \left\vert \sum_{k=1}^{\infty} K(\tilde{T}_k,Z_k) \1_{\tilde{T}_k \leq t\wedge T_n} \right\vert &\leq  \sum_{k=1}^{\infty} |K(\tilde{T}_k,Z_k)| \1_{\tilde{T}_k \leq T}\\ & = \int_0^T \int_{B_1^c} |K(s,z)| \, \NN(ds,dz).\end{align*} The right-hand side term belongs to $L^{\beta}$ reasoning as in \eqref{BDGbigjumpbetaleq1} if $\beta \leq 1$ and as in \eqref{BDGbigjumpsbeta>1} if $\beta >1$ with Assumption \textbf{(J2)}. Eventually, for the compensated Poisson integral, we work in $L^{1+\gamma}(\Omega).$ It follows from BDG's inequality that \begin{align*} &\E \left\vert \int_0^{t\wedge T_n} \int_{B_1} H(s,z) \, \NNN(ds,dz) - \int_0^t \int_{B_1} H(s,z) \, \NNN(ds,dz)\right\vert^{1 + \gamma} \\ &\leq C\, \E \left(\int_t^{t\wedge T_n} \int_{B_1} |H(s,z)|^2 \, \NN(ds,dz) \, ds\right)^{\frac{1+\gamma}{2}} \\ &\leq C\, \E \int_t^{t\wedge T_n} \int_{B_1} |H(s,z)|^{1 + \gamma} \, d\nu(z) \, ds.
\end{align*}
The upper-bound tends to $0$ owing to the dominated convergence theorem and Assumption \textbf{(J1)}.
\end{proof}

We end this preliminary section with a compactness result on the sequence of approximate flows $(\mu^n)_n.$ 

\begin{Prop}\label{compactness}
	The set $\{\mu^n_t, \, n \geq 1,\, t \in[0,T] \}$ is relatively compact in $\PP.$ 
\end{Prop}

\begin{dem}
Let $(\mu^{n_k}_{s_k})_k$ be a subsequence of $\{\mu^n_t, \, n \geq 1,\, t \in[0,T] \}.$ Up to extraction, we can assume that $s_k \rightarrow s \in [0,T].$ If there exists $k_0$ such that $n_k = n_{k_0}$ for infinitely many $k \geq 1.$ The continuity of $t \in [0,T] \mapsto \mu^{n_{k_0}}_t$ allows to conclude that along this subsequence $$ \mu_{s_k}^{n_k} \overset{\PP}{\longrightarrow} \mu^{n_{k_0}}_s.$$  Otherwise, we can assume that $n_k \rightarrow + \infty.$ The triangle inequality yields $$ W_{\beta}(\mu_{s_k}^{n_k},\mu_s) \leq \sup_{n\geq 1} W_{\beta}(\mu^n_{s_k},\mu^n_s) + W_{\beta}(\mu^{n_k}_s,\mu_s).$$ The first term of the right-hand side of the preceding inequality converges to $0$ owing to Proposition \ref{momentlocalizedprocess}, as well as the second one because of Proposition \ref{propapprox}.
\end{dem}

\subsection{Proof of Itô's formulas of Theorem \ref{ito} and Theorem \ref{ito2}}

\begin{proof}[Proof of Theorem \ref{ito}]

In Itô's formula \eqref{itoformula}, there are three types of integrals: the drift term, the small jumps term coming from the compensated Poisson integral and the big jumps term stemming from the Poisson integral. They will be always treated separately in the proof. Let us fix $t \in (0,T]$ and prove Itô's formula at time $t.$\\

\textbf{Step 1: All the terms in Itô's formula \eqref{itoformula} are well-defined.}\\
	
	If $\beta >1$, the drift term is well-defined since Hölder's inequality implies that $$ \int_0^T \E \left\vert \partial_v \del (\mu_s)(X_s)\cdot b_s \right\vert \, ds \leq \left(\E\int_0^T |b_s|^{\beta} \, ds\right)^{1/{\beta}} \left(\E\int_0^T \left\vert \partial_v \del (\mu_s)(X_s) \right\vert^{\beta'} \, ds\right)^{1/{\beta'}}.$$ The right-hand side term is finite thanks to Assumption \textbf{(M)} and Assumption $(3)$ in Theorem \ref{ito} since $(\mu_t)_{t \in [0,T]}$ is compact in $\PP$ by Proposition \ref{momentXprop}. If $ \beta \leq 1,$ Assumption $(3)$ in Theorem \ref{ito} ensures that there exists $C>0$ such that $$ \int_0^T \E \left\vert \partial_v \del (\mu_s)(X_s)\cdot b_s \right\vert \, ds \leq C\, \E \int_0^T |b_s| \, ds,$$ which is finite thanks to Assumption \textbf{(M)}.  We consider now the big jumps term in \eqref{itoformula}. The growth condition in the definition of the linear derivative (see Definition \ref{deflinearderivative}) and Assumption \textbf{(J2)} ensure that $$ \E \left\vert  \del(\mu_s)(X_{s^-}+ K(s,z)) - \del(\mu_s)(X_{s^-}) \right\vert \leq C \left(1 + \E\sup_{s\leq T} |X_s|^\beta + |K(s,z)|^{\beta}\right).$$ The right-hand side term belongs to $L^1([0,T]\times B_1^c, ds\otimes \nu)$ because of the moment estimate \eqref{momentX} in Proposition \ref{momentXprop} and Assumption \textbf{(J2)}. We finally focus on the small jumps term. Note that if $\gamma = 0$, since we have assumed that $\beta \leq \gamma +1$, then $\beta \leq 1$ and Assumption $(3)$ in Theorem \ref{ito} ensures that for all compact $\KK \subset \PP$, there exists $C_{\KK}$ such that for all $v \in \R^d$, we have $$ \displaystyle\sup_{\mu\in \KK} \left\vert \partial_v \del (\mu)(v) \right\vert \leq C_{\KK}.$$ Thus, the mean value theorem gives that for all $s \in [0,T]$ and for all $z \in B_1$ \begin{align}\label{holdergammaequalzero}
		\left\vert \del(\mu_s)(X_{s^-}+H(s,z))  - \del (\mu_s)(X_{s^-}) - H(s,z) \cdot\partial_v\del (\mu_s)(X_{s^-}) \right\vert &\leq C |H(s,z)|\\ & \notag= C |H(s,z)|^{1 + \gamma}.\end{align} If $\gamma >0,$ then \eqref{holdergammaequalzero} holds true. Indeed, it follows from Taylor's formula and Assumption $(2)$ in Theorem \ref{ito} that for all $s \in [0,T]$ and $z \in B_1$ \begin{align}\label{holderbound}
	&\notag\left\vert \del(\mu_s)(X_{s^-}+H(s,z))  - \del (\mu_s)(X_{s^-}) - H(s,z) \cdot\partial_v\del (\mu_s)(X_{s^-}) \right\vert \\ &= \left\vert \int_0^1 \left[\partial_v\del(\mu_s)(X_{s^-}+rH(s,z)) - \partial_v\del (\mu_s)(X_{s^-})\right]\cdot H(s,z) \, dr \right\vert \\ \notag&\leq C|H(s,z)|^{1+\gamma}.
	\end{align}
	 Assumption \textbf{(J1)} allows us to conclude that the small jumps term is well-defined. \\ 
	
	\textbf{Step 2: Itô's formula for $X^n$.}\\
	
	We fix $n \geq 1$ and we aim at proving Itô's formula \eqref{itoformula} for the localized process $X^n$ defined for all $t \in [0,T]$ by $X^n_t = X_{t \wedge T_n}.$ For $m \geq 1,$ we define a subdivision of $[0,t]$ by $$ \forall k \in \{0,\dots,m\},\, t^m_k := \frac{k}{m}t.$$ We will omit the index $m$ and denote by $(t_k)_k$ this subdivision. By definition of the linear derivative, one has for all $k \in \{0, \dots, m\} $ \begin{align}\label{eq1preuve}
	u(\mu^n_{t_{k+1}}) - u(\mu^n_{t_k}) = \int_0^1 \E \left( \del (M^k_r) ( X^n_{t_{k+1}}) -  \del (M^k_r) ( X^n_{t_{k}})\right) \, dr,
\end{align}
where $M_r^k := r \mu^n_{t_{k+1}} + (1-r) \mu^n_{t_k},$ omitting again the dependence in $n$ and $m.$ Let us denote by $F$ the function $\del (M^k_r).$ We can apply the standard Itô formula for $F$ and $X^n$ which yields for all $t\in [0,T]$ \begin{align}\label{eqItoXn}
	\notag F(X^n_t) &= F(X^n_0) + \int_0^{t \wedge T_n} \nabla F (X^n_s)\cdot b_s \, ds \\ &\quad + \int_0^{t\wedge T_n} \int_{B_1^c} \, \left[F(X^n_{s^-} + K(s,z)) - F(X^n_{s^-})\right] \, \NN(ds,dz) \\\notag &\quad + \int_0^{t\wedge T_n} \int_{B_1} \, \left[F(X^n_{s^-} + H(s,z)) - F(X^n_{s^-})\right] \, \NNN(ds,dz) \\ \notag&\quad +  \int_0^{t\wedge T_n} \int_{B_1} \, \left[F(X^n_{s^-} + H(s,z)) - F(X^n_{s^-}) - H(s,z)\cdot \nabla F(X^n_{s^-})\right] \, d\nu(z) \, ds.
\end{align}

Now, we aim at taking the expectation in the previous formula. The compensated Poisson integral is centered. Indeed, almost surely for all $s \in [0,t\wedge T_n)$ and for all $z \in B_1$ $$|X^n_{s^-}| + |H(s,z)| \leq 2n,$$ by definition of $T_n.$ Thus, the mean value theorem ensures that for some constant $C_n>0$ depending only on $n$ and $F,$ we have almost surely for all $s \in [0,t\wedge T_n)$ and for all $z \in B_1$ \begin{align*} |F(X^n_{s^-} + H(s,z)) - F(X^n_{s^-})|^2 &\leq C_n |H(s,z)|^2.\end{align*} Thus, by definition of $T_n$ and by Assumption \textbf{(J1)}, one has almost surely for all $s \in [0,t\wedge T_n)$ and for all $z \in B_1$ $$ |F(X^n_{s^-} + H(s,z)) - F(X^n_{s^-})|^2 \leq n^2C_n|z|^2.$$
From this inequality, we deduce that $$ \E \int_0^{t\wedge T_n} \int_{B_1} \, \left\vert F(X^n_{s^-} + H(s,z)) - F(X^n_{s^-})\right\vert^2 \, d\nu(z) \, ds < + \infty.$$ It follows that the compensated Poisson integral is a true centered martingale. For the Poisson integral, we use the growth property of $F$ coming from the definition of the linear derivative. One has $$ |F(X^n_{s^-} + K(s,z)) - F(X^n_{s^-})| \leq C \left(1+\sup_{s\leq T} |X_s|^{\beta} + |K(s,z)|^{\beta}\right).$$ Thus, Assumption \textbf{(J2)} and the moment estimate \eqref{momentX} in Proposition \ref{momentXprop} prove that $$ \E \int_0^{t\wedge T_n} \int_{B_1^c} \, \left\vert F(X^n_{s^-} + K(s,z)) - F(X^n_{s^-})\right\vert\, d\nu(z) \, ds < + \infty. $$ This yields \begin{align*}
	 &\E \int_0^{t\wedge T_n} \int_{B_1^c} \, \left[F(X^n_{s^-} + K(s,z)) - F(X^n_{s^-})\right] \, \NN(ds,dz) \\ &=  \E\int_0^{t\wedge T_n} \int_{B_1^c} \, \left[F(X^n_{s^-} + K(s,z)) - F(X^n_{s^-})\right] \, d\nu(z)\,ds.\end{align*}  Taking the expectation in the other term of Itô's formula \eqref{eqItoXn} can be done arguing as in Step $1$. It follows that \begin{align}\label{eqItoXnexpectation}
	 \E F(X^n_t) &= \E F(X^n_0) + \E\int_0^{t \wedge T_n} \partial_vF (X^n_s)\cdot b_s \, ds \\ \notag&\quad + \E\int_0^{t\wedge T_n} \int_{B_1^c} \, \left[F(X^n_{s^-} + K(s,z)) - F(X^n_{s^-})\right] \, d\nu(z) \,ds  \\ \notag&\quad +  \E \int_0^{t\wedge T_n} \int_{B_1} \, \left[F(X^n_{s^-} + H(s,z)) - F(X^n_{s^-}) - H(s,z)\cdot \partial_v F(X^n_{s^-})\right] \, d\nu(z) \, ds.
	 \end{align}
	 
Recall that $F$ is equal to $\del (M^k_r).$ Thus, by using \eqref{eq1preuve} and \eqref{eqItoXnexpectation}, one has \begin{align}\label{eq2preuve}
	u(\mu^n_t) - u(\mu^n_0) &= \sum_{k=0}^{m-1} u(\mu^n_{t_{k+1}}) - u(\mu^n_{t_k}) \\\notag &=I_1 + I_2 + I_3,
\end{align}
where 
\begin{align*}
&I_1:= \sum_{k=0}^{m-1} \int_0^1 \E \int_{t_k \wedge T_n}^{t_{k+1} \wedge T_n} \partial_v \del (M^k_r)(X^n_s) \cdot b_s \, ds \, dr, \\ &I_2 := \sum_{k=0}^{m-1} \int_0^1 \E \int_{t_k \wedge T_n}^{t_{k+1} \wedge T_n} \int_{B_1^c} \left[ \del (M_r^k) (X^n_{s^-} + K(s,z)) - \del(M_r^k)(X^n_{s^-})\right] \, d\nu(z) \, ds \, dr,\\ &I_3 := \sum_{k=0}^{m-1} \int_0^1 \E \int_{t_k \wedge T_n}^{t_{k+1} \wedge T_n} \int_{B_1}  \left[\del (M_r^k) (X^n_{s^-} + H(s,z)) - \del(M_r^k)(X^n_{s^-}) \right. \\ &\hspace{5cm}\left. - H(s,z)\cdot \partial_v \del (M^k_r)(X^n_{s^-})\right] \, d\nu(z) \, ds \, dr.
\end{align*} Our goal is now to let $m$ tend to infinity in \eqref{eq2preuve}. For the small jumps term $I_3$, it can be rewritten as \begin{align*}
 I_3 = \int_0^1 \E \int_0^{t\wedge T_n} \int_{B_1}&\left[\del (M_r^{\lfloor ms/t\rfloor}) (X^n_{s^-} + H(s,z)) - \del(M_r^{\lfloor ms/t\rfloor})(X^n_{s^-}) \right. \\ & \left. -H(s,z)\cdot \partial_v \del (M^{\lfloor ms/t\rfloor}_r)(X^n_{s^-})\right] \, d\nu(z) \, ds\,dr,\end{align*} since $\lfloor ms/t\rfloor$ is equal to $k$ if and only if $s \in [t_k,t_{k+1}).$ The continuity of the flow $s \mapsto \mu^n_s \in \PP$ ensures that for all $r \in [0,1]$ and $s \in [0,t]$ $$ M^{\lfloor ms/t\rfloor}_r \underset{m \rightarrow + \infty}{\longrightarrow} \mu^n_s.$$ The continuity of $\del$ and $\partial_v \del$ on $\PP \times \R^d$ implies that \begin{align*}
 	& \del (M_r^{\lfloor ms/t\rfloor}) (X^n_{s^-} + H(s,z)) - \del(M_r^{\lfloor ms/t\rfloor})(X^n_{s^-}) -H(s,z)\cdot \partial_v \del (M^{\lfloor ms/t\rfloor}_r)(X^n_{s^-}) \\ &\underset{m\rightarrow +\infty}{\longrightarrow} \del (\mu^n_s) (X^n_{s^-} + H(s,z)) - \del(\mu^n_s)(X^n_{s^-}) -H(s,z)\cdot \partial_v \del (\mu^n_s)(X^n_{s^-}).
 \end{align*}
 Note that the set $\left\{v\mu_s^n + (1-v)\mu_t^n, \, s,t \in[0,T],\, v \in[0,1]\right\}$ is compact in $\PP$ thanks to Proposition \ref{compactness}. As done at the end of Step 1 in \eqref{holdergammaequalzero} and \eqref{holderbound}, Assumption $(2)$ in Theorem \ref{ito} implies that the integrand in $I_3$ is bounded, up to some multiplicative constant, by $ |\widetilde{H}_s|^{1+\gamma} |z|^{1+\gamma}.$ Assumption \textbf{(J1)} proves that this quantity belongs to $L^1([0,1]\times[0,t]\times B_1 \times \Omega, dr\otimes ds\otimes \nu \otimes \P).$ Thus, the dominated convergence theorem ensures that $I_3$ converges, when $m \rightarrow + \infty,$ towards $$\E \int_0^{t\wedge T_n} \int_{B_1} \left[ \del (\mu^n_s) (X^n_{s^-} + H(s,z)) - \del(\mu^n_s)(X^n_{s^-}) -H(s,z)\cdot \partial_v \del (\mu^n_s)(X^n_{s^-})\right] \, d\nu(z) \, ds.$$ Using the same arguments, we prove that $$ I_2 \underset{m \rightarrow + \infty}{\longrightarrow} \E \int_0^t \int_{B_1^c}  \left[ \del (\mu^n_s) (X^n_{s^-} + K(s,z)) - \del(\mu_s^n)(X^n_{s^-})\right] \, d\nu(z) \, ds.$$ From the growth condition in the definition of the linear derivative, we deduce that almost surely for all $s \in [0,T],\, r \in [0,1],\, z \in B_1,\, n \in \N$ $$ \left\vert \del (M_r^k) (X^n_{s^-} + K(s,z)) - \del(M_r^k)(X^n_{s^-}) \right\vert \leq C \left(1 + \sup_{s \leq T} |X_s|^{\beta} + |K(s,z)|^\beta\right).$$ The upper-bound belongs to $L^1([0,1]\times[0,t]\times B_1^c \times \Omega, dr\otimes ds\otimes \nu \otimes \P)$ by Assumption \textbf{(J2)}. This justifies the use of the dominated convergence theorem, as we did for $I_3.$ Finally, $I_1$ can be handled in the same way, the domination being an immediate consequence of the localization of the process and the continuity of $\partial_v \del.$ Letting $m$ tend to infinity in \eqref{eq2preuve}, one has for all $t \in [0,T]$ 
 
 \begin{align}\label{itoapprox}
 &\notag u(\mu^n_t) - u(\mu^n_0) \\ &=\E\int_0^{t \wedge T_n}  \left( \partial_v \del (\mu^n_s)(X^n_s) \cdot b_s \right) \, ds \\\notag &\quad+  \E\int_0^{t \wedge T_n} \int_{B_1^c}  \left[ \del(\mu^n_s)(X^n_{s^-}+ K(s,z)) - \del(\mu^n_s)(X^n_{s^-})\right] \, d\nu(z) \, ds  \\ \notag&\quad +  \E\int_0^{t \wedge T_n} \int_{B_1}  \left[ \del(\mu^n_s)(X^n_{s^-}+ H(s,z)) - \del(\mu^n_s)(X^n_{s^-}) - \partial_v \del(\mu^n_s)(X^n_{s^-})\cdot H(s,z)\right]  d\nu(z) \, ds  \\  \notag&=: A_1 + A_2 + A_3. \end{align}

 \vspace{2pt}
 
\textbf{Step 3: Itô's formula for $X.$}\\

 Our goal is now to let $n$ tend to infinity in Itô's formula \eqref{itoapprox} for the localized process. For the left-hand side term of \eqref{itoapprox}, we have already seen in Proposition \ref{propapprox} that for all $t \in [0,T]$ $$ \mu^n_t \overset{\PP}{\longrightarrow} \mu_t.$$ The continuity of $u$ on $\PP$ ensures that $u(\mu^n_t) \rightarrow u(\mu_t).$ For the right-hand side of \eqref{itoapprox}, let us start with the limit when $n$ tends to  $+ \infty$ of $A_1.$ It follows from the continuity of $\partial_v \del$ and from Proposition \ref{propapprox} that almost surely for all $s \leq t$ $$ \partial_v \del (\mu^n_s)(X^n_s) \cdot b_s \1_{s \leq T_n} \underset{n \rightarrow
  + \infty}{\longrightarrow} \partial_v \del (\mu_s)(X_s) \cdot b_s.$$
For $\beta\leq 1,$ we use the dominated convergence Theorem justified by Assumption $(3)$ in Theorem \ref{ito} and Assumption \textbf{(M)}. For $\beta >1,$ we conclude with a uniform integrability argument in the space $L^{1}([0,t] \times \Omega, ds \otimes \P).$ Indeed, note that $b \in L^{\beta}([0,t] \times \Omega)$ and $$\sup_{n \geq 1} \E \int_0^t \left\vert \partial_v \del(\mu^n_s)(X^n_s)\1_{s \leq T_n} \right\vert^{\beta '} \, ds <+\infty,$$ thanks to Assumption $(3)$ of Theorem \ref{ito} and the relative compactness of $\{\mu^n_s, \, n \geq 1, \, s \leq T \}$ of Proposition \ref{compactness}. Thus, the sequence $\left( \partial_v \del (\mu^n_s)(X^n_s) \cdot b_s \1_{s \leq T_n}\right)_n$ is uniformly integrable in $L^1([0,t] \times \Omega, ds \otimes \P).$ We have thus proved that $$ A_1 \underset{n \rightarrow +\infty}{\longrightarrow} \E \int_0^t \partial_v \del (\mu_s)(X_s) \cdot b_s \, ds.$$
We now focus on the big jumps term $A_2.$ The continuity of $\del$ and Proposition \ref{propapprox} prove that almost surely for all $s \in [0,t],\, z \in B_1^c$ $$ \del(\mu^n_s)(X^n_{s^-}+ K(s,z)) - \del(\mu^n_s)(X^n_{s^-}) \underset{n\rightarrow + \infty}{\longrightarrow} \del(\mu_s)(X_{s^-} + K(s,z)) - \del (\mu_s)(X_{s^-}).$$
The growth assumption in the definition of the linear derivative together with Assumption \textbf{(J2)} ensure that for all $s \in [0,t],\, z \in B_1^c, \, n \in \N$ $$ \left\vert \del(\mu^n_s)(X^n_{s^-}+ K(s,z)) - \del(\mu^n_s)(X^n_{s^-}) \right\vert  \leq C \left(1+\sup_{s\leq T } |X_s|^{\beta} + |K(s,z)|^{\beta}\right).$$Using Assumption \textbf{(J2)} and the moment estimate \eqref{momentX}, one has $$ \E \int_0^t \int_{B_1^c} \left(1+\sup_{s\leq T } |X_s|^{\beta} + |K(s,z)|^{\beta}\right) \, d\nu(z) \,ds < + \infty.$$ Thus, the dominated convergence theorem yields $$A_2 \underset{n\rightarrow +\infty}{\longrightarrow} \E\int_0^{t} \int_{B_1^c}  \left[ \del(\mu_s)(X_{s^-}+ K(s,z)) - \del(\mu_s)(X_{s^-})\right] \, d\nu(z) \, ds.$$ 

Eventually, we let $n$ tend to infinity in the small jumps term $A_3$ thanks to the dominated convergence theorem. For the domination, we use the compactness of $(\mu^n)_n$ given by Proposition \ref{compactness} and Assumption $(2)$ in Theorem \ref{ito}. Reasoning as in the end of Step 1 in \eqref{holdergammaequalzero} and \eqref{holderbound}, we have almost surely, for all  $s \in [0,t],\, z \in B_1, \, n \in \N$ $$ \left\vert\del(\mu^n_s)(X^n_{s^-}+ H(s,z)) - \del(\mu^n_s)(X^n_{s^-}) - \partial_v \del(\mu^n_s)(X^n_{s^-})\cdot H(s,z) \right\vert \leq C|H(s,z)|^{1+\gamma}.$$
It follows from Assumption \textbf{(J1)} that the dominated convergence theorem can be applied. It yields $$ A_3 \underset{n \rightarrow + \infty}{\longrightarrow}  \E\int_0^{t \wedge T_n} \int_{B_1}  \left[ \del(\mu_s)(X_{s^-}+ H(s,z)) - \del(\mu_s)(X_{s^-}) - \partial_v \del(\mu_s)(X_{s^-})\cdot H(s,z)\right]  d\nu(z) \, ds .$$ We conclude the proof by taking the limit $n \rightarrow + \infty$ in \eqref{itoapprox}.
\end{proof}

\begin{proof}[Proof of Theorem \ref{ito2}]
	It is enough to prove that the map $(t,x) \in [0,T]\times \R^d \mapsto u(t,x,\mu_t)$ is of class $\CC^{1}$. Then, we conclude by applying the standard Itô's formula for this function and for the process $(Y_t)_t $ since $\partial_x u(t,x,\mu_t)$ is $\gamma'$-Hölder continuous uniformly in time. Let us fix $\KK \subset \R^d$ a compact subset, $x \in \KK$, $t \in [0,T]$ and $h \in \R$ such that $t+h \in [0,T].$ We have $$ u(t+h,x,\mu_{t+h}) - u(t,x,\mu_t) = \left[u(t+h,x,\mu_{t+h}) -  u(t,x,\mu_{t+h})\right] + \left[ u(t,x,\mu_{t+h}) - u(t,x,\mu_t)\right].$$ The continuity of $\partial_t u $ on $[0,T]\times \R^d \times \PP$ ensures that \begin{align*}
			\frac1h u(t+h,x,\mu_{t+h}) -  u(t,x,\mu_{t+h}) & = \frac1h \int_{t}^{t+h} \partial_t u(s,x,\mu_{t+h}) \, ds \\ &\underset{h \rightarrow 0}{\longrightarrow} \partial_t u ( t,x,\mu_t). 
		\end{align*}
		
		Then, using Itô's formula of Theorem \ref{ito}, we obtain that \begin{align}\label{extensionproof1}
			\notag&u(t,x,\mu_{t+h}) - u(t,x,\mu_t) \\&= \int_t^{t+h} \E \left( \partial_v \del (t,x,\mu_s)(X_s) \cdot b_s \right) \, ds \\\notag &\quad+ \int_t^{t+h} \int_{B_1^c} \E \left[ \del(t,x,\mu_s)(X_{s^-}+ K(s,z)) - \del(t,x,\mu_s)(X_{s^-})\right] \, d\nu(z) \, ds  \\ \notag&\quad +  \int_t^{t+h} \int_{B_1} \E \left[ \del(t,x,\mu_s)(X_{s^-}+ H(s,z))  \del(t,x,\mu_s)(X_{s^-}) \right.\\ &\left.\notag \hspace{8cm} -\partial_v \del(t,x,\mu_s)(X_{s^-})\cdot H(s,z)\right] \, d\nu(z) \, ds. 
		\end{align}
		
		First, note that the function $(s,t,x) \in [0,T]^2\times \KK\mapsto \E \left( \partial_v \del (t,x,\mu_s)(X_s) \cdot b_s \right)$ is continuous. Indeed, it easily follows from the continuity of $\partial_v \del$ and $b,$ the fact that for all $t \in [0,T],$ $X_{t^-} = X_t$ almost surely, and from a uniform integrability argument. The uniform integrability of the family $\left(\partial_v \del (t,x,\mu_s)(X_s) \cdot b_s \right)_{s,t,x \in [0,T]^2\times\KK}$ comes from the integrability assumption on $b$ and from Assumption $(4)$ in Theorem \ref{ito2}. It follows from the dominated convergence theorem that the function $$(s,t,x)\in [0,T]^2\times \KK \mapsto  \int_{B_1^c} \E \left[ \del(t,x,\mu_s)(X_{s^-}+ K(s,z)) - \del(t,x,\mu_s)(X_{s^-})\right] \, d\nu(z)$$ is continuous. Indeed, the domination is easily deduced from Assumption $(2)$ in Theorem \ref{ito2}, the assumption on $K$ and from \eqref{momentX}. Finally, the function \begin{align*} & (s,t,x)\in [0,T]^2\times \KK \mapsto  \\ &\int_{B_1}\E \left[ \del(t,x,\mu_s)(X_{s^-}+ H(s,z)) - \del(t,x,\mu_s)(X_{s^-}) - \partial_v \del(t,x,\mu_s)(X_{s^-})\cdot H(s,z)\right] \, d\nu(z)\end{align*}
		 is also continuous using the dominated convergence theorem. Indeed by, Assumption $(3),$  we find that almost surely for all $x \in \KK$, $z \in B_1,$ $s,t \in [0,T]$ $$  \left\vert \del(t,x,\mu_s)(X_{s^-}+ H(s,z)) - \del(t,x,\mu_s)(X_{s^-}) - \partial_v \del(t,x,\mu_s)(X_{s^-})\cdot H(s,z)\right\vert \leq C \sup_{s \leq T} (|\widetilde{H}_s||z|)^{1 + \gamma}.$$
		This shows that \begin{align*}
			&\frac{1}{h} (u(t,x,\mu_{t+h}) - u(t,x,\mu_t)) \\ &\underset{h\to 0}{\longrightarrow} \int_{B_1^c} \E \left[ \del(t,x,\mu_t)(X_{t^-}+ K(t,z)) - \del(t,x,\mu_t)(X_{t^-})\right] \, d\nu(z)  \\ &\quad +  \int_{B_1} \E \left[ \del(t,x,\mu_t)(X_{t^-}+ H(t,z)) - \del(t,x,\mu_t)(X_{t^-}) - \partial_v \del(t,x,\mu_t)(X_{t^-})\cdot H(t,z)\right] \, d\nu(z). 
		\end{align*}
	We have thus proved that $t \in [0,T] \mapsto u(t,x,\mu_t)$ is differentiable and that $\frac{d}{dt} u(t,\cdot,\mu_t)$ is continuous on $[0,T]\times \R^d$ and satisfies Itô's formula \eqref{itoformula2}.

	\end{proof}

\section{Backward Kolmogorov PDE on the space of measures and empirical projection}\label{sectionPDE}
 In this section, we use our Itô's formula to derive the backward Kolmogorov PDE on the space of measures associated with a general Lévy-driven McKean-Vlasov SDE. More precisely, it describes the dynamics of the associated semigroup acting on functions defined on the space of measures under regularity assumptions on it. In the Brownian case, it has been done in Chapter $5$ of \cite{CarmonaProbabilisticTheoryMean2018} and in \cite{deraynal2021wellposedness}. Let us introduce $Z=(Z_t)_t$ a Lévy process on $\R^d$ which has the following Lévy-Itô decomposition \begin{equation}\label{deflevynoise}
Z_t = \int_0^t \int_{B_1} z \, \NNN(ds,dz) + \int_0^t \int_{B_1^c} z \, \NN(ds,dz),\end{equation}  where $\NN$ is the associated Poisson random measure with Lévy measure $\nu.$ We fix our assumptions on the Lévy process $Z$. We assume that there exists $\beta \in (0,2]$ and $\gamma \in [0,1]$ with $\beta \leq 1 + \gamma$ and such that the Lévy measure $\nu$ satisfies the following properties. 

\begin{enumerate}
	
	\item[] \textbf{(J1')} We have \begin{equation}\label{assumptionJ1Levymeasure} \int_{B_1} |z|^{1 + \gamma} \, d\nu(z) < + \infty. \end{equation}
	\item[]\textbf{(J2')} For all $t \in [0,T],$ $Z_t$ has a finite moment of order $\beta$ i.e.\ \begin{equation}\label{assumptionJ2Levymeasure}
	\int_{B_1^c} |z|^\beta \, d\nu(z) < + \infty.\end{equation}
\end{enumerate}
Note that we can always take $\gamma =1$ by definition of a Lévy measure. We are interested in the following McKean-Vlasov SDE on the time interval $[0,T]$ \begin{equation}\label{edsmckv}
\left\{  \begin{array}{lll}
&dX_t = b_t(X_t,\mu_t) \,dt + \sigma_t(X_{t^-},\mu_t)\, dZ_t, \\ &\mu_t = \mathcal{L}(X_t),\\ &X_0 = \xi \in L^{\beta}(\Omega,\FF_0),
\end{array}\right.
\end{equation}
where $b:[0,T]\times \R^d \times \PPP_{\beta}(\R^d) \rightarrow \R^d$ and $\sigma : [0,T]\times \R^d \times \PPP_{\beta}(\R^d)\rightarrow \R^{d \times d}$ are measurable.\\

We make the following assumptions on the McKean-Vlasov SDE \eqref{edsmckv}.

\begin{enumerate}
	\item[] \textbf{(H1)} There is weak existence and uniqueness for \eqref{edsmckv} for any initial distribution $\mu_0$ in $\PPP_{\beta}(\R^d).$ 
	
	\item[] \textbf{(H2)} If $\beta \geq 1,$ there exists $C>0$ such that for all $t \in [0,T],$ $x \in \R^d$ and $\mu \in \PPP_{\beta}(\R^d)$ \begin{align}\label{growthassumptionbetageq1} \notag|b_t(x,\mu)| &\leq C \left( 1 + |x| + M_{\beta}(\mu)\right), \\ |\sigma_t(x,\mu)| &\leq C (1 + M_{\beta}(\mu)),\end{align}
	where $M_{\beta}(\mu) = \left(\int_{\R^d} |x|^{\beta} \, d\mu(x)\right)^{\frac{1}{\beta}}.$\\
	
	If $\beta < 1,$ there exists $C>0$ such that for all $t \in [0,T],$ $x \in \R^d$ and $\mu \in \PPP_{\beta}(\R^d)$ \begin{align}\label{growthassumptionbetaleq1} \notag|b_t(x,\mu)| &\leq C, \\ |\sigma_t(x,\mu)| &\leq C.\end{align}
	
	\item[] \textbf{(H3)} The coefficients $b$ and $\sigma$ are continuous on $[0,T]\times \R^d \times \PPP_{\beta}(\R^d).$
\end{enumerate}
Note that under these assumptions, Lemma $4.1$ in \cite{FrikhaWellposednessapproximationonedimensional2020} ensures that $(\mu_t)_t$ belongs to $\CC^0([0,T];\PPP_{\beta}(\R^d))$ and that \begin{equation}
\E \sup_{t \in [0,T]} |X_t|^{\beta} < + \infty.
\end{equation}


The action of the semigroup associated with \eqref{edsmckv} is given in the following definition.

\begin{Def}[Semigroup]\label{defsemigroup}
For a function $u : \PPP_{\beta}(\R^d) \rightarrow \R,$ the action of the semigroup associated with the McKean-Vlasov SDE \eqref{edsmckv} on the map $u$ is given by $\phi_u$ defined by \begin{equation}
	\phi_u : \left\{ \begin{array}{rrl}
	&[0,T] \times \PPP_{\beta}(\R^d) & \rightarrow \R \\ &(t,\mu) &\mapsto u(\mathcal{L}(X_T^{T-t,\mu})),
	\end{array}\right.
\end{equation}
where $\mathcal{L}(X_T^{T-t,\mu})$ denotes the distribution of any solution to \eqref{edsmckv} at time $T$ starting at time $T-t$ from a random variable $\xi$ with distribution $\mu.$ This makes sense by the well-posedness assumption \textbf{(H1)} for \eqref{edsmckv}.
\end{Def}

\begin{Thm}[Backward Kolmogorov PDE]\label{ThmEDP} Assume that \textbf{(H1)}, \textbf{(H2)}, \textbf{(H3)} are satisfied and that the function $\phi_u$ satisfies the following properties.
	
	 \begin{enumerate}
		\item The function $\phi_u$ belongs to $\CC^0([0,T]\times \PPP_{\beta}(\R^d);\R),$  $\partial_t \phi_u $ exists and is continuous on $(0,T] \times \PPP_{\beta}(\R^d),$ and finally $\dell$ and $\partial_v \dell$ exist and are continuous on $(0,T]  \times \PPP_{\beta}(\R^d) \times \R^d.$

		\item If $\gamma>0,$ for any compact $\KK \subset \PP$ and $a>0,$ there exists $C>0$ such that for all $\mu \in \KK,$ $ x,y \in \R^d$ and $t \in [a,T]$ $$ \displaystyle\left\vert \partial_v\dell(t,\mu)(x) - \partial_v\dell(t,\mu)(y)  \right\vert \leq C |x-y|^{\gamma}.$$

		\item For any compact $\KK \subset \PP$ and $a>0,$ we have  $$\left\{ \begin{array}{lll} &\displaystyle\sup_{t \in[a,T]}\sup_{\mu \in \KK} \int_{\R^d} \left\vert \partial_v \dell (t,\mu)(v) \right\vert^{\beta'} \, d\mu(v) < + \infty \quad \text{if ${\beta >1},$}  \\ &  \displaystyle\sup_{t\in [a,T]}\sup_{\mu\in \KK} \sup_{v \in \R^d} \left\vert \partial_v \dell (t,\mu)(v) \right\vert < + \infty \quad  \text{if $\beta \leq1.$} \end{array} \right.$$
	\end{enumerate}
	Then, the function $\phi_u$ satisfies for any $t \in [0,T)$ and $\mu \in \PPP_{\beta}(\R^d)$ \begin{equation}\label{PDE}
		\left\{ \begin{aligned}
		&\partial_t \phi_u (T-t,\mu) = \mathscr{L}_t\phi_u(T-t,\mu), \\ &\phi_u(T-t,\mu)_{|t=T} = u(\mu),
		\end{aligned}\right.
	\end{equation}
	where the operator $\mathscr{L}_t$ is given, for any $s \in (0,T]$ and $ \mu \in \PPP_{\beta}(\R^d)$, by \begin{align}\label{operatormeasure}
\notag	&\mathscr{L}_t \phi_u(s, \mu)\\ \notag& := \int_{\R^d} \partial_v \dell (s,\mu)(v) \, b_t(v,\mu) \, d\mu(v) \\ &\quad + \int_{\R^d} \int_{B_1^c}  \left[\dell (s ,\mu)(v + \sigma_t(v,\mu)z) - \dell (s,\mu)(v)\right]  \, d\nu(z) \, d\mu(v)\\ \notag&\quad + \int_{\R^d} \int_{B_1} \left[\dell (s,\mu)(v + \sigma_t(v,\mu)z) - \dell (s ,\mu)(v) -  \partial_v \dell (s,\mu)(v)\cdot(\sigma_t(v,\mu)z)\right]  \, d\nu(z) \, d\mu(v).
	\end{align}
\end{Thm}

\begin{proof}First, note that by definition of $\phi_u$ and thanks to the well-posedness assumption \textbf{(H1)}, the function $ t \in [0,T] \mapsto \phi_u(T-t,\mu_t)$ is constant.  Differentiating with respect to $t$ the function $t \mapsto \phi(T-t,\mu_t)$ thanks to Itô's formula of Theorem \ref{ito2}, one has for all $t \in [0,T)$  \begin{align*}
		0 &= -\partial_t \phi_u(T-t,\mu_t) \\ &\quad + \int_{\R^d} \partial_v \dell (T-t,\mu_t)(v) \, b_t(v,\mu) \, d\mu_t(v) \\ \notag&\quad + \int_{\R^d} \int_{B_1^c}  \left[\dell (T-t ,\mu_t)(v + \sigma_t(v,\mu_t)z) - \dell (T-t,\mu_t)(v)\right]  \, d\nu(z) \, d\mu_t(v)\\ \notag&\quad + \int_{\R^d} \int_{B_1} \left[\dell (T-t,\mu_t)(v + \sigma_t(v,\mu_t)z) - \dell (T-t ,\mu_t)(v) \right. \\ &\left. \hspace{4cm}-  \partial_v \dell (T-t,\mu_t)(v)\cdot(\sigma_t(v,\mu_t)z)\right]  \, d\nu(z) \, d\mu_t(v).
	\end{align*}
	
	We conclude the proof by noting that the flow of measures $(\mu_t)_t$ associated with the McKean-Vlasov SDE \eqref{edsmckv} can be initialised at time $t$ with any measure $\mu \in \PP.$ 
\end{proof}

Let us recall the definition of empirical projection.

\begin{Def}[Empirical projection]\label{def_empirical_proj}
		Fix $u : \PPP_{\beta}(\R^d) \rightarrow \R$. For all $N \geq 1$, the empirical projection $u^N$ of $u$ is defined for all $\bm{x} = (x_1, \dots,x_N) \in (\R^d)^N$ by $$ u^N(\bm{x}) = u (\muu^N_{\bm{x}}),$$
		where $\muu^N_{\x} = \frac1N \displaystyle\sum_{j=1}^N \delta_{x_j}.$
\end{Def}

We also state the following Lemma on the regularity of the empirical projection. The proof is completely analogous to that of Proposition $5.91$ of \cite{CarmonaProbabilisticTheoryMean2018}.

\begin{Lemme}\label{projempC2}
	Let $ u : \PPP_{\beta}(\R^d) \rightarrow \R$ be a function admitting a linear derivative $ \del $ satisfying the following properties. \begin{enumerate}
		\item For all $\mu \in \PPP_{\beta}(\R^d)$, $\del (\mu) \in \CC^2(\R^d;\R)$.
		
		\item The functions $\partial_v \del$ and $ \partial^2_v \del$ are continuous on $\PPP_{\beta}(\R^d) \times \R^d.$
		
		\item For all $v \in \R^d$, the function $ \mu \in \PPP_{\beta}(\R^d) \mapsto \partial_v \del (\mu)(v)$ has a linear derivative $$ (\mu,v') \mapsto \frac{\delta}{\delta m} \partial_{v}\del (\mu)(v,v') \in \R^{d},$$ which is $\CC^1$ with respect to $v'$ and such that $\partial_{v'} \frac{\delta}{\delta m} \partial_{v}\del$ is continuous on $\PPP_{\beta}(\R^d) \times \R^d \times \R^d$. 
	\end{enumerate}
	
	Then, for all $N \geq 1$, the empirical projection $u^N$ of $u$ is of class $\CC^2$. Moreover, we have for all $\bm{x}=(x_1, \dots, x_N) \in (\R^d)^N$ and $i,j \in \{ 1,\dots N\}$ $$\partial_{x_i} u^N(x_1, \dots, x_N) = \frac{1}{N}  \partial_v \del (\overline{\mu}^N_{\bm{x}})(x_i),$$ and $$ \partial_{x_j} \partial_{x_i} u^N(x_1, \dots, x_N) = \frac{1}{N^2} \partial_{v'} \frac{\delta}{\delta m} \partial_v \del (\overline{\mu}^N_{\bm{x}})(x_i,x_j)  + \1_{i=j} \frac{1}{N} \partial^2_x \del (\overline{\mu}^N_{\bm{x}}) (x_j ).$$
\end{Lemme}

We consider now the mean-field interacting particle system  associated with the McKean-Vlasov SDE \eqref{edsmckv}. The goal is to approximate in some sense $\mathscr{L}_t u (\overline{\mu}^N_{\bm{x}})$ by the generator of the particle system applied to the empirical projection $u^N$ of $u$. This will be crucial to prove quantitative weak propagation of chaos using the Kolmogorov backward PDE, as explained in Chapter $5$ of \cite{CarmonaProbabilisticTheoryMean2018}.\\

 Let us denote by $(Z^n)_n$ an i.i.d.\ sequence of Lévy processes having the same distribution as $Z = (Z_t)_t,$ and by $(X^n_0)_n$ an i.i.d.\ sequence of random variables with common distribution $\mu_0 \in \PPP_{\beta}(\R^d).$ For a fixed integer $N \geq 1,$ the system of $N$ particles associated with \eqref{edsmckv} is defined as the solution to the following classical SDE on $(\R^{d})^N$ \begin{equation}\label{edsparticles}
\left\{  \begin{array}{lll}
&dX^{i,N}_t = b_t(X^{i,N}_t,\muu^N_t) \,dt + \sigma_t(X^{i,N}_{t^-},\muu^N_{t^-})\, dZ^i_t, \quad \forall i \leq N, \\ &\muu^N_t := \frac{1}{N} \displaystyle\sum_{j=1}^N \delta_{X^{j,N}_t},\\ &X^{i,N}_0 = X^i_0.
\end{array}\right.\end{equation}

We can write for all $i \leq N$ and for all $t \in [0,T]$ $$Z^i_t = \int_0^t \int_{B_1} z \, \NNN^i(ds,dz) + \int_0^t \int_{B_1^c} z \, \NN^i(ds,dz),$$ where $\NN^i$ is the Poisson random measure associated with $Z^i.$ Then, we set for all $t \in [0,T]$ $$\textbf{Z}^N_t = \begin{pmatrix} Z^1_t \\ \vdots \\ Z^N_t \end{pmatrix} \in (\R^d)^N.$$  
As the Lévy processes $(Z^n)_n$ are independent, the process $(\textbf{Z}^N_t)_t$ is a Lévy process in $(\R^d)^N.$  Its Poisson random measure $\bm{\NN}^N$ and its Lévy measure $\bm{\nu}^N$ are defined as follows. For all $\phi : [0,T] \times (\R^d)^N \rightarrow \R^+,$ one has \begin{align}\label{Poissonmeasureparticles}
	 \int_{0}^T\int_{(\R^d)^N} \phi(s,\x) \, \bm{\NN}^N(ds,d\x) \notag&= \sum_{i=1}^N \int_0^t \int_{\R^d} \phi(s,0, \dots,0,x_i,0,\dots,0) \, \NN^i(ds,dx_i) \\ &=  \sum_{i=1}^N \int_0^t \int_{\R^d} \phi(s,\bm{\tilde{x}_i}) \, \NN^i(ds,dx),\end{align} where  $\bm{\tilde{x}_i} := (0, \dots,x,\dots,0) \in (\R^d)^N$ for $x \in \R^d,$ where $x$ appears in the $i$-th position. For all $\phi : (\R^d)^N \rightarrow \R^+,$ one has \begin{equation}\label{Levymeasureparticles}
	 \int_{(\R^d)^N} \phi(\x) \, d\bm{\nu}^N(\x)= \sum_{i=1}^N  \int_{\R^d} \phi(\bm{\tilde{x}_i}) \, d\nu(x).\end{equation}  Note that since $(Z^n)_n$ are independent processes, for all $t \in [0,T],$ the support of the random measure $\bm{\NN}^N(t,d\x)$ is contained in $$ \bigcup_{i=0}^{N-1} \{0_{\R^d}\}^i \times \R^d \times \{0_{\R^d}\}^{N-1-i} \subset (\R^d)^N.$$

Let us define for all $t \in [0,T],$ $\x=(x_1, \dots,x_N) \in (\R^d)^N$ \begin{equation*}
\mathbf{b}^N_t(\x) := \begin{pmatrix} b_t(x_1,\muu^N_{\x}) \\ \vdots \\ b_t(x_N,\muu^N_{\x}) \end{pmatrix} \in (\R^d)^N \quad 	\bm{\sigma}^N_t(\x) := \begin{pmatrix} \sigma_t(x_1,\muu^N_{\x})& 0 & 0 \\0& \ddots&0 \\ 0 & 0 &\sigma_t(x_N,\muu^N_{\x}) \end{pmatrix} \in \R^{Nd \times Nd}.
\end{equation*}
 Thus, writing $\bm{X}^N_t = \begin{pmatrix}X^{1,N}_t \\ \vdots \\ X^{N,N}_t\end{pmatrix},$ SDE \eqref{edsparticles} defining the particle system can be rewritten as  \begin{equation}\label{edsparticulesvect}
\left\{  \begin{array}{lll}
&d\bm{X}^N_t = \mathbf{b}^N_t(\bm{X}^N_t) \,dt + \bm{\sigma}^N_t(\bm{X}^N_{t^-})\, d\textbf{Z}^N_t,  \\ &\bm{X}^N_0 = \begin{pmatrix} X^1_0 \\ \vdots \\ X^N_0 \end{pmatrix}.
\end{array}\right.\end{equation}

 Let us recall the definition of the stable operators $(L_t)_t$ associated with the particle system \eqref{edsparticles}. For any $f : (\R^d)^N \rightarrow \R$ of class $\CC^2$ and for any $\bm{x} = (x_1,\dots, x_N) \in (\R^d)^N,$ we define $L_tf$ by \begin{equation}\label{eq:stable_operator}
 L_tf (\bm{x}) := \bm{b}_t(\bm{x})\cdot \nabla f(\bm{x})+ \int_{(\R^d)^N} f(\bm{x} + \bm{\sigma}_t(\bm{x})\bm{z}) - f(\bm{x}) - \nabla f(\bm{x})\cdot (\bm{\sigma}_t(\bm{x})\bm{z}\1_{|\bm{z}|<1}) \, d\bm{\nu}(\bm{z}).\end{equation}

 \begin{Prop}\label{prop_approx_operator}
 	Let $ u : \PP \rightarrow \R$ be a function satisfying the assumptions of Theorem \ref{ito} and Lemma \ref{projempC2}. Assume moreover that $\de \del$ and $ \partial_{v'} \de \partial_v \del$ exist and are uniformly bounded on $\PP \times \R^d \times \R^d.$ Then, we have for all $t \in [0,T]$ and for all $\bm{x} = (x_1,\dots, x_N) \in (\R^d)^N$ \begin{align*}
 	L_t u^N(x_1,\dots,x_N)  &= \mathscr{L}_tu(\muu^N_{\bm{x}}) + \frac{1}{N}\sum_{k=1}^N\int_{\R^d} \int_{0}^1 \del (m^k_{t,z,w})(x_k+\sigma_t(x_k,\muu^N_{\bm{x}})z) -\del(
 	m^k_{t,z,w})(x_k) \\ &\hspace{4cm}- \left[ \del (\muu^N_{\bm{x}})(x_k+\sigma_t(x_k,\muu^N_{\bm{x}})z) - \del(\muu^N_{\bm{x}})(x_k)  \right]  \, dw \, d\nu(z) \\ &= \mathscr{L}_tu(\muu^N_{\bm{x}}) + O\left(\frac1N\right),
 	\end{align*}
 	where $m^k_{t,z,w} = w\muu^N_{\bm{x} + \bm{\sigma}_t(\bm{x},\muu^N_{\bm{x}})\bm{\tilde{z}_k}} + (1-w)\muu^N_{\bm{x}}$.
 \end{Prop}

 The last term is interpreted as an error term and we will have to control it to prove our propagation of chaos estimates in Section \ref{sectionOU}.

 \begin{proof}
 	We easily check that we can apply the operator $L_t$ to the empirical projection $u^N$ thanks to the regularity assumptions. It yields for all $\bm{x}=(x_1,\dots,x_N) \in (\R^d)^N$ \begin{align*}
 	L_t u^N(\bm{x}) &= \frac{1}{N} \sum_{k=1}^N b_t (x_k,\muu^N_{\bm{x}}) \partial_v \del (\muu^N_{\bm{x}})(x_k) \\ &\quad + \sum_{k=1}^N \int_{\R^d} u(\muu^N_{\bm{x} + \bm{\sigma}_t(\bm{x})\tilde{\bm{z_k}}}) - u(\muu^N_{\bm{x}}) - \frac{1}{N} \partial_v \del (\muu^N_{\bm{x}})(x_k) \cdot (\sigma_t(x_k,\muu^N_{\bm{x}})z \1_{|z| \leq 1}) \, d\nu(z)\\ &= \frac{1}{N} \sum_{k=1}^N b_t (x_k,\muu^N_{\bm{x}}) \partial_v \del (\muu^N_{\bm{x}})(x_k) \\ &\quad + \frac1N \sum_{k=1}^N \int_{\R^d} \int_0^1 \del (m^k_{t,z,w})(x_k + \sigma_t(x_k,\muu^N_{\bm{x}})z) - \del(m^k_{t,z,w})(x_k) \\ &\hspace{4cm}- \partial_v \del (\muu^N_{\bm{x}})(x_k) \cdot (\sigma_t(x_k,\muu^N_{\bm{x}})z \1_{|z| \leq 1}) \,dw \, d\nu(z) \\ &= \mathscr{L}_tu(\muu^N_{\bm{x}}) + \frac{1}{N}\sum_{k=1}^N\int_{\R^d} \int_{0}^1 \del (m^k_{t,z,w})(x_k+\sigma_t(x_k,\muu^N_{\bm{x}})z) -\del(
 	m^k_{t,z,w})(x_k) \\ &\hspace{4cm}- \left[ \del (\muu^N_{\bm{x}})(x_k+\sigma_t(x_k,\muu^N_{\bm{x}})z) - \del(\muu^N_{\bm{x}})(x_k)  \right]  \, dw \, d\nu(z).
 	\end{align*}
 	It remains to show the bound on the error term between $L_t u^N(\bm{x})$ and $\mathscr{L}_t u(\muu^N_{\bm{x}}).$ We split the domain of integration $\R^d$ of the last integral into $B_1$ and $B_1^c$ and we first consider the integral over $B_1^{c}.$ Defining $ m^k_{t,z,w,r} := r m^k_{t,z,w} + (1-r) \muu^N_{\bm{x}} =  \muu^N_{\bm{x}} + rw\muu^N_{\bm{x}+ \bm{\sigma}_t(\bm{x}) \bm{\tilde{z_k}}},$ we have
 	
 	\begin{align}\label{eq:error_not_well_defined}
 	&\notag\int_{B_1^c} \int_{0}^1 \del (m^k_{t,z,w})(x_k+\sigma_t(x_k,\muu^N_{\bm{x}})z) -\del(
 	m^k_{t,z,w})(x_k)\\ \notag &\hspace{4cm} - \left[ \del (\muu^N_{\bm{x}})(x_k+\sigma_t(x_k,\muu^N_{\bm{x}})z) - \del(\muu^N_{\bm{x}})(x_k)  \right]  \, dw \, d\nu(z)\\&= \frac{1}{N}\int_{B_1^c} \int_{[0,1]^2} w\left[\de\del(m^k_{t,z,w,r})(x_k+\sigma_t(x_k,\muu^N_{\bm{x}})z,x_k+\sigma_t(x_k,\muu^N_{\bm{x}})z)  \right. \\ \notag& \left. \hspace{3cm} -\de\del (m^k_{t,z,w,r})(x_k + \sigma_t(x_k,\muu^N_{\bm{x}})z,x_k) \right] \, dr\,dw\, d\nu(z)  \\\notag &\quad + \frac{1}{N}\int_{B_1^c} \int_{[0,1]^2} w\left[\de\del(m^k_{t,z,w,r})(x_k,x_k+\sigma_t(x_k,\muu^N_{\bm{x}})z)   -\de\del (m^k_{t,z,w,r})(x_k,x_k) \right] \, dr\,dw\, d\nu(z). 
 \end{align}
 	
 	Since $\de \del$ is uniformly bounded and $\nu (B_1^c) <+ \infty,$ we deduce that $$ \left\vert  \frac{1}{N}\int_{B_1^c} \int_{[0,1]^2} w\left[\de\del(m^k_{s,z,w,r})(x_k,x_k+\sigma_t(x_k,\muu^N_{\bm{x}})z)   -\de\del (m^k_{s,z,w,r})(x_k,x_k) \right] \, dr\,dw\, d\nu(z) \right\vert \leq \frac{C}{N}.$$
 	For the integral on $B_1,$ one has 
 	
 	\begin{align*}
 	&\int_{B_1} \int_{0}^1 \del (m^k_{t,z,w})(x_k+\sigma_t(x_k,\muu^N_{\bm{x}})z) -\del(
 	m^k_{t,z,w})(x_k)\\ &\hspace{4cm} - \left[ \del (\muu^N_{\bm{x}})(x_k+\sigma_t(x_k,\muu^N_{\bm{x}})z) - \del(\muu^N_{\bm{x}})(x_k)  \right]  \, dw \, d\nu(z)\\ &=  \int_{B_1} \int_{[0,1]^2} \left[\partial_v\del( m^k_{t,z,w})(x_k+h\sigma_t(x_k,\muu^N_{\bm{x}})z)   \partial_v\del( \muu^N_{\bm{x}})(x_k+h\sigma_t(x_k,\muu^N_{\bm{x}})z) \right]\cdot(\sigma_t(x_k,\muu^N_{\bm{x}})z) \, dh\,dw\, d\nu(z)  \\ &= \frac{1}{N}\int_{B_1} \int_{[0,1]^3} \left[\de\partial_v\del(m^k_{t,z,w,r})(x_k+h\sigma_t(x_k,\muu^N_{\bm{x}})z,x_k+\sigma_t(x_k,\muu^N_{\bm{x}})z)  \right. \\ & \hspace{4cm}- \left. \de \partial_v\del(m^k_{t,z,w,r})(x_k+h\sigma_t(x_k,\muu^N_{\bm{x}})z,x_k) \right]\cdot(\sigma_t(x_k,\muu^N_{\bm{x}})z)\, dr \, dh\,dw\, d\nu(z) \\ &=\frac{1}{N} \int_{B_1} \int_{[0,1]^4} \left[\partial_{v'}\de\partial_v\del(m^k_{t,z,w,r})(x_k+h\sigma_t(x_k,\muu^N_{\bm{x}})z,x_k+l\sigma_t(x_k,\muu^N_{\bm{x}})z) (\sigma_t(x_k,\muu^N_{\bm{x}})z) \right]  \\ & \hspace{10cm} \cdot(\sigma_t(x_k,\muu^N_{\bm{x}})z)\, dl\, dr \, dh\,dw\, d\nu(z).
 	\end{align*}
 	We conclude the proof using that $\sigma$ and $\partial_{v'}\de\partial_v\del$ are uniformly bounded and that $\int_{B_1} |z|^2 \, d\nu(z)$ is finite.
 	
 \end{proof}

\section{Propagation of chaos for a mean-field system of interacting stable Ornstein-Uhlenbeck processes}\label{sectionOU}

We now study a nonlinear stable Ornstein-Uhlenbeck process driven by an $\alpha$-stable process with $\alpha \in (1,2),$ written, for $t \geq 0$ $$ Z_t = \int_0^t \int_{\R^d} z \, \NNN(ds,dz).$$ We assume that $Z$ is non degenerate in the following sense. Writing $y = r\theta \in \R^d$ with $r\in \R^+$ and $\theta \in \S^{d-1},$ where $\S^{d-1}$ denotes the unit sphere in $\R^d,$  the Lévy measure $\nu$ of $Z$ decomposes as $$ \nu(dy) = d\mu(\theta) \frac{dr}{r^{1+\alpha}},$$ where $\mu$ is a non-zero finite measure on $\S^{d-1}.$  We assume that $\nu $ satisfies the following non-degeneracy assumption. \\

\noindent\textbf{(ND)} There exists $\eta >0$ such that for all $\lambda \in \R^d$ \begin{equation}\label{assumptionND}
	 \eta |\lambda|^{2} \leq \int_{\S^{d-1}} | \langle \lambda,\theta \rangle |^{2} \, d\mu(\theta).\end{equation}

Let $A,A',B \in \MM_d(\R)$ be matrices of size $d \times d$ such that $B$ is invertible, $\beta \in [1,\alpha)$, and $\xi \in L^{\beta}(\Omega, \FF_0)$ with distribution $ \mu_0 \in \PP.$ We study the following nonlinear SDE on the time interval $[0,T]$ \begin{equation}\label{MKVstableOUdef}
\left\{ \begin{array}{lll}
dX_t & = (AX_t + A' \E X_t)\,dt + B\,dZ_t, \\X_0  &= \xi.
\end{array} \right.
\end{equation}
 Note that weak uniqueness holds for SDE \eqref{MKVstableOUdef} i.e.\ the distribution of $(X_t)_t$ depends only on the distribution $\mu_0$ of $X_0$ and not on the random variable $X_0$ chosen.\\
 
  We consider the mean-field system of interacting Ornstein-Uhlenbeck processes associated with \eqref{MKVstableOUdef}. Let us denote by $(Z^n)_n$ an i.i.d.\ sequence of $\alpha$-stable processes having the same distribution as $Z = (Z_t)_t,$ and by $(X^n_0)_n$ an i.i.d.\ sequence of random variables with common distribution $\mu_0 \in \PPP_{\beta}(\R^d).$ For a fixed integer $N \geq 1,$ the system of $N$ particles associated with \eqref{edsmckv} is defined as the solution to the following classical SDE on $(\R^{d})^N$ \begin{equation}\label{edsparticlesOU}
\left\{  \begin{array}{lll}
&dX^{i,N}_t = AX^{i,N}_t\,dt + A' \displaystyle\frac{1}{N} \displaystyle\sum_{j=1}^N X^{j,N}_t \, dt +  B\,dZ^i_t,\quad \forall i \leq N,\\ &X^{i,N}_0 = X^i_0.
\end{array}\right.\end{equation}

We now state our quantitative weak propagation of chaos result in the next theorem. 

\begin{Thm}\label{thmpropchaosOU}
 Let $\mathcal{C}$ be the class of continuous functions $ u : \PP \rightarrow \R$ admitting two linear derivatives $\del$ and $\de \del$ such that $\del(\mu)$ and $\de \del (\mu)$ are Lipschitz continuous with Lipschitz constant smaller than $1$. Then, there exists a constant $C=C_T$ independent of $u \in \mathscr{C}$ and $N \geq 1$ and non-decreasing with respect to $T$ such that for all $u \in \mathscr{C}$ and $N \geq 1$, we have \begin{equation}\label{PropPCOUeq1}
	\E \left\vert u(\muu^N_T) - u(\mu_T) \right\vert \leq C\,\E W_1(\muu^N_0,\mu_0) +C\frac{\ln(N)^{\frac{1}{\alpha}}}{N^{1-\frac{1}{\alpha}}},
	\end{equation} 
	and 
	\begin{equation}\label{PropPCOUeq2}
	|\E (u(\muu^N_T) - u(\mu_T))| \leq  C\,\E W_1(\muu^N_0,\mu_0) + \frac{C}{N^{\alpha -1}}.
	\end{equation}
	
	Note that in the particular case where the initial distribution $\mu_0$ belongs to $\PPP_2(\R^d)$, we have
	\begin{equation}\label{PropPCOUeq2bis}
	|\E (u(\muu^N_T) - u(\mu_T))| \leq  \frac{C}{N^{\alpha -1}}.
	\end{equation}
\end{Thm}

\begin{Rq} \begin{itemize}
		
		\item The initial data terms can be handled using \cite{fournier:hal-00915365}, in particular in the case where $\mu_0$ has more moments than $\beta$. Indeed, one has if $\mu_0 \in \mathcal{P}_q(\R^d)$ \begin{equation*}
			\E W_1(\muu^N_0,\mu_0) \leq C \left\{\begin{array}{lll}
				&N^{-\frac12} + N^{-\left(1-\frac{1}{q}\right)}, \quad &\text{if}\quad d=1\quad \text{and}\quad q \neq 2,\\ 
				&N^{-\frac12}\ln(1+N) +N^{-\left(1-\frac{1}{q}\right)}, \quad &\text{if}\quad d=2 \quad \text{and}\quad q \neq 2, \\  &N^{-\frac1d} + N^{-\left(1-\frac{1}{q}\right)}, \quad &\text{if}\quad d\geq 3\quad \text{and} \quad q \neq \frac{d}{d -1}.
			\end{array}\right.
		\end{equation*}
		\item Let us denote by $\Vert \varphi \Vert_{\text{Lip}} := \sup_{x\neq y} \frac{|\varphi(x) -\varphi(y)|}{|x-y|}$ for $\varphi : \R^d \to \R$. The set $$\left\{ \phi : \PP \to \R, \, \exists \varphi : \R^d \to \R, \, \text{with}\, \Vert\varphi \Vert_{\text{Lip}}  \leq 1,\, \text{and}\, \phi(\mu) = \int_{\R^d} \varphi \,d\mu, \, \forall \mu \in \PP \right\}$$ is contained in $\CC$. This allows quantify the convergence with respect to $W_1$. More precisely, for the mean-field limit, we have  by the Kantorovich-Rubinstein theorem and if $\mu_0 \in \mathcal{P}_{2}(\R^d)$ \begin{align*} 
			\sup_{t \in [0,T]}W_1([X^{1,N}_t],\mu_t) & =\sup_{t \in [0,T]}\sup_{\varphi,\, \|\varphi\|_{\text{Lip}}\leq1} \left\vert 
			\E \varphi(X^{1,N}_t) -  \int_{\R^d} \varphi \, d\mu_t \right\vert \\ &= \sup_{t \in [0,T]}\sup_{\varphi,\, \|\varphi\|_{\text{Lip}}\leq1} \left\vert 
			\E \left(\frac1N \sum_{k=1}^N \varphi(X^{k,N}_t)\right) -  \int_{\R^d} \varphi \, d\mu_t \right\vert\\ &=
			\sup_{t \in [0,T]}\sup_{\varphi,\, \|\varphi\|_{\text{Lip}}\leq 1} \left\vert 
			\E \int_{\R^d} \varphi \, d\muu^N_t - \E \int_{\R^d} \varphi \, d\mu_t \right\vert \\ &=
			\sup_{t \in [0,T]}\sup_{\varphi,\, \|\varphi\|_{\text{Lip}}\leq 1} \left\vert 
			\E \phi(\muu^N_t) - \E \phi (\mu_t) \right\vert\\&\leq \frac{C_T}{N^{\alpha -1}},
		\end{align*}
		
	 since the constant $C$ in Theorem \ref{thmpropchaosOU} is non-decreasing with respect to $T$.

	\end{itemize}
\end{Rq}

\begin{Rq}\label{rqbigjumps}
 In order to obtain our estimates of propagation of chaos, we follow the method mentioned in Chapter $5$ of \cite{CarmonaProbabilisticTheoryMean2018} (pages $506-508$). We need to apply the operator $L_t$ associated with the particle system, defined in \eqref{eq:stable_operator}, for the empirical projection $\bm{x} \in (\R^d)^N \mapsto \phi_u(t,\muu^N_{\bm{x}})$ of $\phi_u$, which was defined in Definition \ref{defsemigroup}. In order to use the backward Kolmogorov PDE satisfied by $\phi_u$, we need to control the difference between $L_t\phi_u(t,\muu^N_{\bm{x}})$ and $\mathscr{L}_t \phi_u(t,\muu^N_{\bm{x}}),$ as done in Proposition \ref{prop_approx_operator}. However, the assumptions on $\de \dell$ are not satisfied. Indeed, for this interacting Ornstein-Uhlenbeck system, we can only show (see Proposition \ref{ThmPCremovjumps}) that for any fixed $t \in (0,T],$ $\mu \in \PP,$ and $v,v' \in \R^d$ \begin{equation}
	\left\vert \de \dell (t,\mu)(v,v') \right\vert \leq C( 1 + |v| + |v'| + |v||v'|).	\end{equation}This is due to the unboundedness of drift which is only at most of linear growth with respect to $x$ and $\mu$. The error term in \eqref{eq:error_not_well_defined} is not directly well-defined since $\int_{\R^d} |z|^2 \, d\nu(z) = + \infty.$ To get around this difficulty, we will consider in a first step truncated versions of the operators $L_t$ and $\mathscr{L}_t$ by removing the big jumps of the noise. Namely we replace $\nu$ by the restriction of $\nu$ on the ball $B_N,$ where $N$ is also the number of particles interacting. We refer to \eqref{ThmPC3proofeq3}, \eqref{thmPC3prooff} and \eqref{ThmPC3proofeq4} in the proof of Theorem \ref{thmpropchaosOU}, where the terms are not well-defined without removing the big jumps.
	
	\end{Rq}

 We thus set for $t \in [0,T]$ \begin{equation}\label{deftroncatednoise}Z_{N,t} := \int_0^t \int_{B_N} z \, \NNN(ds,dz).\end{equation} 
 Then, we consider the McKean-Vlasov SDE \eqref{MKVstableOUdef} driven by $Z_N.$ The solution to this SDE is denoted by $(X_{N,t})_t$ and satisfies 

 \begin{equation}\label{edsmckvtronquee1}
\left\{ \begin{array}{lll}
dX_{N,t} & = (AX_{N,t} + A' \E X_{N,t})\,dt + B\,dZ_{N,t}, \\X_0  &= \xi \in L^{\beta}(\Omega,\FF_0).
\end{array} \right.
\end{equation}

We denote by $\mu_{N,t}$ the distribution of $X_{N,t}$ and by $\phi_{N,u}$ the solution to the backward Kolmogorov PDE associated with the McKean-Vlasov SDE \eqref{edsmckvtronquee1} defined in Definition \ref{defsemigroup}.\\

Before proving the preceding theorem, we state in the following proposition regularity properties and bounds on $\phi_{N,u}$ that we will use.

\begin{Prop}\label{ThmPCremovjumps}
 For all $u \in \mathscr{C}$ and $N \geq 1,$ the function $\phi_{N,u}$ satisfies the following properties. \begin{enumerate}
		\item The function $\phi_{N,u}$ belongs to $\CC^0([0,T]\times \PPP_{\beta}(\R^d);\R),$  $\partial_t \phi_{N,u} $ exists and is continuous on $(0,T] \times \PPP_{\beta}(\R^d)$. Moreover $\dellN,$ $\partial_v \dellN$ and $\partial^2_v \dellN$ exist and are continuous on $(0,T]  \times \PPP_{\beta}(\R^d) \times \R^d.$
		
		\item The functions $ \de \partial_v \dellN $ and $\partial_{v'} \de \partial_v \dellN$ exist and are continuous on $(0,T] \times \PPP_{\beta}(\R^d) \times \R^{d} \times \R^d.$
		\item There exists three continuous functions on $(0,T]$ denoted by $g_1$, $g_2$ and $g_3$ such that $g_1$ is globally bounded on $[0,T]$ and $g_2 \in  L^1(0,T)$ which satisfy the following properties. For all $u \in \mathcal{C}$, $N\geq 1,$ $t \in (0,T],$ $\mu \in \PPP_{\beta}(\R^d)$ and $v,v' \in \R^d$, one has
		\begin{align}\label{ThmPC3eq1}
		\notag& \left\vert \partial_v \dellN (t,\mu)(v) \right\vert \leq g_1(t), \\ & \left\vert \partial_{v'} \de \partial_v\dellN (t,\mu)(v,v') \right\vert \leq g_2(t), \\ \notag& \left\vert \partial^2_v \dellN (t,\mu)(v) \right\vert \leq g_3(t).
		\end{align}

	\end{enumerate}
	
\end{Prop}

The proof is postponed to the Appendix (Section \ref{proofpropositionreg}).\\

We can now prove the propagation of chaos estimate for the interacting Ornstein-Uhlenbeck processes.

\begin{proof}[Proof of Theorem \ref{thmpropchaosOU}]
	
		We assume that the matrix $B$ equal to the identity $I_d$ since it does not change anything to the proof. Let us fix $N \geq 1$. As we did in \eqref{deftroncatednoise} for $Z=(Z_t)_t$, we remove the big jumps of the i.i.d.\ copies $(Z^i)_i$ of $Z$ by defining for all $i \in \N$ and $t \in [0,T] $ $$ Z^i_{N,t} := \int_0^t \int_{B_N} z \, \NNN^i(ds,dz).$$ For the sake of clarity, we set for $t \in [0,T],$ $x \in \R^d$ and $\mu \in \PPP_1(\R^d)$ $$ b(x,\mu):= Ax + A' \int_{\R^d} x \, d\mu(x).$$ Then, we consider $(X_{N,t})_t$ the nonlinear Ornstein-Uhlenbeck process driven by $Z_{N}$ together with $\bm{X}^N_N:=(X^{1,N}_{N,t}, \dots, X^{N,N}_{N,t})_t$ the associated interacting particle system, i.e.\ the solutions to \begin{equation}\label{edsmckvtronquee}
	\left\{  \begin{array}{lll}
	&dX_{N,t} = b(X_{N,t},\mu_{N,t}) \,dt + dZ_{N,t}, \\ &\mu_{N,t} = \mathcal{L}(X_{N,t}),\\ &X_{N,0} = \xi \in L^{\beta}(\Omega,\FF_0),
	\end{array}\right. 
	\end{equation} and \begin{equation}\label{edsparticlestronquee}
	\left\{  \begin{array}{lll}
	&dX^{i,N}_{N,t}= b(X^{i,N}_{N,t},\muu^N_{N,t}) \,dt +  dZ^i_{N,t}, \quad \forall t \in [0,T],\, \forall i \leq N, \\ &\muu^N_{N,t}= \displaystyle\frac{1}{N} \displaystyle\sum_{j=1}^N \delta_{X^{j,N}_{N,t}},\\ &X^{i,N}_{N,0} = X^i_0.
	\end{array}\right.\end{equation}
	Let us emphasize that we remove the jumps that are bigger than the number of particles $N.$ Recall that the operator $\mathscr{L}_{N}$ associated with \eqref{edsmckvtronquee} defined in Theorem \ref{ThmEDP} is given, for all $s \in [0,T]$ and $\mu \in \PPP_{\beta}(\R^d)$, by \begin{align}\label{operatormeasuretronquee}
	\mathscr{L}_{N} \phi_{N,u}(s, \mu) & := \int_{\R^d} \partial_v \dellN (s,\mu)(v) \, b(v,\mu) \, d\mu(v) \\ \notag&\quad + \int_{\R^d} \int_{B_N} \left[\dellN (s ,\mu)(v+z) - \dellN (s ,\mu)(v) -  \partial_v \dellN (s,\mu)(v)\cdot z\right]  \, d\nu(z) \, d\mu(v).
	\end{align}
In order to establish \eqref{PropPCOUeq1}, our aim is to control $\E|\phi_{N,u}(T-t,\muu^N_t) - \phi_u(T-t,\mu_t)|,$ for $t <T$ uniformly in $N.$ We decompose it in the following way \begin{align*}
	&\notag\E|\phi_{N,u}(T-t,\muu^N_t) - \phi_{N,u}(T-t,\mu_t)| \\ &\leq \E|\phi_{N,u}(T-t,\muu^N_t) - \phi_{N,u}(T-t,\muu^N_{N,t})| + \E|\phi_{N,u}(T-t,\muu^N_{N,t}) - \phi_{N,u}(T-t,\mu_{N,t})|\\ &\quad + \E|\phi_{N,u}(T-t,\mu_{N,t}) - \phi_{N,u}(T-t,\mu_t)|.
	\end{align*}
	Since $\partial_v \dellN$ is bounded on $[0,T] \times \PP \times \R^d$ uniformly in $N$ by Proposition \ref{ThmPCremovjumps}, we deduce by the Kantorovich-Rubinstein theorem that for some constant $C>0$, we have for all $N\geq 1,$ $t \in [0,T),$ $ \mu,\nu \in \PP$ $$ |\phi_{N,u}(T-t,\mu) - \phi_{N,u}(T-t,\nu)| \leq C W_1(\mu,\nu).$$ Thus, we obtain \begin{align}\label{ThmPC3proof1}
	\notag\E|\phi_{N,u}(T-t,\muu^N_t) - \phi_{N,u}(T-t,\mu_t)| &\leq C\, \E W_1(\muu^N_t,\muu^N_{N,t}) + \E|\phi_{N,u}(T-t,\muu^N_{N,t}) - \phi_{N,u}(T-t,\mu_{N,t})| \\ &\quad + C\,\E W_1(\mu_{N,t},\mu_t) \\ \notag&=: A_1 + A_2 + A_3.
	\end{align}
	First, we deal with $A_1.$ Note that we can write $$ X^{k,N}_{N,t} - X^{k,N}_t = \int_0^t (b(X^{k,N}_{N,s}, \muu^N_{N,s}) - b(X^{k,N}_s,\muu^N_s)\, ds + \int_0^t \int_{B_N^c} z \, \NNN^k(ds,dz).$$ Using the fact that $b$ is Lipschitz continuous on $\R^d \times \PPP_1(\R^d)$, one has for all $t \in [0,T]$\begin{align*}
	\E  |X^{k,N}_{N,t} - X^{k,N}_t | &\leq C \left[ \int_0^t\E|X^{k,N}_{N,s} - X^{k,N}_s |\, ds + \int_0^t \frac1N \sum_{j=1}^N \E|X^{j,N}_{N,s} - X^{j,N}_s |\, ds  + \int_0^t\int_{B_N^c} |z| \, d\nu(z)\,ds\right] \\ &\leq  C \left[\int_0^t\E|X^{k,N}_{N,s} - X^{k,N}_s |\, ds + \int_0^t \frac1N \sum_{j=1}^N \E|X^{j,N}_{N,s} - X^{j,N}_s |\, ds + \frac{1}{N^{\alpha -1}} \right].
	\end{align*}
	
	Gronwall's inequality ensures that there exists a constant $C=C_T$ depending only on $T$ and such that for any $t \in [0,T]$ \begin{equation}\label{ThmPC3proofeqi}
	\E W_1(\muu^N_t,\muu^N_{N,t}) \leq \frac1N \sum_{k=1}^N \E|X^{k,N}_{N,t} - X^{k,N}_t | \leq \frac{C}{N^{\alpha-1}}
	\end{equation}
	We follow the same lines of reasoning to treat $A_3.$ Indeed, writing for all $t \in [0,T]$ $$ X_{N,t} - X_{t} = \int_0^t (b(X_{N,s}, \mu_{N,s}) - b(X_{s},\mu_s)\, ds + \int_0^t \int_{B_N^c} z \, \NNN(ds,dz),$$ we deduce as previously that there exists $C=C_T>0$ such that for any $t \in [0,T]$ \begin{equation}\label{ThmPC3proofeqiii}
	\E W_1(\mu_{N,t}, \mu_t) \leq \frac{C}{N^{\alpha-1}}.
	\end{equation}
	It remains to treat the term $A_2$ in \eqref{ThmPC3proof1}. Using Lemma \ref{projempC2} and the regularity properties satisfied by $\phi_{N,u}$ stated in Proposition \ref{ThmPCremovjumps}, we obtain that the function $(t,\bm{x}) \in [0,T) \times (\R^d)^N\mapsto \phi_{N,u}(T-t,\muu^N_{\bm{x}})$ belongs to $\CC^{1,2}([0,T) \times (\R^d)^N).$ Moreover, for all $t\in [0,T)$ and $\bm{x} = (x_1, \dots, x_N)\in (\R^d)^N,$ we have $$ \partial_{\bm{x}} \phi_{N,u}(T-t,\muu^N_{\bm{x}}) = \frac1N\begin{pmatrix}
	\partial_v \dellN (T-t,\muu^N_{\bm{x}})(x_1) \\ \vdots \\ \partial_v \dellN (T-t,\muu^N_{\bm{x}})(x_N)
	\end{pmatrix}.$$ 
	
We denote by $\bm{\NN}^N$ the Poisson random measure  associated with $\bm{Z}^N = (Z^1, \dots, Z^N)$ defined in \eqref{Poissonmeasureparticles} and by $\bm{\nu}^N$ its associated Lévy measure defined in \eqref{Levymeasureparticles}. Applying the standard Itô formula for this function and the $(\R^d)^N$-valued process $(\bm{X}^N_{N,t})_t$ and noticing that the map $t\in[0,T] \mapsto \phi_{N,u}(T-t,\mu_{N,t})$ is constant, we obtain that for all $t \in [0,T)$
	
	\begin{align}\label{thmPC2proof1'}
	\notag&\phi_{N,u}(T-t,\muu^N_{N,t}) - \phi_{N,u}(T-t,\mu_{N,t}) - \left(\phi_{N,u}(T,\muu^N_0) - \phi_{N,u}(T,\mu_0)\right) \\\notag&=  - \int_0^t \partial_t \phi_{N,u}(T-s,\muu^N_{N,s}) \, ds \\ \notag&\quad + \frac1N \sum_{i=1}^N \int_0^t \partial_v \dellN (T-s,\muu^N_{N,s})(X^{i,N}_{N,s})\cdot b(X^{i,N}_{N,s}, \muu^N_{N,s}) \, ds \\ \notag&\quad + \int_0^t \int_{(B_N)^N} \left[\phi_{N,u}(T-s,\muu^N_{\bm{X}^N_{N,s^-}+ \bm{z}}) - \phi_{N,u}(T-s,\muu^N_{\bm{X}^N_{N,s^-}}) - \partial_{\bm{x}}\phi_{N,u}(T-s,\muu^N_{\bm{X}^N_{N,s^-}})\cdot \bm{z}\right] \, d\bm{\nu}^N(\bm{z}) \, ds\\ &\quad + \int_0^t \int_{(B_1)^N} \left[\phi_{N,u}(T-s,\muu^N_{\bm{X}^N_{N,s^-}+ \bm{z}}) - \phi_{N,u}(T-s,\muu^N_{\bm{X}^N_{N,s^-}})\right] \, \bm{\NNN}^N(ds,d\bm{z}) \\ \notag&\quad + \int_0^t \int_{(B_N \backslash B_1)^N} \left[\phi_{N,u}(T-s,\muu^N_{\bm{X}^N_{N,s^-}+ \bm{z}}) - \phi_{N,u}(T-s,\muu^N_{\bm{X}^N_{N,s^-}}) \right]\, \bm{\NNN}^N(ds,d\bm{z}) \\ \notag & =: I_1 + I_2 + I_3 + I_4 + I_5.\end{align}
	
	Note that $$\int_0^t \int_{(B_N \backslash B_1)^N} \left[\phi_{N,u}(T-s,\muu^N_{\bm{X}^N_{N,s^-}+ \bm{z}}) - \phi_{N,u}(T-s,\muu^N_{\bm{X}^N_{N,s^-}}) \right]\, d\bm{\nu}^N(\bm{z}) \, ds$$ is well-defined. Indeed, recalling that for $h \in \R^d$, $\bm{\tilde{h_j}} = ( 0, \dots, 0,h,0, \dots, 0) \in (\R^d)^N,$ where $h$ appears in the $j$-th coordinate, one has  \begin{align*}
	&\int_0^t \int_{(B_N \backslash B_1)^N} \left\vert\phi_{N,u}(T-s,\muu^N_{\bm{X}^N_{N,s^-}+ \bm{z}}) - \phi_{N,u}(T-s,\muu^N_{\bm{X}^N_{N,s^-}}) \right\vert\, d\bm{\nu}^N(\bm{z}) \, ds \\ &= \sum_{i=1}^N \int_0^t \int_{B_N \backslash B_1} \left\vert\phi_{N,u}(T-s,\muu^N_{\bm{X}^N_{N,s^-}+ \bm{\tilde{z_i}}}) - \phi_{N,u}(T-s,\muu^N_{\bm{X}^N_{N,s^-}}) \right\vert\, d\nu(z) \, ds \\& \leq \frac{1}{N} \sum_{i=1}^N \int_0^t \int_{B_N \backslash B_1} \int_0^1\left\vert \partial_v\dellN(T-s,m^i_{s,z,w})(X^{i,N}_{N,s^-} + h z) \right\vert \vert z \vert\,dw \, d\nu(z) \, ds \\ &\leq C \int_0^t g_1(T-s) \, ds \int_{B_N \backslash B_1} |z| \, d\nu(z),
	\end{align*}
		where $m^i_{s,z,w} := w \muu^N_{\bm{X}^N_{N,s^-}+ \bm{\tilde{z_i}}} + (1-w)\muu^N_{\bm{X}^N_{N,s^-}}$. We conclude since $g_1$ is bounded by Proposition \ref{ThmPCremovjumps}.\\

	 Using the form of the Lévy measure of the particle system in \eqref{Levymeasureparticles}, we can write \begin{align*}
I_3 &= \sum_{i=1}^N \int_0^t \int_{(B_N)^N} \left[\phi_{N,u}(T-s,\muu^N_{\bm{X}^N_{N,s^-}+ \bm{\tilde{z_i}}}) - \phi_{N,u}(T-s,\muu^N_{\bm{X}^N_{N,s^-}})  \right. \\ & \left.\hspace{5cm}-\frac1N  \partial_v \dellN (T-s,\muu^N_{\bm{X}^N_{N,s^-}})(X^{i,N}_{N,s^-})\cdot z\right]\, d\nu(z) \, ds \\ &= \frac1N\sum_{i=1}^N \int_0^t \int_{B_N} \int_0^1 \left[\dellN(T-s, m^i_{s,z,w})(X^{i,N}_{N,s^{-}} + z)\right. \\ & \left.\hspace{3cm}- \dellN(T-s,m^i_{s,z,w})(X^{i,N}_{N,s^{-}})  -\frac1N  \partial_v \dellN (T-s,\muu^N_{\bm{X}^N_{N,s^-}})(X^{i,N}_{N,s^-})\cdot z\right]\, dw\, d\nu(z) \, ds,
\end{align*}
where $m^i_{s,z,w} = w \muu^N_{\bm{X}^N_{N,s^-}+ \bm{\tilde{z_i}}} + (1-w)\muu^N_{\bm{X}^N_{N,s^-}}.$ In order to make appear the backward Kolmogorov PDE \eqref{PDE}, we decompose $I_3$ in the following way

\begin{align*}
I_3 &=  \int_0^t \int_{\R^d} \int_{B_N} \left[\dellN(T-s, \muu^N_{\bm{X}^N_{N,s^-}})(x + z) - \dellN(T-s,\muu^N_{\bm{X}^N_{N,s^-}})(x) \right. \\ & \hspace{4cm} \left.-  \partial_v \dellN (T-s,\muu^N_{\bm{X}^N_{N,s^-}})(x)\cdot z \right]\, d\nu(z) \,d\muu^N_{\bm{X}^N_{N,s^-}}(x) \, ds \\ &\quad + \frac1N\sum_{i=1}^N \int_0^t \int_{B_N} \int_0^1 \left[\dellN(T-s, m^i_{s,z,w})(X^{i,N}_{N,s^{-}} + z) -  \dellN(T-s,m^i_{s,z,w})(X^{i,N}_{N,s^{-}}) \right. \\ & \hspace{2cm}+\left. \dellN(T-s,\muu^N_{\bm{X}^N_{N,s^-}})(X^{i,N}_{N,s^{-}}) - \dellN(T-s, \muu^N_{\bm{X}^N_{N,s^-}})(X^{i,N}_{N,s^{-}} + z) \right]\,dw\, d\nu(z) \, ds \\ &=: I_{3,A} + I_{3,B}.
\end{align*}

Since $(\bm{X}^N_{N,s})_{s \in [0,T]}$ is càd-làg, we deduce that almost surely for almost all $s \in [0,t]$ $$ \muu^N_{\bm{X}^N_{N,s^{-}}} = \muu^N_{\bm{X}^N_{N,s}} = \muu^N_{N,s}.$$ Thanks to the backward Kolmogorov PDE \eqref{PDE} in Theorem \ref{ThmEDP} applied with $\beta \in [1, \alpha)$ and $\gamma =1$ and justified by Proposition \ref{ThmPCremovjumps}, one has \begin{align*}
I_1 + I_2 + I_{3,A} &=  \int_0^t -\partial_t \phi_{N,u}(T-s,\muu^N_{N,s}) + \mathscr{L}_N \phi_{N,u}(T-s,\muu^N_{N,s}) \, ds \\ &=0.
\end{align*}

Thus, we obtain the following decomposition  \begin{align}\label{ThmPC3proofeq1}
\notag&\phi_{N,u}(T-t,\muu^N_{N,t}) - \phi_{N,u}(T-t,\mu_{N,t}) - \left(\phi_{N,u}(T,\muu^N_{0}) - \phi_{N,u}(T,\mu_{0})\right) \\ \notag&=  \frac1N\sum_{i=1}^N \int_0^t \int_{B_N} \int_0^1 \left[\dellN(T-s, m^i_{s,z,w})(X^{i,N}_{N,s^{-}} + z) -  \dellN(T-s,m^i_{s,z,w})(X^{i,N}_{N,s^{-}}) \right. \\ \notag& \hspace{2cm}+\left. \dellN(T-s,\muu^N_{\bm{X}^N_{N,s^-}})(X^{i,N}_{N,s^{-}}) - \dellN(T-s, \muu^N_{\bm{X}^N_{N,s^-}})(X^{i,N}_{N,s^{-}} + z) \right]\,dw\, d\nu(z) \, ds  \\ &\quad + \int_0^t \int_{B_N\backslash B_1} \left[\phi_{N,u}(T-s,\muu^N_{\bm{X}^N_{N,s^-}+ \bm{z}}) - \phi_{N,u}(T-s,\muu^N_{\bm{X}^N_{N,s^-}}) \right]\, \bm{\NNN}^N(ds,d\bm{z}) \\ \notag&\quad + \int_0^t \int_{B_1} \left[\phi_{N,u}(T-s,\muu^N_{\bm{X}^N_{N,s^-}+ \bm{z}}) - \phi_{N,u}(T-s,\muu^N_{\bm{X}^N_{N,s^-}})\right] \, \bm{\NNN}^N(ds,d\bm{z}) \\\notag &=  I_{3,B}  + I_4 + I_5.\end{align}

It follows that for all $t \in [0,T)$ \begin{equation}\label{thmPC2proof3}
\E|\phi_{N,u}(T-t,\muu^N_{N,t}) - \phi_{N,u}(T-t,\mu_{N,t})| \leq   \E |\phi_{N,u}(T,\muu^N_0) - \phi_{N,u}(T,\mu_0)| +  \E (|I_{3,B}| + |I_4| + |I_5|).
\end{equation}
We treat each term separately. For $I_{3,B},$ we write 
	
	\begin{align}\label{ThmPC3proofeq3}
	\notag&I_{3,B}\\\notag &= \frac1N\sum_{i=1}^N \int_0^t \int_{B_N} \int_{[0,1]^2} \left[\partial_v\dellN(T-s, m^i_{s,z,w})(X^{i,N}_{N,s^{-}} + hz)  \right. \\\notag & \hspace{4cm}-\left. \partial_v\dellN(T-s, \muu^N_{\bm{X}^N_{N,s^-}})(X^{i,N}_{N,s^{-}} + hz) \right]\cdot z \, dh\,dw\, d\nu(z) \, ds \\ &= \frac{1}{N^2}\sum_{i=1}^N \int_0^t \int_{B_N} \int_{[0,1]^3} \left[\de\partial_v\dellN(T-s, m^i_{s,z,w,r})(X^{i,N}_{N,s^{-}} + hz,X^{i,N}_{N,s^{-}} + z)  \right. \\ \notag& \hspace{4cm}\quad-\left. \de \partial_v\dellN(T-s, m^i_{s,z,w,r})(X^{i,N}_{N,s^{-}} + hz,X^{i,N}_{N,s^{-}}) \right]\cdot z\, dr \, dh\,dw\, d\nu(z) \, ds\\ \notag&=\frac{1}{N^2}\sum_{i=1}^N \int_0^t \int_{B_N} \int_{[0,1]^4} \left[\partial_{v'}\de\partial_v\dellN(T-s, m^i_{s,z,w,r})(X^{i,N}_{N,s^{-}} + hz,X^{i,N}_{N,s^{-}} + kz)  z \right]\\ &\notag\hspace{12cm}\cdot z\, dk\, dr \, dh\,dw\, d\nu(z) \, ds,
	\end{align}
	
	where $ m^i_{s,z,w,r} := r m^i_{s,z,w} + (1-r) \muu^N_{\bm{X}^N_{N,s^-}} =  \muu^N_{\bm{X}^N_{N,s^-}} + rw\muu^N_{\bm{X}^N_{N,s^-} +  \bm{\tilde{z_i}}}.$ We deduce that there exists a constant $C>0$ independent of $N$ such that  for all $t \in (0,T]$ \begin{align}\label{thmPC3prooff}
	\notag\E|I_{3,B}| &\leq \frac{1}{N^2}\sum_{i=1}^N \E \int_0^t \int_{B_N}  g_2(T-s)|z|^2 \, d\nu(z) \, ds  \\ &\leq \frac{C}{N^{\alpha -1}}.
	\end{align}

	Let us recall that $m^i_{s,z,w} = w \muu^N_{\bm{X}^N_{N,s^-}+ \bm{\tilde{z_i}}} + (1-w)\muu^N_{\bm{X}^N_{N,s^-}}.$ It follows from the Cauchy-Schwarz inequality and the $L^2$-isometry of compensated Poisson random integrals that there exists a constant $C>0$ independent of $N$ such that we have for all $t \in (0,T]$ \begin{align}\label{thmPC3proof6}
	\notag\E |I_4|  & \leq  \left(\E \left(\int_0^t \int_{(B_1)^N} \left\vert\phi_{N,u}(T-s,\muu^N_{\bm{X}^N_{N,s^-}+ \bm{z}}) - \phi_{N,u}(T-s,\muu^N_{\bm{X}^N_{N,s^-}}) \right]\, \bm{\NNN}^N(ds,d\bm{z}) \right)^2 \right)^{\frac12}\\ &= \left( \sum_{i=1}^N\E \int_0^t \int_{B_1} \left\vert\phi_{N,u}(T-s,\muu^N_{\bm{X}^N_{N,s^-}+ \bm{\tilde{z_i}}}) - \phi_{N,u}(T-s,\muu^N_{\bm{X}^N_{N,s^-}}) \right\vert^2\, d\nu(z) \, ds \right)^{\frac12} \\ \notag&= \left( \sum_{i=1}^N\E \int_0^t \int_{B_1} \left\vert \frac1N \int_{[0,1]^2}\partial_v\dellN(T-s,m^i_{s,z,w})(X^{i,N}_{N,s^-}+ hz)\cdot z\,  dh \,dw \right\vert^2\, d\nu(z) \, ds \right)^{\frac12} \\ \notag&\leq  \left( \sum_{i=1}^N\E \int_0^t \int_{B_1} \frac{1}{N^2}\left\vert  g_1(T-s) \right\vert^2 C^2 |z|^2\, d\nu(z) \, ds \right)^{\frac12}\\ \notag&\leq \frac{C}{\sqrt{N}}.
	\end{align}
	
	Finally, for $I_5,$ BDG's inequalities and the fact that $1<\alpha< 2$ yields for all $t \in [0,T)$ \begin{align}\label{ThmPC3proofeq4}
	\notag\E |I_5| & \leq  \left(\E \left(\int_0^t \int_{(B_N \backslash B_1)^N} \left\vert\phi_{N,u}(T-s,\muu^N_{\bm{X}^N_{N,s^-}+ \bm{z}}) - \phi_{N,u}(T-s,\muu^N_{\bm{X}^N_{N,s^-}}) \right]\, \bm{\NNN}^N(ds,d\bm{z}) \right)^\alpha \right)^{\frac{1}{\alpha}}\\\notag & \leq C\left(\E \left[ \int_0^t \int_{(B_N \backslash B_1)^N} \left\vert\phi_{N,u}(T-s,\muu^N_{\bm{X}^N_{N,s^-}+ \bm{z}}) - \phi_{N,u}(T-s,\muu^N_{\bm{X}^N_{N,s^-}}) \right\vert^2\, \bm{\NN}^N(ds,d\bm{z})\right]^{\frac{\alpha}{2}} \right)^{\frac{1}{\alpha}} \\& \leq C\left(\E \int_0^t \int_{(B_N \backslash B_1)^N} \left\vert\phi_{N,u}(T-s,\muu^N_{\bm{X}^N_{N,s^-}+ \bm{z}}) - \phi_{N,u}(T-s,\muu^N_{\bm{X}^N_{N,s^-}}) \right\vert^\alpha\, \bm{\NN}^N(ds,d\bm{z}) \right)^{\frac{1}{\alpha}} \\ \notag&= C\left( \sum_{i=1}^N\E \int_0^t \int_{B_N \backslash B_1} \left\vert \frac1N \int_{[0,1]^2}\partial_v\dellN(T-s,m^i_{s,z,w})(X^{i,N}_{N,s^-}+ hz)\cdot z \,dh \, dw \right\vert^\alpha\, d\nu(z) \, ds \right)^{\frac1\alpha} \\ \notag&\leq  C\left( \sum_{i=1}^N\E \int_0^t \int_{B_N \backslash B_1} \left\vert \frac1N g_1(T-s) \right\vert^\alpha |z|^\alpha\, d\nu(z) \, ds \right)^{\frac{1}{\alpha}}\\ \notag&\leq C\frac{\ln(N)^{\frac{1}{\alpha}}}{N^{1-\frac{1}{\alpha}}}.
	\end{align}

	As a consequence of \eqref{ThmPC3proofeq1}, \eqref{thmPC3prooff}, \eqref{thmPC3proof6},  \eqref{ThmPC3proofeq4}, and the fact that $(\alpha -1) \wedge \frac12 > 1 -\frac{1}{\alpha},$ we have for a constant $C>0$ that for all $ N\geq 1$ and $t \in [0,T)$ \begin{equation}\label{ThmPC3proofeq5}
	\E|\phi_{N,u}(T-t,\muu^N_{N,t}) - \phi_{N,u}(T-t,\mu_{N,t})| \leq   \E |\phi_{N,u}(T,\muu^N_0) - \phi_{N,u}(T,\mu_0)| +C\frac{\ln(N)^{\frac{1}{\alpha}}}{N^{1-\frac{1}{\alpha}}}.
	\end{equation}
	Note that our assumptions ensure that the function $\phi_{N,u}(T, \cdot)$ is Lipschitz continuous with respect to the Wasserstein metric $W_1$ uniformly with respect to $N$. It is a consequence of the Kantorovich-Rubinstein theorem since $\dellN (T, \mu)$ is Lipschitz continuous uniformly with respect to $\mu$ and $N$ by Proposition \ref{ThmPCremovjumps}. We thus obtain for all $ N\geq 1,$ $t \in [0,T)$ \begin{equation}\label{ThmPC3proofeq6}
	\E|\phi_{N,u}(T-t,\muu^N_{N,t}) - \phi_{N,u}(T-t,\mu_{N,t})| \leq   C\, \E W^1(\muu^N_0,\mu_0) +C\frac{\ln(N)^{\frac{1}{\alpha}}}{N^{1-\frac{1}{\alpha}}}.
	\end{equation} Using \eqref{ThmPC3proof1}, \eqref{ThmPC3proofeqi}, \eqref{ThmPC3proofeqiii} and \eqref{ThmPC3proofeq6}, we deduce that for all $t \in [0,T)$ $$\E|\phi_{N,u}(T-t,\muu^N_t) - \phi_{N,u}(T-t,\mu_t)| \leq  C\, \E W^1(\muu^N_0,\mu_0) +C\frac{\ln(N)^{\frac{1}{\alpha}}}{N^{1-\frac{1}{\alpha}}}.$$ We conclude the proof of \eqref{PropPCOUeq1} by letting $t$ tend to $T$ thanks to the continuity of $\phi_{N,u}$. The second estimate \eqref{PropPCOUeq2} is also proved. Indeed, coming back to \eqref{ThmPC3proofeq1}, we see with the previous computations that $I_4$ an $I_5$ are centered because of the martingale property of compensated Poisson random integrals. Thus, we obtain that $$ |\E (\phi_{N,u}(T-t,\muu^N_{N,t}) - \phi_{N,u}(T-t,\mu_{N,t}))| \leq |\E (\phi_{N,u}(T,\muu^N_0) - \phi_{N,u}(T,\mu_0))| + \E|I_{3,B}|.$$ The initial data term is handled in the same way as previously and the control of  $\E|I_{3,B}|$ has already been done in \eqref{thmPC3prooff}. We thus obtain that for all $t \in [0,T)$ $$|\E(\phi_{N,u}(T-t,\muu^N_t) - \phi_{N,u}(T-t,\mu_t))| \leq  C\, \E W^1(\muu^N_0,\mu_0) + \frac{C}{N^{\alpha-1}},$$ and we conclude by letting $t$ tend to $T.$ Finally, to prove \eqref{PropPCOUeq2bis} in the case where the initial condition belongs to $\PPP_2(\R^d)$, we use the same argument as at the end of the proof of Theorem $3.6$ page $37$ of \cite{frikha2021backward}. To do this, we note that for a constant $C>0$ independent of $N$, we have for all $\mu \in \PPP_{\beta}(\R^d)$ and for all $v,v' \in \R^d$ $$ \left\vert\de\dellN (T,\mu)(v,v')\right\vert \leq C (1+|v|)(1+|v'|),$$ thanks to Proposition \ref{ThmPCremovjumps}. It yields $$ |\E (\phi_{N,u}(T,\muu^N_0) - \phi_{N,u}(T,\mu_0))| \leq \frac{C}{N},$$ which ends the proof. 

\end{proof}

\appendix

\section{Estimates on the density of a stable Ornstein-Uhlenbeck process}\label{appendixdensityestimates}

Fix $Z:=(Z_t)_t$ an $\alpha$-stable process on $\R^d,$ with $\alpha>1$ and having a non degenerate Lévy measure. More precisely, writing $y = r\theta \in \R^d$ with $r\in \R^+$ and $\theta \in \S^{d-1},$ where $\S^{d-1}$ denotes the unit sphere in $\R^d,$  the Lévy measure of $Z$ decomposes as $$ \nu(dy) = d\mu(\theta) \frac{dr}{r^{1+\alpha}},$$ for $\mu$ a non-zero finite measure on $\S^{d-1}.$ In Appendix \ref{appendixdensityestimates}, we denote by $\langle \cdot, \cdot \rangle$ the usual scalar product on $\R^d$. The Lévy symbol associated with $Z$ is given by $$ \psi (\lambda) =  \int_{\R^d} \left[ e^{i\langle \lambda, y \rangle} - 1 - i\langle \lambda, y \rangle   \right]\, d\nu(y).$$  We introduce, for $\delta \in (0,+\infty],$ the truncated symbol $ \psi_{\delta}$ defined for $\lambda \in \R^d$ by $$ \psi_{\delta} (\lambda) : = \int_{|y| < \delta}   \left[ e^{i\langle \lambda, y \rangle} - 1  -i\langle \lambda,y \rangle \right]\, d\nu(y).$$  We assume that $\nu $ satisfies the following non-degeneracy assumption \\

\noindent\textbf{(ND)} There exists $\eta >0$ such that for all $\lambda \in \R^d$ $$ \eta |\lambda|^{2} \leq \int_{\S^{d-1}} | \langle \lambda,\theta \rangle |^{2} \, d\mu(\theta).$$

Let us denote by $\NN$ the Poisson random measure associated with $Z$ and by $\NNN$ its compensated Poisson random measure. For any $t\in [0,T]$, we can write $$ Z_t = \int_0^t \int_{\R^d} z \, \NNN(ds,dz).$$ Let us fix $t\in (0,T].$ We define $\tilde{Z}^1 = (\tilde{Z}^1_s)_s$ by $$ \tilde{Z}^1_s := \int_0^s\int_{|z| \geq t^{\frac{1}{\alpha}} }z \, \NNN(du,dz),$$  and $\tilde{Z}^2_s:= Z_s - \tilde{Z}^1_s.$ Thus, one has $$ Z_t = \tilde{Z}^1_t + \tilde{Z}^2_t,$$ where $\tilde{Z}^1$ and $\tilde{Z}^2$ are independent Lévy processes with Lévy symbols respectively given, for all $\lambda \in \R^d$, by $$ \psi^1 (\lambda) = \psi_{\infty}(\lambda) - \psi_{t^{\frac{1}{\alpha}}}(\lambda) \quad\text{and} \quad \psi^2 (\lambda) =\psi_{t^{\frac{1}{\alpha}}}(\lambda). $$  

We also introduce, for $N \in \N$ such that $N \geq T^{\frac{1}{\alpha}},$ the truncated process defined for $s \in [0,T]$ by $$ Z_{N,s}:= Z_s -\int_0^s \int_{|z| \geq  N} z \, \NNN(du,dz).$$ 

We set for $s \in [0,T]$ 
$$\tilde{Z}^{3}_{N,s} := \int_0^s \int_{t^{\frac{1}{\alpha}} \leq|z| <N} z \, \NN(du,dz),$$ and $$ \tilde{Z}^{4}_{N,s}= - \int_0^s\int_{t^{\frac{1}{\alpha}} \leq |z| < N } z \, d\nu(z) \,du.$$ We have the following decomposition for $Z_{N,t}$ $$ Z_{N,t} = \tilde{Z}^2_t + \tilde{Z}^{3}_{N,t}  + \tilde{Z}^{4}_{N,t},$$ where the three random variables are mutually independent.\\

Let $A,B \in \MM_d(\R)$ be two matrices such that $B$ is invertible. The Ornstein-Uhlenbeck process $Y:=(Y_t)_t$ associated with $Z$ is defined as the solution to the following well-posed SDE $$ dY_t = AY_t dt + B\,dZ_t,$$ starting at $t=0$ from $0.$ Using Itô's formula, one can see that $Y$ can be expressed as \begin{equation}\label{stableintegral}
	 Y_t = \int_0^t e^{(t-s)A}B \, dZ_s. \end{equation} We also consider the Ornstein-Uhlenbeck process $(Y_{N,t})_t$ defined in the same way but with $Z$ replaced by the truncated stable process $Z_N.$ It writes, for $t \leq T$ \begin{align*}
Y_{N,t} &= \int_0^t e^{(t-s)A}B \,d\tilde{Z}^2_s + \int_0^t e^{(t-s)A}B \,d\tilde{Z}^{3}_{N,s} + \int_0^t e^{(t-s)A}B \,d\tilde{Z}^{4}_{N,s} \\ &=: \tilde{Y}^2_t + \tilde{Y}^{3}_{N,t} + \tilde{Y}^{4}_{N,t}. \end{align*}

We begin with the existence of a density for the small jumps term $\tilde{Y}^2$ of the truncated stable Ornstein-Uhlenbeck process $Y_N$.

\begin{Prop}\label{densityOU}
	For all $t\in(0,T]$ the distribution of $\tilde{Y}^2_t$ has a density with respect to the Lebesgue measure $\tilde{p}^2(t,\cdot)$ which belongs to the Schwartz space $\mathcal{S}(\R^d)$ of regular and rapidly decreasing functions. Moreover, for any $ m \geq 0,$ and for any multi-index $\beta \in \N^d,$ there exists a constant $C_{T,m,\beta}>0$ such that for all $t \in (0,T], \, x \in \R^d$, one has \begin{equation}\label{estimatesmalljump} \left\vert\partial^\beta_x \tilde{p}^2(t,x)\right\vert \leq \frac{C_{T,m,\beta} }{t^{\frac{d+|\beta|}{\alpha}}} \left(1+\frac{|x|}{t^{\frac{1}{\alpha}}}\right)^{-m}.\end{equation}
	
	Moreover, for all $\gamma \geq 0,$ there exists a constant $C>0$ depending only on $T,d,\alpha,\gamma,\beta,\eta$ such that \begin{equation}\label{momentsmalljump} \int_{\R^d} |x|^{\gamma} |\partial_x^{\beta}\tilde{p}^2(t,x)| \, dx \leq Ct^{\frac{\gamma - |\beta|}{\alpha}}. \end{equation}
	
\end{Prop}

\begin{proof} We fix $t >0.$ Reasoning as page $105$ of \cite{SatoLevyprocessesinfinitely1999}, the characteristic function of $\tilde{Y}^2_t$ is given by the function $$ p \in \R^d \mapsto \phi_t(p):=\exp\left(\int_0^t \psi^2(B^*e^{sA^*}p) \, ds\right).$$  Changing variables in $v:=\frac{s}{t}$ and $\rho := \frac{r}{t^{\frac1\alpha}}$, one has for $p\in \R^d$
	\begin{align*}
	\phi_t(p) &=  \exp \left( \int_0^t \int_0^{ t^{\frac1\alpha}} \int_{\S^{d-1}} e^{i\langle B^*e^{sA^*}p,r\theta \rangle} -1 - i\langle B^*e^{sA^*}p,r\theta \rangle  \, d\mu(\theta)\, \frac{dr}{r^{1+\alpha}}\, ds\right) \\ &=  \exp \left( \int_0^1 \int_0^1 \int_{\S^{d-1}} e^{i\langle B^*e^{tvA^*}p,t^{\frac1\alpha}\rho\theta \rangle} -1 - i\langle B^*e^{tvA^*}p,t^{\frac1\alpha}\rho\theta \rangle  \, d\mu(\theta)\, \frac{d\rho}{\rho^{1+\alpha}}\, dv\right) \\ &= \exp\left(\int_0^1 \psi_1(B^*e^{tvA^*}t^{\frac1\alpha}p) \, dv\right).
	\end{align*}
	
	Reasoning as in Lemma $3.2$ in \cite{Chenolderregularitygradient2020}, there exists $\eta >0$ such that for any $q \in \R^d$  \begin{equation}\label{estimatesmalljumps}
	\text{Re}(\psi_1(q)) \leq -\eta (|q|^{\alpha} \wedge |q|^2).
	\end{equation}
	We have thus for any $q \in \R^d$ \begin{align*}
	\left\vert \exp\left(\int_0^1 \psi_1(B^*e^{tvA^*}t^{\frac1\alpha}q) \, dv\right) \right\vert &=   \exp \left(\int_0^1 \text{Re}(\psi_1( B^*e^{tvA^*}t^{\frac1\alpha}q)) \,  dv\right) \\ &\leq \exp \left( -\eta\int_0^1 |B^*e^{tvA^*}t^{\frac1\alpha}q|^{\alpha} \wedge |B^*e^{tvA^*}t^{\frac1\alpha}q|^2 dv\right).
	\end{align*}

	We now prove that there exists $\eta_T >0$ such that for all $s\in [0,T]$ and for any $q \in \R^d$ \begin{equation}\label{coercivity}
	\vert B^*e^{sA^*}q \vert  \geq \eta_T |q|
	\end{equation}
	
	First, note that \begin{align*}
	\vert  B^*e^{sA^*}q \vert^2 & =\langle BB^* e^{sA^*} q,e^{sA^*}q\rangle \\ &\geq \eta_1 \vert e^{sA^*}q \vert^2,
	\end{align*} for $\eta_1>0$ since $BB^*$ is positive-definite because $B$ is invertible. We prove that there exists $\eta_2>0$ such that for any $q \in \R^d$ and $s \in [0,T]$  $$ \vert e^{sA^*}q \vert^2 \geq \eta_2 |q|^2,$$ which will conclude the proof of \eqref{coercivity}. Reasoning by contradiction, there exists $(s_n)_n \in [0,T]^{\N}$ and $(q_n)_n \in (\R^d)^{\N}$ with $ |q_n|=1$ for all $n$ and such that for all $n$ \begin{equation}\label{inegvp}
	\langle e^{s_nA^*}e^{s_nA}q_n,q_n \rangle \leq \frac1n.	\end{equation}  By compactness, one can assume that $(s_n)_n$ converges to $s \in [0,T]$ and that $(q_n)_n$ converges to $q\in \R^d$, with $|q|=1.$ Letting $n$ tend to infinity in \eqref{inegvp}, we obtain that $q = 0$ since $e^{sA}$ is invertible. This is a contradiction.\\
	
	Thus, there exists a constant $\eta_T>0$ such that for any $q \in \R^d$ and $t \in [0,T]$ \begin{equation}\label{integrabiltyfourier}
	\left\vert \exp\left(\int_0^1 \psi_1(B^*e^{tvA^*}q) \, dv\right) \right\vert \leq \exp \left( -\eta_T |q|^{\alpha} \wedge |q|^2 \right).
	\end{equation}
	This proves that the characteristic function of $\tilde{Y}^2_t$ is integrable. It follows that the distribution of $\tilde{Y}^2_t$ has a density with respect to the Lebesgue measure denoted by $\tilde{p}^2(t,\cdot)$ and given, for all $x \in \R^d$, by \begin{align*}\tilde{p}^2(t,x) & =\frac{1}{(2\pi)^d} \int_{\R^d} e^{-i\langle x, p \rangle}\exp\left(\int_0^1 \psi_1(B^*e^{tvA^*}t^{\frac1\alpha}p) \, dv\right) \, dp\end{align*}
	
	Changing variables in $q := t^{\frac1\alpha} p$ yields $$
	\tilde{p}^2(t,x) =\frac{t^{-\frac{d}{\alpha}}}{(2\pi)^d} \int_{\R^d} e^{-i\left\langle q,\frac{x}{t^{\frac1\alpha}} \right\rangle} \exp \left( \int_0^1 \psi_1( B^*e^{tvA^*}q) \,  dv\right) \, dq.$$
	
	Let us define for all $t >0$ and $q \in \R^d$ $$ g_t(q) :=  \exp \left(\int_0^1 \psi_1( B^*e^{tvA^*}q) \,  dv\right).$$

	We have proved in \eqref{integrabiltyfourier} that for any $q \in \R^d$ \begin{align}\label{decroissancefourier}
	|g_t(q)| \leq  e^{-\eta_T(|q|^\alpha\wedge |q|^2)}.
	\end{align}
	
	It follows that for any $t >0,$ $\tilde{p}^2(t,\cdot)$ belongs to $\CC^{\infty}(\R^d)$ and that for any multi-index $\beta \in \N^d$ \begin{equation}\label{fourierderivativedensity}
	\partial^\beta \tilde{p}^2(t,x)= \frac{t^{-\frac{d+|\beta|}{\alpha}}}{(2\pi)^d} \int_{\R^d} e^{-i\left\langle q,\frac{x}{t^{\frac1\alpha}} \right\rangle} (-iq)^\beta g_t(q) \, dq.\end{equation}
	
	Now, we prove that for all $t \in (0,T]$ $g_t \in \CC^{\infty}(\R^d).$  For the derivatives of order $1,$ we differentiate under the integral to obtain that for $ j \in \{1, \dots,d\}$ $$\partial_{q_j} g_t(q) =  \left(\int_0^1 \int_0^1 \int_{\S^{d-1}} i\rho \left[ e^{tvA}B\theta \right]_j \right(e^{i\langle B^*e^{tvA^*}q,\rho\theta \rangle} -1\left)\, d\mu(\xi)\, \frac{d\rho}{\rho^{1 + \alpha}} \, dv \right)g_t(q),$$  where $[v]_j$ denotes the $j-$th component of a vector $v \in \R^d.$ This term is well-defined since the continuity of $t \in [0,T] \mapsto e^{tA}$ ensures that there exists $C_T>0$ such that for all $q \in \R^d, \, v \in [0,1],\, \theta \in \S^{d-1}$ \begin{align}\label{ineqexponential}
	\left\vert e^{i\langle q, \rho e^{tvA}B \theta  \rangle} -1\right\vert \leq C_T |q| \rho.
	\end{align}We can easily prove by induction that for any $b \in \N^d,$ $\partial^b g_t(q)$ exists and can be expressed as a linear combination of terms of the form $$ \left(\prod_{l=1}^N  \int_0^1 \int_0^1 \int_{\S^{d-1}} \prod_{j=1}^{M_l} \left[\rho e^{tvA}B\theta\right]_{k_l}\left(e^{i\langle q, \rho e^{tvA}B \theta  \rangle} -\1_{M_l=1} \right)  \, d\mu(\theta)\, \frac{d\rho}{\rho^{1 + \alpha}} \, dv \right)g_t(q),$$ where $N \geq 1,$ $M_l \geq 1$ and $k_l \in \{ 1, \dots,d\}$. \\
	
	Using \eqref{decroissancefourier}, we deduce that there exists two  constants $C>0$ and $D\geq 1$ depending on $T$ and $b$ such that for all $t \in (0,T]$ \begin{align*}
	\left\vert\prod_{l=1}^N  \int_0^1 \int_0^1 \int_{\S^{d-1}} \prod_{j=1}^{M_l} \left[\rho e^{tvA}B\theta\right]_{k_l}\left(e^{i\langle q, \rho e^{tvA}B \theta  \rangle} -\1_{M_l=1} \right)  \, d\mu(\theta)\, \frac{d\rho}{\rho^{1 + \alpha}} \, dv \right\vert |g_t(q)| \leq  C (1+|q|^D)e^{-\eta_T(|q|^\alpha\wedge |q|^2)}.
	\end{align*}Indeed, the only difficulty appears when $M_i=1$ since $\int_0^1 \rho\, \frac{d\rho}{\rho^{1+\alpha}}$ is equal to infinity. However in this case, we can use \eqref{ineqexponential} and the fact that $$ \int_0^1 \rho^2 \, \frac{d\rho}{\rho^{1+ \alpha}} < + \infty.$$
	This shows that for all multi-index $a,b \in \N^d$ \begin{equation}\label{estimateSchwartz}
	\sup_{t \leq T} \Vert q^a\partial^b g_t(q) \Vert_{L^1_q(\R^d)} < + \infty. \end{equation}We easily see that the same controls hold replacing $g_t(q)$ by $(-iq)^{\beta} g_t(q).$ Thus, one has for $m \geq 0$ \begin{align*}
	\sup_{t\leq T} \sup_{x \in \R^d} (1+|x|)^m\left\vert \mathcal{F}_q ((-iq)^\beta g_t(q))(x) \right\vert &\leq C \sum_{\gamma \in \N^d, \, |\gamma| \leq m} \sup_{t\leq T} \left\Vert \partial^\gamma ((-iq)^{\beta}g_t(q)) \right\Vert_{L^1_q(\R^d)} \\ &< + \infty,
	\end{align*}
	where $\mathcal{F}_q$ denotes the Fourier transform with respect to the variable $q \in \R^d.$ 
	Using \eqref{fourierderivativedensity}, it follows that for a constant $C=C_{T,m,\beta}>0,$ we have for all $t \in (0,T]$ and for all $x \in \R^d$ $$ \left\vert \partial^\beta \tilde{p}^2(t,x) \right\vert \leq \frac{C}{t^{\frac{d+|\beta|}{\alpha}}} \left( 1 + \frac{|x|}{t^{\frac{1}{\alpha}}}\right)^{-m}.$$

	This ends the proof of \eqref{estimatesmalljump}. The moment estimate \eqref{momentsmalljump} directly follows from  \eqref{estimatesmalljump} choosing $m$ large enough and from a change of variables.

\end{proof}

We can now prove the same results for the Ornstein-Uhlenbeck process $Y_{N}$.

\begin{Prop}\label{densityestimateOUtruncated}
	For all $N \geq 1$ and for all $t \in(0,T]$ the distribution of $Y_{N,t}$ has a density with respect to the Lebesgue measure denoted by $p^N(t,\cdot) \in \CC^{\infty}(\R^d;\R^+).$ Moreover, for any $\beta \in \N^d$ and $\gamma \in [0,\alpha)$  there exists a constant $C>0$ depending only on $T,d,\alpha, \eta, \beta,\gamma$ such that for any $N \in \N$ \begin{equation}\label{densitymomentestimate}
		 \int_{\R^d} |x|^{\gamma} |\partial^{\beta}_x p^N(t,x)| \, dx \leq C t^{\frac{\gamma - |\beta|}{\alpha}}.\end{equation} The same results hold for $(Y_t)_t$ which has a density $p(t,\cdot)$ satisfying the moment estimate \eqref{densitymomentestimate}.\\
	 
	 Moreover, for any $N \geq 1$ and $t_0 \in (0,T]$, $p_N$ belongs to $\CC^{1,\infty}([t_0,T]\times \R^d)$ and $p^N(t,\cdot)$, $\partial_t p^N(t,\cdot)$ belong to $\mathcal{S}(\R^d)$ uniformly in time on $[t_0,T]$. More precisely, for all $ m \geq 0,$ and for any multi-index $\beta \in \N^d,$ there exists a constant $C_{t_0,T,m,\beta,N}>0$ such that for all $ x \in \R^d$, one has \begin{equation}\label{estimatesschwartzpN} \sup_{t\in [t_0,T]}\left\vert\partial^\beta_x p^N(t,x)\right\vert \leq C_{t_0,T,m,\beta,N} \left(1+|x|\right)^{-m},\end{equation} 
	 and 
	 \begin{equation}\label{estimatesschwartztimepN} \sup_{t\in [t_0,T]}\left\vert\partial^\beta_x \partial_t p^N(t,x)\right\vert \leq C_{t_0,T,m,\beta,N} \left(1+|x|\right)^{-m},\end{equation}.
\end{Prop}

\begin{proof}
	\textbf{Proof of \eqref{densitymomentestimate}.} Let us fix $t \in (0,T].$ We recall that we have decomposed $Y_{N,t}$ as a sum of three independent random variables \begin{align*}Y_{N,t} &= \tilde{Y}^2_t + \tilde{Y}^{3}_{N,t} + \tilde{Y}^{4}_{N,t} \\ &= \tilde{Y}^2_t+  \int_0^t \int_{t^{\frac{1}{\alpha}} \leq |z| < N} e^{(t-s)A}Bz \, \NN(ds,dz) - \int_0^t \int_{t^{\frac{1}{\alpha}} \leq |z|<N} e^{(t-s)A}B z \, d\nu(z) \, ds.\end{align*} 
	The density of $Y_{N,t}$ is thus given for all $x \in \R^d$ by $$ p^N(t,x) = \E( \tilde{p}^2(t,x- \tilde{Y}^{3}_{N,t} - \tilde{Y}^{4}_{N,t})).$$ The regularity of $p^N(t,\cdot)$ follows from Proposition \ref{densityOU} by differentiation under the integral. We now fix $\gamma \in [0,\alpha).$ Then, we write \begin{align*}
	\int_{\R^d} |x|^{\gamma} |\partial^{\beta}_x p^N(t,x)| \, dx &=  \int_{\R^d} |x|^{\gamma} |\E( \partial_x^{\beta}\tilde{p}^2(t,x- \tilde{Y}^{3}_{N,t} - \tilde{Y}^{4}_{N,t}))| \, dx \\ &\leq  C\int_{\R^d} \left(|x|^{\gamma} + \E| \tilde{Y}^{3}_{N,t}|^{\gamma} + \E |\tilde{Y}^{4}_{N,t}|^{\gamma}\right) | \partial_x^{\beta}\tilde{p}^2(t,x)| \, dx.
	\end{align*}
	Using the moment estimate \eqref{momentsmalljump}, it remains to prove that for a constant $C>0$ independent of $N$, one has for all $t \in [0,T]$ $$\E| \tilde{Y}^{3}_{N,t}|^{\gamma} + \E |\tilde{Y}^{4}_{N,t}|^{\gamma} \leq C t^{\frac{\gamma}{\alpha}}.$$ For $\tilde{Y}^{3}_{N,t},$ note that $\int_0^t\int_{t^{\frac{1}{\alpha}} \leq |z| } |z| \, \NN(ds,dz)$ has the same distribution as $t^{\frac1\alpha} \int_0^1\int_{1\leq |z| } |z| \, \NN(ds,dz).$ This can be seen easily using the characteristic functions. Moreover, since $\gamma < \alpha,$ we have $$ \E \left\vert \int_0^1\int_{1\leq |z| } |z| \, \NN(ds,dz) \right\vert^{\gamma} < + \infty.$$ We thus obtain that  \begin{align*}
	\E| \tilde{Y}^{3}_{N,t}|^{\gamma} & = \E \left\vert  \int_0^t \int_{t^{\frac{1}{\alpha}} \leq |z| < N} e^{(t-s)A}Bz \, \NN(ds,dz) \right\vert^{\gamma} \\ &\leq C_T \E \left\vert  \int_0^t \int_{t^{\frac{1}{\alpha}} \leq |z| } |z| \, \NN(ds,dz) \right\vert^{\gamma} \\ &\leq C_T t^{\frac{\gamma}{\alpha}}. 
	\end{align*}
	Finally, for $\tilde{Y}^{4}_{N,t},$ one has \begin{align*}
	\E|\tilde{Y}^{4}_{N,t}|^{\gamma} &= \E \left\vert\int_0^t \int_{t^{\frac{1}{\alpha}} \leq r <N}\int_{\S^{d-1}} e^{(t-s)A}B r\theta\, \frac{dr}{r^{1+\alpha}} \, d\mu(\theta)\, ds\right\vert^{\gamma}\\ &\leq \E \left\vert \frac{t^{1/\alpha -1}}{\alpha -1}\int_0^t \int_{\S^{d-1}} e^{(t-s)A}B \theta \, d\mu(\theta)\, ds\right\vert^{\gamma} \\&\leq C_T t^{\frac{\gamma}{\alpha}}.
	\end{align*}
	
	This concludes the proof of \eqref{densitymomentestimate}. The same reasoning holds for the moment estimate of $(Y_t)_t$.\\
	
	\textbf{Proof of \eqref{estimatesschwartzpN} and \eqref{estimatesschwartztimepN}.}
	As in the proof of Proposition \ref{densityOU}, we obtain the following expression for the density of $Y^N_t$ 
	
	\begin{equation*}
		p^N(t,x)  =\frac{1}{(2\pi)^d} \int_{\R^d} e^{-i\langle x, p \rangle}\exp\left(\int_0^t \psi_N(B^*e^{sA^*}p) \, ds\right) \, dp.
	\end{equation*} We follow the same lines as in the proof of Proposition \ref{densityOU}. Indeed, reasoning as in Lemma $3.2$ in \cite{Chenolderregularitygradient2020}, there exists $\eta >0$ such that for any $q \in \R^d$  \begin{equation}\label{estimatesPsiN}
	\text{Re}(\psi_N(q)) \leq -\eta (|q|^{\alpha} \wedge |q|^2),
\end{equation}

It follows that for some constant $\eta_T>0$, one has for all $ p\in \R^d$, $t \in [0,T]$. \[ \exp\left(\int_0^t \psi_N(B^*e^{sA^*}p) \, ds\right) \leq  \exp(-t\eta_T(|q|^\alpha \wedge |q|^2) ).\]
Moreover, we prove that there exists $C >0$ such that for any $q \in \R^d$ 
\begin{equation}\label{estimatesPsiN2}
	|\psi_N(q)| \leq C|q|^{\alpha},
\end{equation}
Indeed, using Taylor formula for the first integral, one has \begin{align*}
	|\psi_N(q)| & \leq  \int_{|y| \leq |q|^{-1}} |e^{i\langle q,y\rangle} - 1 -i \langle q,y\rangle| \, d\nu(y) + \int_{ |q|^{-1} \leq |y| \leq N} |e^{i\langle q,y\rangle} - 1 -i \langle q,y\rangle| \, d\nu(y)\\& \leq  \int_{|y| \leq |q|^{-1}}C|y|^2|q|^2 \, d\nu(y) + \int_{ |q|^{-1} \leq |y| \leq N} 2 + |q||y| \, d\nu(y) \\ &\leq C|q|^2 (|q|^{-1})^{2 - \alpha} + |q|^\alpha + C|q| (|q|^{-1})^{\alpha -1} \\ &\leq C|q|^\alpha.
\end{align*}

Thanks to \eqref{estimatesPsiN2}, we prove by differentiation under the integral that $p^N \in \CC^{1,\infty}((0,T] \times \R^d)$ and that \[\partial_t p^N(t,x) =\frac{1}{(2\pi)^d} \int_{\R^d} e^{-i\langle x, p \rangle} \psi_N(B^*e^{tA^*}p)\exp\left(\int_0^t \psi_N(B^*e^{sA^*}p) \, ds\right) \, dp. \]

As for \eqref{estimateSchwartz}, we prove that for all multi-index $a,b \in \N^d$, we have 

\begin{align*}&\sup_{t\in[t_0,T]}\left\Vert p^a\partial_p^b \left(\psi_N(B^*e^{tA^*}p)\exp\left(\int_0^t \psi_N(B^*e^{sA^*}p) \, ds\right)\right)\right\Vert_{L^1_p(\R^d)} \left\Vert p^a\partial_p^b \exp\left(\int_0^t \psi_N(B^*e^{sA^*}p) \, ds\right)\right\Vert_{L^1_p(\R^d)} \\ &=: C_{t_0,T,a,b,N}  < +\infty. \end{align*}

The same reasoning used in the proof of Proposition \ref{densityOU} allows to deduce \eqref{estimatesschwartzpN} and \eqref{estimatesschwartztimepN}. 
	
\end{proof}

\section{Proof of Proposition \ref{ThmPCremovjumps}}\label{proofpropositionreg}
\begin{proof}[Proof of Proposition \ref{ThmPCremovjumps}]
	
	We recall that $X_N$ denotes the solution to SDE \eqref{edsmckvtronquee1} driven by the truncated stable process $(Z_{N,t})_t.$ We fix $(t,\mu) \in (0,T] \times \PPP_{\beta}(\R^d)$ and $X_0$ a random variable with distribution $\mu$. Applying Itô's formula for the function $f : (t,x) \in [0,+\infty) \times \R^d \mapsto e^{-tA}x$, we obtain that for all $t \geq 0 $ $$ X_{N,t} = e^{tA} X_0 + \int_0^t e^{(t-s)A}A' \E X_{N,s}\, ds + \int_0^t e^{(t-s)A} B \, dZ_{N,s}.$$ By differentiating the map $t \mapsto e^{-tA}\E X_{N,t},$ we deduce that for all $t \geq 0 $ $$ \E X_{N,t} = e^{t(A+A')}\E X_0.$$ It yields for all $t \geq 0$ \begin{align}\label{eq:dec_OU}
	X_{N,t} \notag&= e^{tA} X_0 + \left(\int_0^t e^{(t-s)A}A' e^{s(A+A')}\, ds \right) \E X_0 + \int_0^t e^{(t-s)A} B \, dZ_{N,s} \\ &=:  e^{tA} X_0 + K_t\E X_0 + Y_{N,t}.\end{align} The process $(Y_{N,t})_t$ has been studied in Appendix \ref{appendixdensityestimates}. By Proposition \ref{densityestimateOUtruncated}, for $t>0,$ it admits a density $p^N(t,\cdot)$ satisfying the moment estimates \eqref{densitymomentestimate} and the Schwartz estimates \eqref{estimatesschwartzpN} and \eqref{estimatesschwartztimepN}. Since the moment estimates \eqref{densitymomentestimate} are uniform in $N,$ we denote $X_N$ by $X$, for the sake of clarity, and we omit all the subscripts $N$. We deduce that $X_t$ as a density with respect to the Lebesgue measure depending on the initial data $\mu \in \PP$ and denoted by $q(\mu,t,\cdot).$ Moreover, using  \eqref{eq:dec_OU}, we have for all $y \in \R^d$ \begin{equation}\label{fundsolution}
	q(\mu,t,y) = \int_{\R^d} p(t,y-e^{tA}x - K_t M(\mu)) \, d\mu(x),
	\end{equation} 
	where $ p=p^N$ is given by Proposition \ref{densityestimateOUtruncated}, and $M(\mu) := \int_{\R^d} x \, d\mu(x),$ for $\mu \in \PPP_1(\R^d).$\\

We deduce using Proposition \ref{densityestimateOUtruncated} that $q(\mu,\cdot,y)$ is of class $\CC^1$ on $(0,T]$ and that $q(\cdot,t,y)$ admits a linear derivative given for all $v \in \R^d$ by \begin{align*}
	\frac{\delta}{\delta m} q(\mu,t,y)(v) = p(t,y-e^{tA}v - K_tM(\mu)) - \int_{\R^d} \partial_y p(t,y-e^{tA}x - K_t M(\mu))\cdot (K_t v) \, d\mu(x).
	\end{align*}It is easy to see by differentiation under the integral that for all $v,y \in \R^d$, one has  \begin{align*}\partial_v \frac{\delta}{\delta m} q(\mu,t,y)(v) =  -e^{tA^*}\partial_y p(t,y-e^{tA}v - K_tM(\mu)) - \int_{\R^d} K_t^*\partial_y p(t,y-e^{tA}x - K_t M(\mu)) \, d\mu(x).  \end{align*} 

Using Proposition \ref{densityestimateOUtruncated}, we see that $\partial_t q$ is continuous on $\PP \times (0,T] \times \R^d$ and that $\frac{\delta}{\delta m} q$ and $\partial_v \frac{\delta}{\delta m} q$ are continuous on $\PP \times (0,T] \times \R^d \times \R^d.$\\

	Let us fix $u\in \mathcal{C}.$ Since SDE \eqref{MKVstableOUdef} is time-homogeneous, the associated function $\phi_u$ is thus given, for any $\mu \in \PP$ and for any $t \in (0,T]$ by $ \phi_u(t,\mu) = u(\theta(\mu,t)),$ where $\theta (\mu,t):= q(\mu,t,y)\,dy$ is the distribution of any solution to SDE \eqref{MKVstableOUdef} starting as $t=0$ from the distribution $ \mu \in \PP.$ 

Moreover, reasoning exactly as in the proof on Proposition 2.3 of \cite{frikha2021backward}, we prove using Proposition \ref{densityestimateOUtruncated} that $\phi_u$ satisfies the regularity results $(1)$ and $(2)$ of Proposition \ref{ThmPCremovjumps}. For the time derivative, one has $$ \partial_t \phi_u (\mu,t) = \int_{\R^d} \del (\theta(\mu,t))(y) \partial_t q(\mu,t,y) \,dy.$$ Moreover, it ensures that for all $v \in \R^d$  $$ \dell (\mu,t)(v) = \int_{\R^d}  \del (\theta(\mu,t))(y) \frac{\delta q}{\delta m} (\mu,t,y)(v) \, dy,$$  
	
	and that \begin{align}\label{expression_lions_derivative} \partial_v\dell (t,\mu)(v) &= -\int_{\R^d}  \del (\theta(\mu,t))(y) e^{tA^*} \partial_y p(t,y-e^{tA}v - K_tM(\mu))\, dy \\ \notag&\quad - \int_{\R^d}\int_{\R^d}  \del (\theta(\mu,t))(y) K_t^*\partial_y p(t,y-e^{tA}x - K_t M(\mu)) \, d\mu(x) \, dy, 
	\end{align}
	where $A^*$ denotes the transpose matrix of a matrix $A$. Using Proposition \ref{densityestimateOUtruncated}, we obtain that for all $x \in \R^d$  \begin{equation}\label{integralderivative}
	\int_{\R^d} \partial_y p(t,y-e^{tA}x - K_tM(\mu)) \, dy =0.\end{equation} We can thus write \eqref{expression_lions_derivative} in the following centered form \begin{align*} &\partial_v\dell (t,\mu)(v)\\ &= -\int_{\R^d}  \left(\del (\theta(\mu,t))(y) - \del (\theta(\mu,t))(e^{tA}v + K_t M(\mu)) \right) e^{tA^*} \partial_y p(t,y-e^{tA}v - K_tM(\mu))\, dy \\ &\quad- \int_{\R^d}\int_{\R^d}  \left(\del (\theta(\mu,t))(y) - \del (\theta(\mu,t))(e^{tA}x + K_t M(\mu)) \right) K_t^*\partial_y p(t,y-e^{tA}x - K_t M(\mu)) \, d\mu(x) \, dy. \end{align*}
 Using the Lipschitz assumption on $\del$ as well as the moment estimate in Proposition \ref{densityestimateOUtruncated}, we deduce that there exists $C>0$ depending only on $T,$ which may change from line to line, such that for any  $v \in \R^d$ \begin{align}\label{exempleOUeq1}
	\left\vert \partial_v\dell (t,\mu)(v)\right\vert\notag &\leq C \int_{\R^d} |y||\partial_ yp(t,y)| \, dy \\ &\leq C.
	\end{align}
	 Similarly, one has  for all $v \in \R^d$  \begin{align*} &\partial^2_v\dell (t,\mu)(v)\\ &= -\int_{\R^d}  \left(\del (\theta(\mu,t))(y) - \del (\theta(\mu,t))(e^{tA}v + K_t M(\mu)) \right) e^{tA^*} \partial_y^2 p(t,y-e^{tA}v- K_tM(\mu))e^{tA}\, dy.  \end{align*} 
	 Using the Lipschitz assumption on $\del$ as well as the moment estimate in Proposition \ref{densityestimateOUtruncated}, we deduce that there exists $C>0$ depending only on $T,$ such that for any $v \in \R^d$ \begin{align}\label{exempleOUeq2}
	\left\vert \partial_v^2\dell (t,\mu)(v)\right\vert\notag &\leq C \int_{\R^d} |y| |\partial^2_ yp(t,y)| \, dy \\ &\leq Ct^{-\frac{1}{\alpha}}.
	\end{align}
	Following again the same lines as in Proposition 2.3 of \cite{frikha2021backward}, we obtain that for all $v\in \R^d$, the function $ \partial_v\dell (t,\cdot)(v) $ admits a linear derivative given for all $v' \in \R^d$ by
	
	\begin{align*} &\de \partial_v\dell (t,\mu)(v,v')\\ &= -\int_{\R^{2d}}  \de\del (\theta(\mu,t))(y,y')  e^{tA^*} \partial_y p(t,y-e^{tA}v - K_tM(\mu)) p(t,y'-e^{tA}v'-K_tM(\mu))\, dy\, dy' \\ &\quad + \int_{\R^{3d}} \de\del (\theta(\mu,t))(y,y')e^{tA^*}\partial_y p(t,y-e^{tA}v - K_t M(\mu)) \\ & \hspace{8cm} \partial_y p(t,y' - e^{tA}x - K_tM(\mu))\cdot (K_t v')  \, dy \, dy' \, d\mu(x) \\ &\quad + \int_{\R^{d}}  \del (\theta(\mu,t))(y)  e^{tA^*} \partial^2_y p(t,y-e^{tA}v - K_tM(\mu)) K_tv'\, dy \\ &\quad- \int_{\R^{d}}  \del (\theta(\mu,t))(y) K_t^* \partial_y p(t,y-e^{tA}v' - K_tM(\mu)) \, dy \\ &\quad-\int_{\R^{3d}}  \de\del (\theta(\mu,t))(y,y')  K_t^* \partial_y p(t,y-e^{tA}x - K_tM(\mu)) p(t,y'-e^{tA}v'-K_tM(\mu))\, dy\, dy'\, d\mu(x)\\ &\quad + \int_{\R^{4d}}  \de\del (\theta(\mu,t))(y,y')  K_t^* \partial_y p(t,y-e^{tA}x - K_tM(\mu)) \\ & \hspace{8cm} \left(\partial_y p(t,y'-e^{tA}x'-K_tM(\mu))\cdot K_tv' \right)\, dy\, dy'\, d\mu(x')\, d\mu(x)  \\ &\quad +  \int_{\R^{2d}}  \del (\theta(\mu,t))(y) K_t^* \partial^2_y p(t,y-e^{tA}x - K_tM(\mu))K_tv' \, dy \, d\mu(x). 
	\end{align*}
	By differentiation under the integral, we easily deduce that for all $v,v' \in \R^d$
	
	\begin{align*} &\partial_{v'}\de \partial_v\dell (t,\mu)(v,v')\\ &= \int_{\R^{2d}}  \de\del (\theta(\mu,t))(y,y')  e^{tA^*} \partial_y p(t,y-e^{tA}v - K_tM(\mu))\otimes e^{tA^*}\partial_y p(t,y'-e^{tA}v'-K_tM(\mu))\, dy\, dy' \\ &\quad + \int_{\R^{3d}} \de\del (\theta(\mu,t))(y,y')e^{tA^*}\partial_y p(t,y-e^{tA}v - K_t M(\mu)) \otimes K_t^*\partial_y p(t,y' - e^{tA}x - K_tM(\mu)) \, dy \, dy' \, d\mu(x) \\ &\quad + \int_{\R^{d}}  \del (\theta(\mu,t))(y)  e^{tA^*} \partial^2_y p(t,y-e^{tA}v - K_tM(\mu)) K_t\, dy \\ &\quad + \int_{\R^{d}}  \del (\theta(\mu,t))(y) K_t^* \partial^2_y p(t,y-e^{tA}v' - K_tM(\mu))e^{tA} \, dy \\ &\quad +\int_{\R^{3d}}  \de\del (\theta(\mu,t))(y,y')  K_t^* \partial_y p(t,y-e^{tA}x - K_tM(\mu)) \otimes e^{tA^*} \partial_y p(t,y'-e^{tA}v'-K_tM(\mu))\, dy\, dy'\, d\mu(x)\\ &\quad + \int_{\R^{4d}}  \de\del (\theta(\mu,t))(y,y')  K_t^* \partial_y p(t,y-e^{tA}x - K_tM(\mu)) \\ & \hspace{8cm} \otimes \left(K_t^*\partial_y p(t,y'-e^{tA}x'-K_tM(\mu)) \right)\, dy\, dy'\, d\mu(x)\, d\mu(x')  \\ &\quad +  \int_{\R^{2d}}  \del (\theta(\mu,t))(y) K_t^* \partial^2_y p(t,y-e^{tA}x - K_tM(\mu))K_t \, dy \, d\mu(x). 
	\end{align*}
	
	It follows from \eqref{integralderivative} that we can write 
	
	\begin{align*} &\partial_{v'}\de \partial_v\dell (t,\mu)(v,v')\\ &= \int_{\R^{2d}}  \left(\de\del (\theta(\mu,t))(y,y')  - \de\del (\theta(\mu,t))(e^{tA}v + K_t M(\mu),e^{tA}v' + K_t M(\mu)) \right) \\ & \hspace{4cm}e^{tA^*} \partial_y p(t,y-e^{tA}v - K_tM(\mu))\otimes e^{tA^*}\partial_y p(t,y'-e^{tA}v'-K_tM(\mu))\, dy\, dy' \\ &\quad + \int_{\R^{3d}}  \left(\de\del (\theta(\mu,t))(y,y')  - \de\del (\theta(\mu,t))(e^{tA}v + K_t M(\mu),e^{tA}x + K_t M(\mu)) \right) \\ & \hspace{4cm}e^{tA^*}\partial_y p(t,y-e^{tA}v - K_t M(\mu)) \otimes K_t^*\partial_y p(t,y' - e^{tA}x - K_tM(\mu)) \, dy \, dy' \, d\mu(x) \\ &\quad + \int_{\R^{d}} \left( \del (\theta(\mu,t))(y) - \del (\theta(\mu,t))(e^{tA}v + K_tM(\mu)) \right)  e^{tA^*} \partial^2_y p(t,y-e^{tA}v - K_tM(\mu)) K_t\, dy \\ &\quad + \int_{\R^{d}}  \left( \del (\theta(\mu,t))(y) - \del (\theta(\mu,t))(e^{tA}v' + K_tM(\mu)) \right)  K_t^* \partial^2_y p(t,y-e^{tA}v' - K_tM(\mu))e^{tA} \, dy \\ &\quad +\int_{\R^{3d}}  \left(\de\del (\theta(\mu,t))(y,y')  - \de\del (\theta(\mu,t))(e^{tA}x + K_t M(\mu),e^{tA}v' + K_t M(\mu)) \right) \\ & \hspace{4cm} K_t^* \partial_y p(t,y-e^{tA}x - K_tM(\mu)) \otimes e^{tA^*} \partial_y p(t,y'-e^{tA}v'-K_tM(\mu))\, dy\, dy'\, d\mu(x)\\ &\quad + \int_{\R^{4d}}  \left(\de\del (\theta(\mu,t))(y,y')  - \de\del (\theta(\mu,t))(e^{tA}x + K_t M(\mu),e^{tA}x' + K_t M(\mu)) \right) \\ & \hspace{2cm}  K_t^* \partial_y p(t,y-e^{tA}x - K_tM(\mu))  \otimes \left(K_t^*\partial_y p(t,y'-e^{tA}x'-K_tM(\mu)) \right)\, dy\, dy'\, d\mu(x)\, d\mu(x')  \\ &\quad +  \int_{\R^{2d}}  \left( \del (\theta(\mu,t))(y) - \del (\theta(\mu,t))(e^{tA}x + K_tM(\mu)) \right) K_t^* \partial^2_y p(t,y-e^{tA}x - K_tM(\mu))K_t \, dy \, d\mu(x). 
	\end{align*}
	
	Using the Lipschitz assumption on $\del$ and $\de \del$ as well as the moment estimate in Proposition \ref{densityestimateOUtruncated}, we deduce that there exists $C>0$ depending only on $T$ such that for any $v,v' \in \R^d$  \begin{align}\label{exempleOUeq3}
	\left\vert \partial_{v'}\de \partial_v\dell (t,\mu)(v,v')\right\vert\notag &\leq C \left(\int_{\R^d} |y| |\partial^2_yp(t,y)| \, dy  + \int_{\R^d} |y| |\partial_yp(t,y)| \, dy \int_{\R^d}|\partial_y p(t,y)| \, dy \right)\\ &\leq Ct^{-\frac{1}{\alpha}}.
	\end{align}Note that the right-hand side term belongs to $L^1(0,T)$ since  $\alpha \in (1,2)$. Thus, \eqref{exempleOUeq1}, \eqref{exempleOUeq2} and \eqref{exempleOUeq3} ensure that the point $(3)$ is Proposition \ref{ThmPCremovjumps} is satisfied. It concludes the proof.

\end{proof}

\subsection*{Acknowledgements}

I would like to thank Paul-Eric Chaudru de Raynal for his supervision, advices and for his careful reading of the paper.

	\bibliographystyle{alpha}
\bibliography{bibli}

\end{document}